\documentclass[10pt,a4paper]{article}
\usepackage{graphicx}
\usepackage{dcolumn}
\usepackage{bm}
\usepackage{graphics}
\usepackage[centerlast]{caption2}
\usepackage{epsfig}
\usepackage{color}
\usepackage{amssymb}
\usepackage{amsthm}
\usepackage{lineno}
\usepackage{mathrsfs}
\usepackage{amsmath}
\usepackage{amsfonts}
\usepackage{bm,cite}
\numberwithin{equation}{section}
\usepackage{tabularx}
\usepackage[top=2cm, bottom=2cm, left=3cm, right=3cm] {geometry}
\usepackage{hyperref}
\hypersetup{hypertex=true, pagebackref=true, colorlinks=true, linkcolor=blue, anchorcolor=blue, citecolor=blue}
\newtheorem{thm}{Theorem}[section]
\newtheorem{cor}[thm]{Corollary}
\newtheorem{lem}[thm]{Lemma}
\newtheorem{prop}[thm]{Proposition}
\theoremstyle{definition}

\theoremstyle{definition}
\newtheorem{rem}[thm]{Remark}
\theoremstyle{definition}

\def\R{{\mathbb{R}}}
\usepackage{lineno}

\begin{document}
	\title{Twice Epi-Differentiability of Orthogonally Invariant Matrix Functions and Application}
	\author{Jiahuan He\\School of Science,
		Harbin University of Science and Technology,\\
		Harbin, 150080, P.R. China\\
		Chao Kan and Wen Song\footnote{Corresponding author, E-mail address: wsong@hrbnu.edu.cn}\\School of Mathematical and Sciences, Harbin Normal University,\\ Harbin, 150025, P.R. China}
	\date{}
    \maketitle
	\par\noindent
    \textrm{\textbf{Abstract.}}
	In this paper, our focus lies on the study of the second-order variational analysis of orthogonally invariant matrix functions. It is well-known that an orthogonally invariant matrix function is an extended-real-value function defined on ${\mathbb M}_{m,n}\,(n \leqslant m)$ of the form $f \circ \sigma$ for an absolutely symmetric function $f \colon \R^n \rightarrow [-\infty,+\infty]$ and the singular values $\sigma \colon {\mathbb M}_{m,n} \rightarrow \R^{n}$. We establish several second-order properties of orthogonally invariant matrix functions, such as parabolic epi-differentiability, parabolic regularity, and twice epi-differentiability when their associated absolutely symmetric functions enjoy some properties. 
	Specifically, we show that the nuclear norm of a real $m \times n$ matrix is twice epi-differentiable and we derive an explicit expression of its second-order epi-derivative. 
	Moreover, for a convex orthogonally invariant matrix function, we calculate its second subderivative and present sufficient conditions for twice epi-differentiability. This enables us to establish second-order optimality conditions for a class of matrix optimization problems.
	\vspace*{2mm}
	\vspace*{2mm}
	
	\renewcommand{\thefootnote}{\fnsymbol{footnote}}
	
	\section{Introduction}
	 Consider convex matrix optimization problem
	 $$
	 \min_{X \in {\mathbb M}_{m,n}} \psi(X) + \theta(X),
	 $$
	 where ${\mathbb M}_{m,n}$ stands for the Euclidean space of $m \times n$ real matrices, with inner product $\langle X,Y \rangle = {\rm tr}(X^T Y)$, $\psi \colon {\mathbb M}_{m,n} \rightarrow \R$ is a twice continuously differentiable convex function, and $\theta \colon {\mathbb M}_{m,n} \rightarrow [-\infty,+\infty]$ is a proper closed convex function. The problems described above represent a broad and significant class of convex matrix optimization problems with extensive applications in many fields, including matrix completion, rank minimization, graph theory, and machine learning. The examples of $\theta$ include the indicator function over the positive semidefinite cone \cite{VB}, the nuclear norm (defined as the sum of all singular values of a matrix) \cite{Watson  1992}, the spectral norm (which represents the largest singular value of a matrix) \cite{Mathias 1990}, and the matrix Ky Fan $k$-norm ($1 \leqslant k \leqslant n$) function (which corresponds to the sum of $k$ largest singular values of a matrix)\cite{Ding 2017}. The functions $\theta$ in the aforementioned examples belong to a distinguished class, which can be expressed either in the form of
	 \begin{equation*}
	 	\theta(X)=(g \circ \lambda)(X),\,\, X \in {\mathbb S}^n
	 \end{equation*}
	 with the function $g \colon \R^n \rightarrow [-\infty,+\infty]$ being proper closed convex and symmetric, and $\lambda$ is a function, which assigns to each matrix $ X\in {\mathbb S}^n$ its eigenvalue vector $(\lambda_1(X), \ldots, \lambda_n(X))$ arranged in nonincreasing order, referred to as spectral functions \cite{Lewis 1999}, or in the form of
	 \begin{equation}\label{composite form of orthogonally invariant matrix function}
	 	\theta(X)=(f \circ \sigma)(X),\,\, X\in {\mathbb M}_{m,n}
	 \end{equation}
	 with the function $f \colon \R^n \rightarrow [-\infty,+\infty]$ being proper closed convex and absolutely symmetric, and $\sigma$ is a function, which assigns to each matrix $ X\in {\mathbb M}_{m,n}$ its singular value vector $(\sigma_1(X), \ldots, \sigma_n(X))$ arranged in nonincreasing order, known as orthogonally invariant matrix functions \cite{Lewis 1995} or singular value functions \cite{LSI}. Recall that $F \colon {\mathbb M}_{m,n} \rightarrow [-\infty,+\infty] $ is an orthogonally invariant matrix function if for any matrix $X$ in ${\mathbb M}_{m,n}$, any $m \times m$ orthogonal matrix $U_m$ and any $n \times n$ orthogonal matrix $U_n$, one has
	 $$
	 F(U_m^T X U_n)=F(X).
	 $$
	 The work of Lewis and Sendov \cite{LSI, LSII} serves as a foundational contribution to the variational analysis of this class of functions. Lewis \cite{Lewis 1995} showed that convexity, lower semi-continuity, differentiability, and essentially smoothness of the absolutely symmetric function $f$ in \eqref{composite form of orthogonally invariant matrix function} are inherited by the orthogonally invariant matrix function $\theta$. 
	 A similar observation was made by Lewis and Sendov \cite{LSI} about Fréchet differentiability, regularity, and strictly differentiability, and by 
	 Cui et al. \cite{CDZ} about the ${\mathcal C}^2$-cone reducibility and the metric subregularity of their subdifferentials. The calculation of various notions of subdifferentials for orthogonally invariant matrix functions, including the limiting subdifferential, the Clarke subdifferential, and the proximal subdifferential, which play an important role in second-order variational analysis, was demonstrated in Lewis and Sendov \cite{LSI, LSII}. 
	 The central question addressed in this paper is whether the remarkable pattern observed can be extended to other significant second-order variational properties. 
	 
	 It was proved in Mohammadi and Sarabi \cite{MS} that important second-order variational properties of a composite function $g \circ \varphi$, where $g \colon {\mathbb Y} \rightarrow [-\infty,+\infty]$ is convex and $\varphi \colon {\mathbb X} \rightarrow {\mathbb Y}$ is twice differentiable with ${\mathbb X}$ and ${\mathbb Y}$ being finite-dimensional Hilbert space, can be established at any $\bar x \in {\rm dom}(g \circ \varphi)$, provided that $g$ is parabolically regular and that the following metric subregularity constraint qualification is satisfied: There exists a constant $\kappa \geq 0$ such that the estimate
	 $$
	 {\rm dist}(x, {\rm dom}(g \circ \varphi)) \leq \kappa {\rm dist}(\varphi(x), {\rm dom}g),
	 $$
	 holds for all $x$ in a neighborhood of $\bar x$. It was observed in \cite{DLMS} that the metric subregularity constraint qualification, when applied to spectral functions, is inherently satisfied. In addition, because a quadratic expansion of parabolic nature suffices to substitute for the twice differentiability of the inner function, Mohammadi and Sarabi \cite{MS2024} characterized parabolic regularity for spectral functions and calculated their second subderivative when the symmetric functions associated with them are convex. In the spirit of  Mohammadi and Sarabi \cite{MS,MS2024}, we study second-order variational properties, including parabolic epi-differentiability and twice epi-differentiability, of orthogonally invariant matrix functions. 
	 
	 The first-order directional derivatives of the singular values, as developed in \cite{DST}, serve as a powerful tool for characterizing the first-order necessary optimality conditions in matrix cone optimization problems. 
	 To derive the second-order necessary or sufficient conditions for the matrix optimization problem, it is essential to examine the second-order directional derivatives of the singular values. 
	 Zhang et al. \cite{ZZX} proposed a direct method to derive the formula for the second-order directional derivative of any eigenvalue of a symmetric matrix, as presented in Torki \cite{Torki 2001}. From this, they established a corresponding formula for the second-order directional derivative of any singular value of a matrix. 
	 
	 The structure of the paper is organized as follows. In Section 2, we review key concepts from variational analysis, along with important notions related to eigenvalues and singular values, which are central to the discussions throughout the paper. Section 3 presents explicit expressions for the first- and second-order epi-derivatives of the nuclear norm. In Section 4, we establish a chain rule for the subderivatives of orthogonally invariant matrix functions in \eqref{composite form of orthogonally invariant matrix function}. Additionally, we derive the tangent cone and second-order tangent sets associated with orthogonally invariant matrix sets. It is further demonstrated that the subderivative is a symmetric function with respect to a subset of the space of signed permutation matrices. Finally, in Section 5, we compute the second subderivative of orthogonally invariant matrix functions under the assumption that the corresponding absolutely symmetric functions are convex. We also present sufficient conditions for the twice epi-differentiability of orthogonally invariant matrix functions. As an application, we derive second-order optimality conditions for a class of matrix optimization problems.

    \section{Preliminaries}
	\subsection{Notation}
	 Throughout the whole paper, we will assume that $m$ and $n$ are natural numbers and $n \leqslant m$. 
	 The notations and concepts of convex analysis that we employ are standard \cite{H-UL, Rockafellar 1970, Rockafellar Wets}. Given a nonempty set $C \subset \mathbb{R}^n$, the indicator function of the set $C$, denoted by $\delta_C$, is defined by $\delta_C(x) = 0$ if $x \in C$ and $+ \infty$ otherwise. 
	 The support function of the set $C$ is defined by $\sigma_C(\cdot) := \sup_{s\in C}\langle s, \cdot \rangle$.  
	 The distance function from a point $x$ to the set $C$, denoted by ${\rm dist}(x,C)$ or ${\rm d}_C(x)$, is defined by ${\rm dist}(x,C) = {\rm d}_C(x):=\inf\{\|x-y\| \colon y \in C\}$. The tangent cone to $C$ at $\bar x$ is defined by
	 $$
	 T_C(\bar x) = \{ w \in \mathbb{R}^n \colon \exists t_k \downarrow 0, w_k \to w \,\,\mbox{as}\,\, k \to \infty \,\,\mbox{with}\,\, \bar x +t_kw_k \in C\}.
	 $$
	 The second-order tangent set to $C \subset \mathbb{R}^n$ at $\bar x \in C$ for a tangent vector $w \in T_C(\bar x)$ is given by
	 $$
	 T_C^2(\bar x,w) = \{u \in \mathbb{R}^n \colon \exists t_k \downarrow 0, u_k \to u \,\,\mbox{as}\,\, k \to \infty \,\,\mbox{with}\,\, \bar x +t_k w + \frac{1}{2}t_k^2 u_k \in C \}.
	 $$
	 A set $C \subset \mathbb{R}^n$ is called parabolically derivable at $\bar x $ for $w$ if $T_C^2(\bar x,w)$ is nonempty, and for each $u \in T_C^2(\bar x,w)$, there exists $\epsilon >0$ and an arc $\xi \colon [0,\epsilon] \to C$ with $\xi(0)=\bar x$, $\xi^{\prime}_+(0)=w$, and $\xi^{\prime \prime}_+(0)=u$, where $\xi^{\prime \prime}_+(0) = \lim_{t \downarrow 0} [\xi(t)-\xi(0)-t\xi^{\prime}_+(0)]/ \frac{1}{2}t^2$.
	 
	 Let $g \colon \mathbb{R}^n \to \mathbb{R} \cup \{\pm \infty\}$ be a function. The domain and epigraph of $g$ are defined as
	 ${\rm dom}g:= \{ x \in \mathbb{R}^n \colon g(x) < \infty\}$ and ${\rm epi}g := \{ (x,r) \in \mathbb{R}^n \times \mathbb{R} \colon g(x) \le r\}.$ $g$ is proper if ${\rm dom}g \ne \emptyset$  and $g(x) > -\infty $ for all $x \in {\rm dom}g$. $g$ is convex and lower semi-continuous (lsc) if ${\rm epi}g$ is convex and closed, respectively. 
	 The function $g$ is called locally Lipschitz continuous around $\bar x$ relative to $C \subset {\rm dom}g $ with constant $l \geq 0$ if $\bar x \in C$ with $g(\bar x)$ finite, and there exists a neighborhood $U$ of $\bar x$ such that 
	 \begin{equation*}
	 	|g(x)-g(y)| \leq l\|x-y\| \,\, \mbox{for all} \,\, x,y \in U \cap C.
	 \end{equation*}
	 We say that  $g$ is locally Lipschitz continuous relative to $C$ if it is locally Lipschitz continuous around every $\bar x \in C$ relative to $C$. 
	 The directional derivative of $g$ at $\bar x$ in the direction $w$ is defined as
	 \begin{equation*}
	 	g^{\prime}(\bar x ; w) := \lim_{t \downarrow 0} \frac{g(\bar x + t w) -g(\bar x)}{t},
	 \end{equation*}
	 and the subderivative of $g$ at $\bar x$ in the direction $w$ is defined as
	 \begin{equation*}
	 	{\rm d} g(\bar x)(w) := \liminf_{ {t \downarrow 0} \atop{ w^{\prime} \to w}} \frac{g(\bar x + t w^{\prime})-g(\bar x)}{t}.
	 \end{equation*}
     It is known that if $g$ is locally Lipschitz continuous around $\bar x$ relative to its domain, then ${\rm dom}\,{\rm d}g(\bar x) = T_{{\rm dom}g} (\bar x)$(cf. Mohammadi and Sarabi \cite[Proposition 2.2]{MS}).

     In what follows, we'll review the classical notion of second-order variational analysis. Let's first recall the concepts of the second subderivative and parabolic regularity for functions, respectively.
     Given a function $g \colon \mathbb{R}^n \to \mathbb{R} \cup \{\pm \infty\}$ and a point $\bar x \in \R^n$ with $g(\bar x)$ finite, we say $g$ is (parabolic) second order directionally differentiable at $\bar x$ if $g$ is directionally differentiable at $\bar x$ and for any $w,z \in \mathbb{R}^n$,
     $$
     \lim_{\tau \downarrow 0} \frac{g(\bar x+ \tau w + \frac{1}{2} \tau^2 z)-g(\bar x)-\tau g^{\prime}(\bar x;w)}{\frac{1}{2} \tau^2} \,\, \mbox{exists};
     $$
     and the above limit is said to be the (parabolic) second order directional derivative of $g$ at $\bar x$ along the directions $w$ and $z$, denoted by $g^{\prime \prime}(x;w,z)$.
  
     Define the parametric family of second-order difference quotients for $g$ at $\bar x$ for $v \in \R^n$ by
     $$
     \bigtriangleup_{\tau}^2 g(\bar x \mid v)(w) := \frac{g(\bar x+\tau w)-g(\bar x)-\tau \langle v,w\rangle}{\frac{1}{2} \tau^2} \,\,\mbox{with}\,\, w \in \R^n, \tau>0.
     $$
     The second subderivative of $g$ at $\bar x$ for $v \in \R^n$  is defined by
     $$
     {\rm d}^2 g(\bar x \mid v)(w) = \liminf_{ {\tau \downarrow 0} \atop{ w^{\prime} \to w}} \bigtriangleup_{\tau}^2 g(\bar x \mid v)(w^{\prime})  \,\,\mbox{with}\,\, w \in \R^n.
     $$
     We say that $g$ is twice epi-differentiable at $\bar x$ for $v \in \R^n$ if the functions $\bigtriangleup_{\tau}^2 g(\bar x \mid v)$ epi-converge to ${\rm d}^2 g(\bar x \mid v)$ as $\tau \downarrow 0$, see \cite[Definition 13.6]{Rockafellar Wets}.
     
     Now consider another kind of second-order difference quotient, which is called a parabolic difference quotient, defined by 
	 $$
	 \bigtriangleup_{\tau}^2 g(\bar x)(w \mid z) := \frac{g(\bar x+ \tau w + \frac{1}{2} \tau^2 z)-g(\bar x)-\tau {\rm d}g(\bar x)(w)}{\frac{1}{2} \tau^2} \,\,\mbox{with}\,\, w,z \in \R^n, \tau>0.
	 $$
	 The parabolic subderivative of $g$ at $\bar x$ for $w$  with ${\rm d}g(\bar x)(w)$ finite with respect to $z$ is defined by
	 $$
	 {\rm d}^2 g(\bar x)(w  \mid z) := \liminf_{ {\tau \downarrow 0} \atop{ z^{\prime} \to z}} \bigtriangleup_{\tau}^2 g(\bar x)(w \mid z^{\prime}).
	 $$
	 $g$ is said to be parabolically epi-differentiable at $\bar x$ for $w$ if ${\rm dom}\,{\rm d}^2 g(\bar x)(w  \mid \cdot) \neq \emptyset$ and for every $z \in \R^n$ and every sequence $\tau_k \downarrow 0$ there exist sequences $z_k \to z$ such that 
	 \begin{equation*}\label{parabolically epi-differentiable of g}
		{\rm d}^2 g(\bar x)(w  \mid z) = \lim_{k \to \infty} \frac{g(\bar x+ \tau_k w + \frac{1}{2} \tau_k^2 z_k)-g(\bar x)-\tau_k {\rm d}g(\bar x)(w)}{\frac{1}{2} \tau_k^2}.
	 \end{equation*}
     We say that $g$ is parabolically epi-differentiable at $\bar x$ if it satisfies the above condition at $\bar x$ for any $w$ with ${\rm d}g(\bar x)(w)$ finite. 
     It has been demonstrated in Mohammadi and Sarabi \cite[Proposition 4.1]{MS} that if $g$ is locally Lipschitz continuous around $\bar x$ relative to its domain and parabolically epi-differentiable at $\bar x$ for $w \in T_{{\rm dom}g} (\bar x)$, then ${\rm dom}\,{\rm d}^2 g(\bar x) (w \mid \cdot)= T_{{\rm dom}g}^2 (\bar x,w)$ and ${\rm dom}g$ is parabolically derivable at $\bar x$ for $w$. 
	 Below, we record an important relationship between the second subderivative and the parabolic subderivative of functions, used extensively in our paper (see \cite[Proposition 13.64]{Rockafellar Wets}).
	
	 \begin{prop}\label{prop rs of ss and ps}
	     For $g \colon \mathbb{R}^n \to \mathbb{R} \cup \{\pm \infty\}$, any point $\bar x \in \R^n$ with $g(\bar x)$ finite and any vector $w$  with ${\rm d}g(\bar x)(w)$ finite, let $v$ be such that ${\rm d}g(\bar x)(w) = \langle v,w \rangle$. Then
	   \begin{equation}\label{rs of ss and ps}
	 	{\rm d}^2 g(\bar x \mid v)(w) \leq \inf_z \{ {\rm d}^2 g(\bar x)(w  \mid z) - \langle v, z \rangle \}.
	   \end{equation}
	 \end{prop}

	 We say that a function $g \colon \mathbb{R}^n \to \mathbb{R} \cup \{\pm \infty\}$ is parabolically regular at a point $\bar x$ for a vector $v$ if $g(\bar x)$ is finite and the inequality in \eqref{rs of ss and ps} holds with equality for every $w$ satisfying ${\rm d}g(\bar x)(w) = \langle v,w \rangle$.
     The critical cone of a function $g \colon \mathbb{R}^n \to \mathbb{R} \cup \{\pm \infty\}$ at a point $\bar x$ for a vector $v$ is defined by
     $$
       K_g(\bar x,v) := \{w \in \R^n \colon {\rm d}g(\bar x)(w) = \langle v,w \rangle\}.
     $$

    \subsection{Eigenvalues and Singular Values}
	 Given an $m \times n$ matrix $X$ and index sets $I \subset \{1, \ldots,m\},\, J \subset \{1, \ldots ,n\}$,  denote by $X_{ij}$ the $(i,j)$th entry of $X$ and denote by $X_{IJ}$ the submatrix of $X$ obtained by removing all the rows of $X$ not in $I$ and all the columns of $X$ not in $J$. The matrix $X_{I}$ is the submatrix of $X$ with columns specified by $I$, unless otherwise specified. For any vector $x \in \R^n$, let ${\rm diag} (x)$ denote the matrix where $({\rm diag} (x))_{ii} = x_i$ for all $i$, and $({\rm diag} (x))_{ij} = 0$ for $i \neq j$. It's worth noting that ${\rm diag} (x)$ may represent an $m \times n$, $n \times n$, or $m \times m$ matrix (the latter when $x \in \R^m$). However, the context will always clarify which dimension is applicable. 
	 Let ${\mathbb S}^n$ be the space of all real $n \times n$ symmetric matrices. Let ${\cal O}^n$ be the set of all $n \times n$ orthogonal matrices and ${\cal O}^{m,n}$ denote the Cartesian product ${\cal O}^m \times {\cal O}^n$. 
	 The induced Frobenius norm of $X  \in {\mathbb M}_{m,n}$ is defined via the trace inner product by $\|X\| =  \sqrt{{\rm tr}(X^T X)}$. 
	 
     If $A \in {\mathbb S}^n$, then we can arrange its $n$ real eigenvalues in the decreasing order:
	 $$
	 \lambda_1(A) \geq \lambda_2(A) \geq \cdots \geq \lambda_s(A) \geq \cdots \geq \lambda_n(A),
	 $$
	 where $\lambda_s(A)$ is the sth largest eigenvalue of $A$ (counting multiplicity of each of them). 
	 For any $A \in {\mathbb S}^n$, there exists $U \in {\cal O}^n$ for which we have
	 \begin{equation}\label{SD of A}
	 	A = U \big({\rm diag}(\lambda(A))\big) U^{T}\,\, \mbox{with}\,\, \lambda(A)= \big( \lambda_1(A), \ldots, \lambda_n(A) \big).
	 \end{equation}
     For a given matrix $A \in {\mathbb S}^n$, the set of such orthogonal matrices $U$  in \eqref{SD of A} is denoted by ${\cal O}^{n}(A)$. 
     We say that two matrices $A$ and $B$ in ${\mathbb S}^n$ have a simultaneous ordered spectral decomposition if there exists $U \in {\cal O}^{n}$ such that $A = U \big({\rm diag}(\lambda(A)) \big) U^{T}$ and $B = U \big({\rm diag}(\lambda(B)) \big) U^{T}$. 
     It is well-known that any two matrices $A$ and $B$ in ${\mathbb S}^n$ satisfy the inequality
	 \begin{equation}\label{Fan’s inequality}
	 	\langle A,B \rangle \leq \langle \lambda(A), \lambda(B) \rangle,
	 \end{equation}
	 which is known as Fan’s inequality. Moreover, equality in \eqref{Fan’s inequality} amounts to $A$ and $B$ admitting a simultaneous ordered spectral 
	 decomposition. 
	 We denote by $l_s^{\lambda}$ the number of eigenvalues, ranking before $s$, which are equal to $\lambda_s(A)$ (including $\lambda_s(A)$) and $j_s^{\lambda}$ the number of eigenvalues, ranking strictly after $s$, which are equal to $\lambda_s(A)$. According to \cite[Proposition 1.4, Lemma 1.1]{Torki 2001}, we have the following result.
	 
	 \begin{lem}\label{first order diff of ED}
	 	Let $A \in {\mathbb S}^n$ and $U =[u_1, \ldots, u_n] \in {\cal O}^n$ such that 
	 	$$
	 	U^{T} A U = {\rm diag}(\lambda_1(A),  \ldots, \lambda_n(A)).
	 	$$
	 	If we set 
	 	\begin{eqnarray*}
	 		U_s &=& [u_{s-l_s^{\lambda}+1}, \ldots, u_{s+j_s^{\lambda}}],\\
	 		U_s^c &=& [u_1, \ldots, u_{s-l_s^{\lambda}}, u_{s+j_s^{\lambda}+1}, \ldots, u_n],
	 	\end{eqnarray*}
	 	then for a small perturbation matrix $E\in {\mathbb S}^n$,
	 	$$
	 	\lambda_s(A+E) = \lambda_s(A) + \lambda_{l_s^{\lambda}} \big(U_s^{T} E U_s + U_s^{T} E U_s^c(\lambda_s(A)I-\Lambda_s)^{-1} {U_s^c}^{T} E U_s \big) + O(\|E\|^3),
	 	$$
	 	where $\Lambda_s = {\rm diag} \big(\lambda_1(A), \ldots, \lambda_{s-l_s^{\lambda}}(A), \lambda_{s+j_s^{\lambda}+1}(A), \ldots, \lambda_n(A) \big)$. 
	 \end{lem}

     For any $X \in {\mathbb M}_{m,n}$ with rank $r$, there exists $(U,V) \in {\cal O}^{m,n}$ for which we have a singular value decomposition:
     \begin{equation}\label{SVD of X}
     	X = U \Sigma(X) V^{T},
     \end{equation}
     where $\Sigma(X):={\rm diag}(\sigma(X)) \in {\mathbb M}_{m,n}$ is a diagonal matrix with singular values in descending order:
     $$
     \sigma_1(X) \ge \sigma_2(X) \ge \ldots \ge \sigma_r(X) > 0 = \sigma_{r+1}(X) = \ldots = \sigma_n(X)
     $$
     on the diagonal. For a given matrix $X \in {\mathbb M}_{m,n}$, the set of such orthogonal matrices $U$ and $V$  in \eqref{SVD of X} is denoted by ${\cal O}^{m,n}(X)$. 
     We say that two matrices $X$ and $Y$ in ${\mathbb M}_{m,n}$ have a simultaneous ordered singular value decomposition if there exists $(U,V) \in {\cal O}^{m,n}$ such that $X = U \Sigma(X) V^{T}$ and $Y = U \Sigma(Y) V^{T}$. The next lemma, known as Von Neumann’s trace theorem, shows precisely when two matrices $X$ and $Y$ admit simultaneous ordered singular value decomposition (see \cite[Theorem 4.6]{LSI}).
	 
	 \begin{lem}\label{Von Neumann’s Trace Theorem}
	  Any matrices $X$ and $Y$ in ${\mathbb M}_{m,n}$ satisfy the inequality 
	 	\begin{equation}\label{Von Neumann’s Trace Inequality}
	 		 \langle X,Y \rangle \leq \langle \sigma(X), \sigma(Y) \rangle.
	 	\end{equation}
	  Equality holds if and only if $X$ and $Y$ have a simultaneous ordered singular value decomposition.
	 \end{lem}

	 It is not hard to see that any two matrices $X$ and $Y$ in ${\mathbb M}_{m,n}$, the estimate
	 \begin{equation}\label{Von Neumann’s Trace inequality}
	 	\|\sigma(X)-\sigma(Y)\| \leq \|X-Y\|
	 \end{equation}
	 always holds. 
	 Denote the three index sets $\alpha, \beta, \beta_0$ by
	 \begin{equation*}
		\alpha:=\{1, \cdots, r\}, \beta:=\{r+1,\cdots, n\},\,\mbox{and}\, \beta_0=\{n+1, \cdots, m\}.
	 \end{equation*}
     Assume that $\mu_1 > \cdots > \mu_t >\mu_{t+1}=0$ are the distinct singular values of $X \in {\mathbb M}_{m,n}$, and define the index sets
     \begin{equation}\label{alpha_i}
     	\alpha_i:=\{s \in \alpha \colon \sigma_s(X) = \mu_i \} \,\, \mbox{for all} \,\,  i=1, \ldots,t.
     \end{equation}
	 Obviously, $\alpha=\cup_{i=1}^t \alpha_i$. Define $\hat{\beta} :=\beta \cup \beta_0$.
	 Partition $U \in {\cal O}^m$ as 
	 $$U=
       \begin{bmatrix}
        	U_{\alpha_1}   &   U_{\alpha_2}   &   \cdots   &   U_{\alpha_t}  &  U_{\hat{\beta}}  
       \end{bmatrix},
     $$
     where $U_{\alpha_i} \in \R^{m\times {|\alpha_i|}}$ for $i=1, \ldots,t$, and $U_{\hat{\beta}} \in \R^{m\times {|\hat{\beta}|} }$. Similarly, 
     $$V= 
       \begin{bmatrix}
        	V_{\alpha_1}   &   V_{\alpha_2}   &   \cdots   &   V_{\alpha_t}  &  V_{\alpha_{t+1}}  
       \end{bmatrix} \in {\cal O}^n,
     $$
     where $V_{\alpha_i} \in \R^{n\times {|\alpha_i|}}$ for $i=1, \ldots,t,t+1$ with $\alpha_{t+1}:=\beta$.
	 
	 For a matrix $X\in {\mathbb M}_{m,n},$ its nuclear norm is given by
	 \begin{equation*}
	 	\|X\|_* := \sum_{i=1}^n\sigma_i(X).
	 \end{equation*}
	 Denote	${\mathbb B} := \{ X \colon \|X\|_* \le 1\}$ as the unit sublevel set of the nuclear norm, and then its polar is the closed convex set 
	 \begin{equation*}
	 	{\mathbb B}^o = \Big\{Z \colon \sigma_1(Z):=\max\limits_{1\le i \le n} \sigma_i(Z)\le 1 \Big\}.
	 \end{equation*} 
     From \cite[Example 2]{Watson  1992}, we know that the subdifferential of the nuclear norm at $X$ is
	 \begin{eqnarray*}\label{subdiff of NN}
		\partial \|X\|_* \!\!\!\!&=&\!\!\!\! \big\{ U {\rm diag}(s) V^{T} \colon X = U \Sigma(X) V^{T}, s\in \partial \|\sigma(X)\|_1 \big\}  \nonumber\\
		&=&\!\!\!\!  \big\{ U {\rm diag}(s) V^{T} \colon X = U \Sigma(X) V^{T},  s_i = 1, i\in \alpha; |s_i|\le 1, i \in \beta 
		\big\}  \nonumber\\
		&=&\!\!\!\! \bigg\{ U_{\alpha} V_{\alpha}^{T} + U_{\hat{\beta}} Z {V_{\beta}}^{T} \colon 
		X = \begin{bmatrix}
			U_{\alpha}   &   U_{\hat{\beta}}  
		    \end{bmatrix} 
	        \Sigma(X) 
	        \begin{bmatrix}
	     	V_{\alpha}   &    V_{\beta}  
	        \end{bmatrix}^{T},\,\,
	      \sigma_1(Z) \le 1 \bigg\},
	 \end{eqnarray*}
	 where ${\rm diag}(s)$ is an $m \times n$ diagonal matrix with $s$ on its diagonal and $Z \in \R^{(m-r) \times (n-r)}$.

     Define the linear operator ${\cal B}(\cdot):{\mathbb M}_{m,n} \rightarrow {\mathbb S}^{m+n}$ by
	 \begin{equation*}
		{\cal B}(X):=
		\begin{bmatrix}
			0        &  X \\
			X^{T} &  0
		\end{bmatrix}, \,\,
		X \in {\mathbb M}_{m,n}.
	 \end{equation*}
	 It follows from \cite[Theorem 7.3.3]{Horn} that 
	 \begin{equation}\label{eigen-singular}
		P^{T} {\cal B}(X) P= 
		\begin{bmatrix}
			\Sigma_n(X) &   0         &  0  \\
		    	0       &   0_{m-n}   &  0  \\
			    0       &   0         &  -\Sigma_n(X)
		\end{bmatrix},
	 \end{equation}
	 where the orthogonal matrix $P \in {\cal O}^{m+n}$ is given by
	 \begin{equation*}
		P:=\frac{1}{\sqrt{2}}
		\begin{bmatrix}
			U_{\alpha}   &   U_{\beta}   &   \sqrt{2} U_{\beta_0}   &   -U_{\alpha}  &  -U_{\beta}   \\
			V_{\alpha}   &   V_{\beta}   &            0             &  	V_{\alpha}   &   V_{\beta} 
		\end{bmatrix}.
	 \end{equation*}
	 Obviously, for each $s \in \alpha \cup \beta$, we have that $\sigma_s(X) = \lambda_s({\cal B}(X))$. Therefore, by applying the definition of directional derivative, we derive that $\sigma_s^{\prime}(X;H) = \lambda_s^{\prime}({\cal B}(X); {\cal B}(H))$ for any $H \in {\mathbb M}_{m,n}$.
	 
     Similar to the symmetric case, for matrix $X_0 \in {\mathbb M}_{m,n}$, we denote by $l_s(X_0)$ the number of singular values, ranking before $s$, which are equal to $\sigma_s(X_0)$ (including $\sigma_s(X_0)$) and $j_s(X_0)$ the number of singular values, ranking strictly after $s$, which are equal to $\sigma_s(X_0)$, respectively, i.e., we define $l_s(X_0)$ and $j_s(X_0)$ such that
	 \begin{align*}
		\sigma_1(X_0) \geq \cdots \geq \sigma_{s-l_s(X_0)}(X_0) &> \sigma_{s-l_s(X_0)+1}(X_0) = \cdots = \sigma_s(X_0) = \cdots = \sigma_{s+j_s(X_0)}(X_0)  \\
		&> \sigma_{s+j_s(X_0)+1}(X_0) \geq \cdots \geq \sigma_n(X_0).
	 \end{align*}
	 We use $r_s(X_0)$ to denote the multiplicity of $\sigma_s(X_0)$. In later discussions, when the dependence of $l_s$, $j_s$ and $r_s$, $s \in \alpha \cup \beta$, on $X_0$ can be seen clearly from the context, we often drop $X_0$ from these notations. The following proposition on the directional derivatives of the singular value of a matrix was explored in \cite[Section 5.1]{LSII}.
	
	 \begin{prop}\label{fdd of SV M}
		Let $X_0 \in {\mathbb M}_{m,n}$ of rank $r$ be given and have the singular value decomposition \eqref{SVD of X}. Let $H \in {\mathbb M}_{m,n}$ be a small perturbation matrix. Then for each $s \in \alpha_i, i=1,\ldots,t$, 
		\begin{equation*}\label{fdd of alpha}
			\sigma^{\prime}_s(X_0;H) = \frac{1}{2} \lambda_{l_s} \big(U_{\alpha_i}^T H V_{\alpha_i} + V_{\alpha_i}^T H^T U_{\alpha_i} \big).
		\end{equation*}
		For each $s \in \beta$, one has that
		\begin{equation*}\label{fdd of beta}
			\sigma^{\prime}_s(X_0;H) = \sigma_{l_s} \big(U_{\hat{\beta}}^T H V_{\beta} \big).
		\end{equation*}
	 \end{prop}
	
	 Assume that $\eta^i_1> \cdots> \eta^i_{N_i}$ are the distinct eigenvalues of $\frac{1}{2} \big(U_{\alpha_i}^T H V_{\alpha_i} + V_{\alpha_i}^T H^T U_{\alpha_i} \big)$, and define the index sets
	 \begin{equation}\label{def of beta^i_j}
		\beta^i_j := \bigg\{s \in \alpha_i \colon \lambda_s \Big(\frac{1}{2} \big(U_{\alpha_i}^T H V_{\alpha_i} + V_{\alpha_i}^T H^T U_{\alpha_i} \big) \Big)=\eta^i_j \bigg\} \,\,\mbox{for all}\,\, i=1, \ldots,t,\, j=1,\ldots,N_i.
	 \end{equation}
	 Denote the distinct singular values of $U_{\hat{\beta}}^T H V_{\beta}$ by $\eta^{t+1}_1> \cdots> \eta^{t+1}_{N_{t+1}}> \eta^{t+1}_{N_{t+1}+1}=0$, and define
	 \begin{equation}\label{def of beta^{t+1}_j}
		\beta^{t+1}_k :=\bigg\{s \in \beta \colon \sigma_s \big( U_{\hat{\beta}}^T H V_{\beta} \big)=\eta^{t+1}_k \bigg\} \,\,\mbox{for all}\,\, k=1, \ldots, N_{t+1}, N_{t+1}+1.
	 \end{equation}
     For each $s \in \{1, \ldots, n\}$, there exists $i \in \{ 1, \ldots, t, t+1 \}$ such that $s \in \alpha_i$ and $l_s(X_0) \in \{1, \ldots, |\alpha_i|\}$.  Furthermore, we can find $j \in \{1, \ldots, N_i, N_{t+1}+1\}$ such that $l_s(X_0) \in \beta^i_j$. Define now the integer $\tilde{l}_s(X_0,H)$ by
     \begin{align*}
     \tilde{l}_s(X_0, H) := \left \{
     	\begin{array}{lc}
        l_{l_s(X_0)} \Big( \frac{1}{2} \big(U_{\alpha_i}^T H V_{\alpha_i} + V_{\alpha_i}^T H^T U_{\alpha_i} \big) \Big)  &\,\,\mbox{if}\,\,s \in \alpha,\\
        l_{l_s(X_0)} \big( U_{\hat{\beta}}^T H V_{\beta} \big) & \,\,\mbox{if}\,\,s \in \beta,
        \end{array}\right.
     \end{align*}
     which, in fact, signifies the number of eigenvalues of $\frac{1}{2} \big(U_{\alpha_i}^T H V_{\alpha_i} + V_{\alpha_i}^T H^T U_{\alpha_i} \big)$  that are equal to $\lambda_{l_s(X_0)} \Big( \frac{1}{2} \big(U_{\alpha_i}^T H V_{\alpha_i} + V_{\alpha_i}^T H^T U_{\alpha_i} \big) \Big)$, but are ranked before $\lambda_{l_s(X_0)} \Big( \frac{1}{2} \big(U_{\alpha_i}^T H V_{\alpha_i} + V_{\alpha_i}^T H^T U_{\alpha_i} \big) \Big)$, including $\lambda_{l_s(X_0)} \Big( \frac{1}{2} \big(U_{\alpha_i}^T H V_{\alpha_i} + V_{\alpha_i}^T H^T U_{\alpha_i} \big) \Big)$ for $s \in \alpha$, and the number of singular values of $U_{\hat{\beta}}^T H V_{\beta}$  that are equal to $\sigma_{l_s(X_0)} \big( U_{\hat{\beta}}^T H V_{\beta} \big)$, but are ranked before $\sigma_{l_s(X_0)} \big( U_{\hat{\beta}}^T H V_{\beta} \big)$, including $\sigma_{l_s(X_0)} \big( U_{\hat{\beta}}^T H V_{\beta} \big)$ for $s \in \beta$. 
     As before, we often drop $X_0$ and $H$ from $\tilde{l}_s(X_0, H)$ when the dependence of $\tilde{l}_s$ on $X_0$ and $H$ can be seen clearly from the context. 

	 Similarly, \cite[Theorem 3.1]{ZZX} derives the following explicit formulas of the (parabolic) second order directional derivatives of the singular values. 

     \begin{prop}\label{sdd of SV M}
	 Let $X_0 \in {\mathbb M}_{m,n}$ of rank $r$ be given and have the singular value decomposition \eqref{SVD of X}. Suppose that the direction $H,W \in {\mathbb M}_{m,n}$. Then the following assertions hold:
	 
	 {\rm (i)} For each $s \in \alpha$, there exists $i \in \{ 1, \ldots, t\}$ and $j \in \{1, \ldots,N_i\}$ such that  $s \in \alpha_i$ and $l_s(X_0) \in \beta^i_j$, and there exists $Q^i \in {\cal O}^{|\alpha_i|} \Big( \frac{1}{2} \big(U_{\alpha_i}^T H V_{\alpha_i} + V_{\alpha_i}^T H^T U_{\alpha_i} \big) \Big)$ such that 
	 \begin{equation}\label{sdd of alpha}
		\sigma^{\prime\prime}_s(X_0;H,W) = \lambda_{\tilde{l}_s} \bigg( \Big(Q^i_{\beta^i_j} \Big)^T P_{\alpha_i}^T \Big[{\cal B}(W)-2{\cal B}(H) P_{\alpha_i}^c (\varLambda_s- \mu_i I)^{-1} {P_{\alpha_i}^c}^{T} {\cal B}(H) \Big] P_{\alpha_i} Q^i_{\beta^i_j} \bigg),
	 \end{equation}
     where $\varLambda_s \in \R^{(m+n-r_s)\times (m+n-r_s)}$ is a diagonal matrix whose diagonal elements are the eigenvalues of $P^{T} {\cal B}(X_0) P$ that are not equal to $\sigma_s(X_0)$. 
     
     {\rm (ii)} For each $s \in \beta$, if $l_s(X_0) \in \beta^{t+1}_k$ for some $k \in \{1, \ldots,N_{t+1}\}$, then there exists $( Q_{\hat{\beta} \hat{\beta}}, \widehat{Q}_{\beta \beta} ) \in {\cal O}^{|\hat{\beta}|, |\beta|} \big(U_{\hat{\beta}}^T H V_{\beta} \big)$ such that
	 \begin{align}\label{sdd of beta >0}
		\sigma^{\prime\prime}_s(X_0;H,W) 
	  = &\frac{1}{2} \lambda_{\tilde{l}_s} \bigg( \Big(Q_{\hat{\beta} \beta^{t+1}_k} \Big)^T \Big[ U_{\hat{\beta}}^T W V_{\beta}-2 U_{\hat{\beta}}^{T} H V_{\alpha} \Sigma_{\alpha}^{-1}(X_0) U_{\alpha}^{T} H V_{\beta} \Big] \widehat{Q}_{\beta \beta^{t+1}_k} \nonumber\\
		& + \Big(\widehat{Q}_{\beta \beta^{t+1}_k} \Big)^T \Big[ U_{\hat{\beta}}^T W V_{\beta}-2 U_{\hat{\beta}}^{T} H V_{\alpha} \Sigma_{\alpha}^{-1}(X_0) U_{\alpha}^{T} H V_{\beta} \Big]^T Q_{\hat{\beta} \beta^{t+1}_k} \bigg),
	 \end{align}
     and if $l_s(X_0) \in \beta^{t+1}_{N_{t+1}+1}$, then there exists $( Q_{\hat{\beta} \hat{\beta}}, \widehat{Q}_{\beta \beta} ) \in {\cal O}^{|\hat{\beta}|, |\beta|} \big(U_{\hat{\beta}}^T H V_{\beta} \big)$ such that
     \begin{equation}\label{sdd of beta =0}
     	\sigma^{\prime\prime}_s(X_0;H,W) = \sigma_{\tilde{l}_s} \left( 
     	\begin{bmatrix}
     		Q_{\hat{\beta} \beta^{t+1}_{N_{t+1}+1}}  & Q_{\hat{\beta} \beta_0}
     	\end{bmatrix}^T
      \Big[ U_{\hat{\beta}}^T W V_{\beta}-2 U_{\hat{\beta}}^{T} H V_{\alpha} \Sigma_{\alpha}^{-1}(X_0) U_{\alpha}^{T} H V_{\beta} \Big] \widehat{Q}_{\beta \beta^{t+1}_{N_{t+1}+1}}
     	 \right).
     \end{equation}
    \end{prop}

    Combining Proposition \ref{fdd of SV M} with Proposition \ref{sdd of SV M}, we obtain the following result, which is important for our development in this paper.
    \begin{cor} 
     Let $X_0 \in {\mathbb M}_{m,n}$ of rank $r$ be given and have the singular value decomposition \eqref{SVD of X}. Suppose that the direction $H,W \in {\mathbb M}_{m,n}$. Then for any $t>0$ sufficiently small, we have
	 \begin{equation}\label{estimate for the SVF}
		\sigma \left(X_0+tH+\frac{1}{2}t^2 W +o(t^2) \right) = \sigma(X_0) + t\sigma^{\prime}(X_0;H) + \frac{1}{2}t^2 \sigma^{\prime\prime}(X_0;H,W) +o(t^2).
	 \end{equation}
    \end{cor}

	\section{Twice Epi-Differentiability of the Nuclear Norm}
	 In this section, we study the first- and second-order epi-differentiability of the nuclear norm. The subsequent proposition is a direct consequence of Lemma \ref{first order diff of ED}.
	 \begin{prop}\label{fdd of SV}
		Let $X_0 \in {\mathbb M}_{m,n}$ of rank $r$ be given and have the singular value decomposition \eqref{SVD of X}. Let $H \in {\mathbb M}_{m,n}$ be a small perturbation matrix. Then for $s \in \alpha \cup \beta$, we have
		$$
		\sigma_s(X_0+H) = \sigma_s(X_0) + \lambda_{l_s} \big(P_s^{T} {\cal B}(H) P_s + P_s^{T} {\cal B}(H) P_s^c (\sigma_s(X_0)I-\varLambda_s)^{-1} {P_s^c}^{T} {\cal B}(H) P_s \big) + O \big(\|{\cal B}(H)\|^3 \big),
		$$
		where the columns of $P_s$ form an orthonormal basis of eigenvectors of ${\cal B}(X_0) $ associated with $\sigma_s(X_0)$, $P_s^c$ is the submatrix of $P$ obtained by removing all the columns of $P_s$, and $\varLambda_s \in \R^{(m+n-r_s)\times (m+n-r_s)}$ is a diagonal matrix whose diagonal elements are the eigenvalues of $P^{T} {\cal B}(X_0) P$ that are not equal to $\sigma_s(X_0)$.
		
		Moreover, for any $\tau >0$ we have the following results:
		
		{\rm (i)} For any $s \in \alpha$ there exists $i \in \{ 1, \ldots, t\}$ and $j \in \{1, \ldots,N_i\}$ such that  $s \in \alpha_i$ and $l_s(X_0) \in \beta^i_j$, and we have
		\begin{align}\label{second order expansion}
			\sigma_{s}(X_0+\tau H) =&\sigma_{s}(X_0) + \tau \lambda_{l_s} \big(P_{\alpha_i}^{T} {\cal B}(H) P_{\alpha_i} \big) + \nonumber\\
			&\frac{{\tau}^2}{2}\lambda_{\tilde{l}_s} \bigg(2 \Big(Q^i_{\beta^i_j} \Big)^{T} P_{\alpha_i}^{T} {\cal B}(H) P_{\alpha_i}^c (\sigma_s(X_0)I-\varLambda_s)^{-1} {P_{\alpha_i}^c}^{T} {\cal B}(H) P_{\alpha_i} Q^i_{\beta^i_j} \bigg) + O({\tau}^3),
		\end{align}
		where  $Q^i \in {\cal O}^{|\alpha_i|} \big(P_{\alpha_i}^{T} {\cal B}(H) P_{\alpha_i} \big)$. 
		
		{\rm (ii)} If $n>r$, then for $s \in \beta$ we have
		\begin{equation*}
			\sigma_{s}(X_0+H) = \sigma_{s}(X_0) + \sigma_{l_s} \big(U_{\hat{\beta}}^{T} H V_{\beta} - U_{\hat{\beta}}^{T} H V_{\alpha} \Sigma_{\alpha}^{-1}(X_0) U_{\alpha}^{T} H V_{\beta} \big)+ O(\|{\cal B}(H)\|^3).
		\end{equation*}
        Furthermore, if $l_s(X_0) \in \beta^{t+1}_k$ for some $k \in \{1, \ldots,N_{t+1}\}$, then there exists $( Q_{\hat{\beta} \hat{\beta}}, \widehat{Q}_{\beta \beta} ) \in {\cal O}^{|\hat{\beta}|, |\beta|} \big(U_{\hat{\beta}}^T H V_{\beta} \big)$ such that
		\begin{align*}
		\sigma_{s}(X_0+\tau H) = \sigma_{s}(X_0) + &\tau\sigma_{l_s} \big(U_{\hat{\beta}}^{T} H V_{\beta}\big) + \frac{{\tau}^2}{4} \lambda_{\tilde{l}_s} \bigg( \Big(Q_{\hat{\beta} \beta^{t+1}_k} \Big)^T \Big[ - U_{\hat{\beta}}^{T} H V_{\alpha} \Sigma_{\alpha}^{-1}(X_0) U_{\alpha}^{T} H V_{\beta} \Big] \widehat{Q}_{\beta \beta^{t+1}_k} \nonumber\\
			+& \Big(\widehat{Q}_{\beta \beta^{t+1}_k} \Big)^T \Big[- U_{\hat{\beta}}^{T} H V_{\alpha} \Sigma_{\alpha}^{-1}(X_0) U_{\alpha}^{T} H V_{\beta}\Big]^T Q_{\hat{\beta} \beta^{t+1}_k} \bigg)+ O({\tau}^3),
		\end{align*}
		and if $l_s(X_0) \in \beta^{t+1}_{N_{t+1}+1}$, then there exists $( Q_{\hat{\beta} \hat{\beta}}, \widehat{Q}_{\beta \beta} ) \in {\cal O}^{|\hat{\beta}|, |\beta|} \big(U_{\hat{\beta}}^T H V_{\beta} \big)$ such that
		\begin{align*}
		\sigma_{s}(X_0+\tau H) = &\sigma_{s}(X_0) + \tau\sigma_{l_s} \big(U_{\hat{\beta}}^{T} H V_{\beta}\big)  \\
		&+ \frac{{\tau}^2}{2}\sigma_{\tilde{l}_s} \left( 
		\begin{bmatrix}
		Q_{\hat{\beta} \beta^{t+1}_{N_{t+1}+1}}  & Q_{\hat{\beta} \beta_0}
    	\end{bmatrix}^T
        \Big[- U_{\hat{\beta}}^{T} H V_{\alpha} \Sigma_{\alpha}^{-1}(X_0) U_{\alpha}^{T} H V_{\beta}\Big] 
		\widehat{Q}_{\beta \beta^{t+1}_{N_{t+1}+1}}
		 \right) + O({\tau}^3). 
		\end{align*}
	 \end{prop}
	
	 {\it Proof} It follows from Lemma \ref{first order diff of ED} and \eqref{eigen-singular} that
	  \begin{eqnarray*}
		\!\!\!\!& &\!\!\!\!	\sigma_s(X_0+H) = \lambda_s({\cal B}(X_0+H)) = \lambda_s({\cal B}(X_0) + {\cal B}(H)) \\
		\!\!\!\!&=&\!\!\!\! \lambda_s({\cal B}(X_0)) + \lambda_{l_s} \big(P_s^{T} {\cal B}(H) P_s + P_s^{T} {\cal B}(H) {P_s^c} (\lambda_s({\cal B}(X_0))I-\varLambda_s)^{-1} {P_s^c}^{T} {\cal B}(H) P_s \big) + O(\|{\cal B}(H)\|^3)  \\
		\!\!\!\!&=&\!\!\!\!  \sigma_s(X_0) + \lambda_{l_s} \big(P_s^{T} {\cal B}(H) P_s + P_s^{T} {\cal B}(H) P_s^c (\sigma_s(X_0)I-\varLambda_s)^{-1} {P_s^c}^{T} {\cal B}(H) P_s \big) + O(\|{\cal B}(H)\|^3).
	  \end{eqnarray*}

	 For any $s \in \alpha$, there exists $i \in \{ 1, \ldots, t\}$ and $j \in \{1, \ldots,N_i\}$ such that  $s \in \alpha_i$ and $l_s(X_0) \in \beta^i_j$. By using Lemma \ref{first order diff of ED} with $A =P_{\alpha_i}^{T} {\cal B}(H) P_{\alpha_i}$ and $E =P_{\alpha_i}^{T} {\cal B}(H) P_{\alpha_i}^c (\sigma_s(X_0)I-\varLambda_s)^{-1} {P_{\alpha_i}^c}^{T} {\cal B}(H) P_{\alpha_i}$, we have
	 \begin{eqnarray*}
		\sigma_{s}(X_0+\tau H) 
		\!\!\!\!&=&\!\!\!\!  \sigma_s(X_0) + \lambda_{l_s} \big(\tau P_{\alpha_i}^{T} {\cal B}(H) P_{\alpha_i} + {\tau}^2 P_{\alpha_i}^{T} {\cal B}(H) P_{\alpha_i}^c (\sigma_s(X_0)I-\varLambda_s)^{-1} {P_{\alpha_i}^c}^{T} {\cal B}(H) P_{\alpha_i} \big) + O({\tau}^3)  \\
		\!\!\!\!&=&\!\!\!\!  \sigma_s(X_0) + \tau \lambda_{l_s} \big( P_{\alpha_i}^{T} {\cal B}(H) P_{\alpha_i} + {\tau} P_{\alpha_i}^{T} {\cal B}(H) P_{\alpha_i}^c (\sigma_s(X_0)I-\varLambda_s)^{-1} {P_{\alpha_i}^c}^{T} {\cal B}(H) P_{\alpha_i} \big) + O({\tau}^3)  \\
		\!\!\!\!&=&\!\!\!\!  \sigma_s(X_0) + \tau \bigg[\lambda_{l_s}\big(P_{\alpha_i}^{T} {\cal B}(H) P_{\alpha_i} \big) + \\
		\!\!\!\!& &\!\!\!\!  {\tau}\lambda_{\tilde{l}_s} \bigg( \Big(Q^i_{\beta^i_j} \Big)^{T} P_{\alpha_i}^{T} {\cal B}(H) P_{\alpha_i}^c (\sigma_s(X_0)I-\varLambda_s)^{-1} {P_{\alpha_i}^c}^{T} {\cal B}(H) P_{\alpha_i} Q^i_{\beta^i_j} \bigg) \bigg] + O({\tau}^3)  \\
		\!\!\!\!&=&\!\!\!\!  \sigma_{s}(X_0) + \tau \lambda_{l_s} \big(P_{\alpha_i}^{T} {\cal B}(H) P_{\alpha_i} \big) + \nonumber\\
		\!\!\!\!& &\!\!\!\!  \frac{{\tau}^2}{2}\lambda_{\tilde{l}_s} \bigg(2 \Big(Q^i_{\beta^i_j} \Big)^{T} P_{\alpha_i}^{T} {\cal B}(H) P_{\alpha_i}^c (\sigma_s(X_0)I-\varLambda_s)^{-1} {P_{\alpha_i}^c}^{T} {\cal B}(H) P_{\alpha_i} Q^i_{\beta^i_j} \bigg) + O({\tau}^3).
	 \end{eqnarray*}
	In particular, we have 
	$$
	\sigma^{\prime}_s(X_0;H) = \lambda_{l_s} \big(P_{\alpha_i}^{T} {\cal B}(H) P_{\alpha_i} \big) = \frac{1}{2} \lambda_{l_s} \big(U_{\alpha_i}^T H V_{\alpha_i} + V_{\alpha_i}^T H^T U_{\alpha_i} \big).
	$$
	
	 If $n>r$, then $\beta \neq \emptyset$. For any $s \in \beta$, we have $\sigma_s(X_0)=0$,
	 \begin{equation*}
		P_s=\frac{1}{\sqrt{2}}
		\begin{bmatrix}
			U_{\beta}   &   \sqrt{2} U_{\beta_0}   &    -U_{\beta}      \\
			V_{\beta}   &            0             & 	V_{\beta}   
		\end{bmatrix}, \,\,\mbox{and}\,\,
		\varLambda_s = 
		\begin{bmatrix}
			\Sigma_{\alpha}(X_0)    &           0            \\
			0                 &     -\Sigma_{\alpha}(X_0)
		\end{bmatrix}.
	 \end{equation*}
	 By some elementary calculations, we can obtain that 
	 $$
	   P_s^{T} {\cal B}(H) P_s = S {\cal B}(V_{\beta}^{T} H^{T} U_{\hat{\beta}}) S^{T},
	 $$
	 $$
	   P_s^{T} {\cal B}(H) P_s^c (\sigma_s(X_0)I-\varLambda_s)^{-1} {P_s^c}^{T} {\cal B}(H) P_s = 
	   -S{\cal B}(V_{\beta}^{T} H^{T} U_{\alpha}  \Sigma_{\alpha}^{-1}(X_0)  V_{\alpha}^{T} H^{T} U_{\hat{\beta}}) S^{T},
	 $$
	 where 
	  \begin{equation*}
		S=\frac{1}{2}
		\begin{bmatrix}
			\sqrt{2}I_{|\beta|}   &            0             &   \sqrt{2}I_{|\beta|}      \\
			0                     &     2 I_{|\beta_0|}      & 	          0                \\
			\sqrt{2}I_{|\beta|}   &            0             & 	  -\sqrt{2}I_{|\beta|}
		\end{bmatrix} \in {\cal O}^{2|\beta|+|\beta_0|}.
	  \end{equation*}
	 Therefore, we have
	  \begin{align*}
		  &  \lambda_{l_s} \big(P_s^{T} {\cal B}(H) P_s + P_s^{T} {\cal B}(H) P_s^c (\sigma_s(X_0)I-\varLambda_s)^{-1} {P_s^c}^{T} {\cal B}(H) P_s \big)  \\
	    = & \lambda_{l_s} \Big(S{\cal B} \big(V_{\beta}^{T} H^{T} U_{\hat{\beta}} \big)S^{T} -S{\cal B}\big(V_{\beta}^{T} H^{T} U_{\alpha}  \Sigma_{\alpha}^{-1}(X_0)  V_{\alpha}^{T} H^{T} U_{\hat{\beta}} \big) S^{T} \Big)   \\
		= &\lambda_{l_s} \Big(S{\cal B} \big(V_{\beta}^{T} H^{T} U_{\hat{\beta}} - V_{\beta}^{T} H^{T} U_{\alpha}  \Sigma_{\alpha}^{-1}(X_0)  V_{\alpha}^{T} H^{T} U_{\hat{\beta}} \big) S^{T} \Big)  \\
		= & \lambda_{l_s} \Big({\cal B} \big(V_{\beta}^{T} H^{T} U_{\hat{\beta}} - V_{\beta}^{T} H^{T} U_{\alpha}  \Sigma_{\alpha}^{-1}(X_0)  V_{\alpha}^{T} H^{T} U_{\hat{\beta}} \big) \Big)  \\
		= & \sigma_{l_s} \Big(V_{\beta}^{T} H^{T} U_{\hat{\beta}} - V_{\beta}^{T} H^{T} U_{\alpha}  \Sigma_{\alpha}^{-1}(X_0)  V_{\alpha}^{T} H^{T} U_{\hat{\beta}} \Big)  \\
		= & \sigma_{l_s} \Big(U_{\hat{\beta}}^{T} H V_{\beta} - U_{\hat{\beta}}^{T} H V_{\alpha} \Sigma_{\alpha}^{-1}(X_0) U_{\alpha}^{T} H V_{\beta} \Big),
	  \end{align*}
	 and then $\sigma_{s}(X_0+H) = \sigma_{s}(X_0) + \sigma_{l_s} \big(U_{\hat{\beta}}^{T} H V_{\beta} \big)+ O(\|{\cal B}(H)\|^2)$. In particular, we have $\sigma^{\prime}_s(X_0;H) = \sigma_{l_s} \big(U_{\hat{\beta}}^T H V_{\beta} \big)$. 
	 Similar to the discussions of first-order directional derivatives of singular values, we can obtain their second-order directional derivatives.
	 \qed
	
     Define the sum of all the singular values that are equal to $\sigma_s(X_0)$ but are ranked before $\sigma_s(X_0)$ by 
     $$
       \Psi_s(X_0) = \sigma_{s-l_s+1}(X_0) + \cdots +\sigma_s(X_0),
     $$
     and the sum of the first $s$-largest singular values of $X_0$ by
     $$
      \Phi_s(X_0) = \sigma_1(X_0) + \cdots +\sigma_s(X_0).
     $$
     	
	 \begin{lem}
		Let $X_0 \in {\mathbb M}_{m,n}$ of rank $r$ satisfy $n>r$ and $H \in {\mathbb M}_{m,n}$ be a small perturbation matrix. Then, we have
		\begin{equation}\label{alpha_t 15}
			\Psi_n(X_0+H) = \sigma_{{\mathbb B}^o} \left(U_{\hat{\beta}}^{T} H V_{\beta} - U_{\hat{\beta}}^{T} H V_{\alpha} \Sigma_{\alpha}^{-1}(X_0) U_{\alpha}^{T} H V_{\beta} \right)+ O \big(\|{\cal B}(H)\|^3 \big),
		\end{equation}
		where $\sigma_{{\mathbb B}^o}$ is the support function of the set 
		${\mathbb B}^{o} \subset {\mathbb R}^{|\beta| \times |\hat{\beta}|}$. 
		As a consequence,
		\begin{equation}\label{alpha_t 16}
			\Psi_n(X_0+H) = \sigma_{{\mathbb B}^{o}} \big(U_{\hat{\beta}}^{T} H V_{\beta} \big)+ O \big(\|{\cal B}(H)\|^2 \big).
		\end{equation}	
	 \end{lem}
	
	 {\it Proof} From Proposition \ref{fdd of SV}, we obtain
	  \begin{align*}
		\Psi_n(X_0+H) 
		= & \sigma_{r+1}(X_0+H) + \cdots + \sigma_n(X_0+H)   \\
		= & \sigma_{r+1}(X_0) + \sigma_1 \left(U_{\hat{\beta}}^{T} H V_{\beta} - U_{\hat{\beta}}^{T} H V_{\alpha} \Sigma_{\alpha}^{-1}(X_0) U_{\alpha}^{T} H V_{\beta} \right)  + \cdots + \\
		  & \sigma_n(X_0) + \sigma_{l_n} \left(U_{\hat{\beta}}^{T} H V_{\beta} - U_{\hat{\beta}}^{T} H V_{\alpha} \Sigma_{\alpha}^{-1}(X_0) U_{\alpha}^{T} H V_{\beta} \right ) + O \big(\|{\cal B}(H)\|^3 \big)  \\
		= & (\sigma_1 + \cdots +\sigma_{l_n}) \left(U_{\hat{\beta}}^{T} H V_{\beta} - U_{\hat{\beta}}^{T} H V_{\alpha} \Sigma_{\alpha}^{-1}(X_0) U_{\alpha}^{T} H V_{\beta} \right )+ O \big(\|{\cal B}(H)\|^3 \big)  \\
		= & \sigma_{{\mathbb B}^{o}} \left(U_{\hat{\beta}}^{T} H V_{\beta} - U_{\hat{\beta}}^{T} H V_{\alpha} \Sigma_{\alpha}^{-1}(X_0) U_{\alpha}^{T} H V_{\beta} \right ) + O \big(\|{\cal B}(H)\|^3 \big).
	 \end{align*}
	 \qed

	 \begin{lem}\label{second-order differential of Phi_r}
		The sum of the first $r$-largest singular values $\Phi_r $ is $\mathcal{C}^2$ on the set of $m \times n$  matrices of rank $r$. Let $X_0 \in {\mathbb M}_{m,n}$ of rank $r$ be given and $H \in {\mathbb M}_{m,n}$ be a matrix, the second-order differential of $\Phi_r $ at $X_0$ in the direction $H$ is given by
		\begin{equation*}
			D^2 \Phi_r(X_0)(H,H) = \sum_{k_l \leq r} {\rm tr} \Big(2 P_{k_l}^{T} {\cal B}(H) P_{k_l}^c (\sigma_{k_l}(X_0)I-\varLambda_{k_l})^{-1} {P_{k_l}^c}^{T} {\cal B}(H)P_{k_l} \Big).
		\end{equation*}
	 \end{lem}
	
	 {\it Proof} It follows from \cite[Proposition 8]{DST} that $\Phi_r $ is $\mathcal{C}^2$ on the set of $m \times n$  matrices of rank $r$.
	  For any $s \in \alpha$, it follows from \eqref{second order expansion} that the second-order directional derivative of $\sigma_s$ at $X_0$ in the direction $H$ is
	  \begin{align*}
		\sigma_{s}^{\prime\prime}(X_0; H)
		&=\lim_{\tau \downarrow 0} \frac{\sigma_{s}(X_0+\tau H)-\sigma_{s}(X_0)-\tau \sigma_s^{\prime}(X_0; H)}{\frac{\tau^2}{2}}\\
		&=\lambda_{\tilde{l}_s} \Big(2 Q_{l_s}^{T} P_s^{T} {\cal B}(H) P_s^c (\sigma_s(X_0)I-\varLambda_s)^{-1} {P_s^c}^{T} {\cal B}(H) P_s Q_{l_s} \Big).
	  \end{align*}
	  Denote by $k_1, \ldots, k_q$ the set of the indices of all the leading singular values of $X_0$, i.e.
	  $$
	   \sigma_{k_1}(X_0) > \sigma_{k_2}(X_0) > \cdots > \sigma_{k_q}(X_0).
	  $$
	  Assume that there exists $k_l \leq r$ such that
	  $$
	   \sigma_{k_l-1}(X_0) > \sigma_{k_l}(X_0) =\cdots = \sigma_{k_l+j_{k_l}}(X_0)>  \sigma_{k_{l+1}}(X_0).
	  $$
	  Set $G_{k_l} := P_{k_l}^{T} {\cal B}(H) P_{k_l}^c (\sigma_{k_l}(X_0)I-\varLambda_{k_l})^{-1} {P_{k_l}^c}^{T} {\cal B}(H)P_{k_l}$. 
	  Then for any $q \in \{k_l, k_l+1, \ldots,k_l+j_{k_l}\}$, we have
	  $$
	   \Big(2 Q_{l_q}^{T}G_{k_l} Q_{l_q} \Big) v_{l_q} = \lambda_{\tilde{l}_q} \Big(2 Q_{l_q}^{T}G_{k_l} Q_{l_q} \Big) v_{l_q},
	  $$
	  where $v_{l_q}$ is a unit eigenvector of $2 Q_{l_q}^{T}G_{k_l} Q_{l_q}$ associated with eigenvalue $\lambda_{\tilde{l}_q}(2 Q_{l_q}^{T}G_{k_l} Q_{l_q})$.  Therefore, by some elementary calculations, we obtain that 
	  \begin{align*}
		 & 	D^2 \Phi_r(X_0)(H,H) = \sum_{k_l \leq r} (\sigma_{k_l}+\sigma_{k_l+1}+\ldots+ \sigma_{k_l+j_{k_l}})^{\prime\prime}(X_0; H)\\
		=&  \sum_{k_l \leq r} \lambda_{\tilde{l}_{k_l}} \Big(2 Q_{l_{k_l}}^{T}G_{k_l} Q_{l_{k_l}}\Big)+
		\lambda_{\tilde{l}_{k_l+1}} \Big(2 Q_{l_{k_l+1}}^{T} G_{k_l} Q_{l_{k_l+1}} \Big) +\cdots + \lambda_{\tilde{l}_{k_l+j_{k_l}}} \Big(2 Q_{l_{k_l+j_{k_l}}}^{T} G_{k_l} Q_{l_{k_l+j_{k_l}}} \Big) \\
		=&  \sum_{k_l \leq r} {\rm tr} \Big(2Q^{T} G_{k_l} Q \Big) = \sum_{k_l \leq r} {\rm tr} \Big(2 G_{k_l} \Big),
	  \end{align*}
	  and the proof is complete.
	  \qed

	 Recall from Rockafellar and Wets \cite[Definition 7.20]{Rockafellar Wets} that a function $g \colon \mathbb{R}^n \to \mathbb{R} \cup \{\pm \infty\}$ is said to be semidifferentiable at $\bar x \in \R^n$ with $g(\bar x)$ finite for $w \in \R^n$, if the (possibly infinite) limit
	 \begin{equation*}
	 	\lim_{ {t \downarrow 0} \atop{ w^{\prime} \to w}} \frac{g(\bar x + t w^{\prime})-g(\bar x)}{t}
	 \end{equation*}
     exists. This limit is the semiderivative of $g$ at $\bar x$ for $w$.
	
	 \begin{lem}
		Let $X_0 \in {\mathbb M}_{m,n}$ of rank $r$ satisfy $n>r$ and $H \in {\mathbb M}_{m,n}$ be a small perturbation matrix. The function $\Psi_n$ is semidifferentiable at $X_0$ for $H$ and its semiderivative coincides with subderivative, i.e. $\Psi_n^{\prime}(X_0; H) = {\rm d}	\Psi_n(X_0)(H).$      
		The regular subdifferential of $\Psi_n$ at $X_0$ is given by
		\begin{equation*}\label{subdiff of alpha_t}
			\hat{\partial} \Psi_n(X_0) =\left\{ U_{\hat{\beta}} Z V_{\beta}^{T} \colon 
			X_0 = \begin{bmatrix}
				U_{\alpha}   &   U_{\hat{\beta}}  
			\end{bmatrix} 
		    \Sigma(X_0) 
		    \begin{bmatrix}
			    V_{\alpha}   &   V_{\beta}  
		    \end{bmatrix} ^{T},\, Z \in {{\mathbb B}^{o}}\subset {\mathbb R}^{|\hat{\beta}| \times |\beta|} \right\}.
		\end{equation*}
	 \end{lem}
	
	 {\it Proof} Note that $\Psi_n(X_0)=0$ and $\sigma_{{\mathbb B}^{o}}(U_{\hat{\beta}}^{T} H V_{\beta})$ is continuous in $H$. From \eqref{alpha_t 16}, we have 
	  \begin{equation}\label{subderivative of Psi}
		\Psi_n^{\prime}(X_0; H)=\lim_{\tau \downarrow 0 \atop H^{\prime} \to H}  \frac{\Psi_n(X_0+ \tau H^{\prime}) - \Psi_n(X_0)}{\tau} = \sigma_{{\mathbb B}^{o}} \big(U_{\hat{\beta}}^{T} H V_{\beta} \big).
	  \end{equation}
	  Then $\Psi_n$ is semidifferentiable at $X_0$ for $H$ and also epi-differentiable. Furthermore, we have $\Psi_n^{\prime}(X_0; H)= {\rm d}\Psi_n(X_0)(H)$ and
	  \begin{align*}
		\hat{\partial} \Psi_n(X_0) 
		=& \left\{ U_{\hat{\beta}} Z V_{\beta}^{T} \colon 
		   X_0 = \begin{bmatrix}
			       U_{\alpha}   &   U_{\hat{\beta}}  
	  	         \end{bmatrix} 
	             \Sigma(X_0)
	             \begin{bmatrix}
	             	V_{\alpha}   &   V_{\beta}  
	             \end{bmatrix}^{T}, \sigma_1(Z) \le 1 \right\},   \\
		=& \left\{ U_{\hat{\beta}} Z V_{\beta}^{T} \colon 
		   X_0 =\begin{bmatrix}
		   	      U_{\alpha}   &   U_{\hat{\beta}}  
		        \end{bmatrix} 
	            \Sigma(X_0) 
	            \begin{bmatrix}
	            	V_{\alpha}   &   V_{\beta}  
	            \end{bmatrix}^{T}, Z \in {{\mathbb B}^{o}}\subset {\mathbb R}^{ |\hat{\beta}|  \times|\beta| } \right\}.
	  \end{align*}
	  \qed

	 Next, we present the main conclusion of this section. From this, the second-order epi-derivative of the nuclear norm can be derived.
	 \begin{thm}\label{twice epi-differentiable of Psi_t}
		Let $X_0 \in {\mathbb M}_{m,n}$ of rank $r$ satisfy $n>r$ and $H \in {\mathbb M}_{m,n}$ be a small perturbation matrix. The function $\Psi_n$ is twice epi-differentiable at $X_0$. 
		The second-order epi-derivative of $\Psi_n$ at $X_0$ relative to any $\Omega \in \hat{\partial} \Psi_n(X_0)$ is given explicitly by 
		\begin{equation*}
			{\rm d}^2\Psi_n(X_0 \mid \Omega)(H) = \left\{
			\begin{array}{lr}
				-2\langle \Omega, H V_{\alpha} \Sigma_{\alpha}^{-1} U_{\alpha}^{T} H \rangle 
				& \mbox{if}\,\, \Psi_n^{\prime}(X_0; H) = \langle \Omega,H \rangle,   \\
				+\infty                                                                                                                                           & \mbox{if}\,\, \Psi_n^{\prime}(X_0; H) > \langle \Omega,H\rangle.     \\
			\end{array}
			\right.
		\end{equation*}
	 \end{thm}
 
	 {\it Proof} Let $\Omega \in \hat{\partial} \Psi_n(X_0)$. Given $X_0, H \in {\mathbb M}_{m,n}$, $H_k\rightarrow H$, and $t_k \rightarrow 0^{+}$, we focus our attention on the second-order difference quotient
	  \begin{equation}\label{second-order DQ of alpha_t}
		\frac{\Psi_n(X_0 + t_k H_k)- \Psi_n(X_0)- t_k \langle\Omega, H_k \rangle } {\frac{1}{2}t_k^2}.
	  \end{equation}
	  To begin with, assume that $\Psi_n^{\prime}(X_0; H) > \langle \Omega, H\rangle$. By using the relation \eqref{alpha_t 16} with $H=t_k H_k$, \eqref{second-order DQ of alpha_t} becomes
	  \begin{equation*}
		\frac{\Psi_n^{\prime}(X_0; H_k)- \langle\Omega, H_k \rangle } {\frac{1}{2}t_k}+ O(\|{\cal B}(H_k)\|^2),
	  \end{equation*}
	  which goes to $+\infty$ as $k \rightarrow +\infty$.

	  Assume now that $\Psi_n^{\prime}(X_0; H) = \langle \Omega, H\rangle$. Let $H_k\rightarrow H$ and $t_k \rightarrow 0^{+}$. Since $\Omega \in \hat{\partial} \Psi_n(X_0)$, it is possible to write $\Omega = U_{\hat{\beta} }Z V_{\beta}^{T}$ for some $Z \in {{\mathbb B}^{o}}$. By using the development \eqref{alpha_t 15}, the difference quotient \eqref{second-order DQ of alpha_t} can be written as
	  \begin{equation}\label{second-order DQ of alpha_t 18}
		\frac{\sigma_{{\mathbb B}^{o}} \left(U_{\hat{\beta}}^{T} H_k V_{\beta} - t_k U_{\hat{\beta}}^{T} H_k V_{\alpha} \Sigma_{\alpha}^{-1}(X_0) U_{\alpha}^{T} H_k V_{\beta} \right)- \left\langle Z, U_{\hat{\beta}}^{T} H_k V_{\beta} \right\rangle } {\frac{1}{2}t_k}+ O(t_k).
	  \end{equation}
	  Moreover, since $Z \in {{\mathbb B}^{o}}$, it is possible to bound from below the previous quantity by
	  \begin{align*}
		& \frac{ \left \langle Z, U_{\hat{\beta}}^{T} H_k V_{\beta} - t_k U_{\hat{\beta}}^{T} H_k V_{\alpha} \Sigma_{\alpha}^{-1}(X_0) U_{\alpha}^{T} H_k  V_{\beta} \right\rangle - \left\langle Z, U_{\hat{\beta}}^{T} H_k V_{\beta} \right\rangle } {\frac{1}{2}t_k}+ O(t_k)   \\
		=& -2 \left\langle Z,  U_{\hat{\beta}}^{T} H_k V_{\alpha} \Sigma_{\alpha}^{-1}(X_0) U_{\alpha}^{T} H_k  V_{\beta} \right\rangle + O(t_k) \\
		=& -2 \left\langle U_{\hat{\beta}} Z V_{\beta}^{T}, H_k V_{\alpha} \Sigma_{\alpha}^{-1}(X_0) U_{\alpha}^{T} H_k \right\rangle + O(t_k)
	  \end{align*}
	  which converges to $-2 \langle \Omega, H V_{\alpha} \Sigma_{\alpha}^{-1}(X_0) U_{\alpha}^{T} H \rangle$.

	  Let $t_k \rightarrow 0^{+}$. We have to exhibit a sequence $\{H_k\}_k$ which converges to $H$ such that the limit of \eqref{second-order DQ of alpha_t} equals $-2 \langle \Omega, H V_{\alpha} \Sigma_{\alpha}^{-1}(X_0) U_{\alpha}^{T} H \rangle$. 
	  Let us consider $H_k= H+ t_k H_k V_{\alpha} \Sigma_{\alpha}^{-1}(X_0) U_{\alpha}^{T} H_k$. 
	  The difference quotient \eqref{second-order DQ of alpha_t}, which is equal to \eqref{second-order DQ of alpha_t 18}, therefore becomes
	  \begin{align*}
		&\frac{\sigma_{{\mathbb B}^{o}} \big(U_{\hat{\beta}}^{T} H V_{\beta} \big)- \langle \Omega , H \rangle -t_k \langle \Omega, H_k V_{\alpha} \Sigma_{\alpha}^{-1}(X_0) U_{\alpha}^{T} H_k \rangle} {\frac{1}{2}t_k}+ O(t_k)  \\
		= &\frac{\Psi_n^{\prime}(X_0; H) -\langle \Omega , H \rangle}{\frac{1}{2}t_k} -2 \langle \Omega, H_k V_{\alpha} \Sigma_{\alpha}^{-1}(X_0) U_{\alpha}^{T} H_k \rangle+ O(t_k),
	  \end{align*}
	  which converges to $-2 \langle \Omega, H V_{\alpha} \Sigma_{\alpha}^{-1}(X_0) U_{\alpha}^{T} H \rangle$, as required.
	  \qed

	 This, combined with Lemma \ref{second-order differential of Phi_r} and \cite[Proposition 2.10]{Rockafellar 1988}, implies the following corollary.
	 \begin{cor}\label{twice epi-differentiable of nuclear norm}
		Let $X_0 \in {\mathbb M}_{m,n}$ of rank $r$ satisfy $n>r$ and $H \in {\mathbb M}_{m,n}$ be a small perturbation matrix. The nuclear norm is twice epi-differentiable at $X_0$. 
		The second-order epi-derivative of the nuclear norm at $X_0$ relative to any $\Omega \in \hat{\partial} \|X_0\|_*$ is given explicitly by 
		\begin{align*}
			{\rm d}^2 \|\cdot\|_*(X_0 \mid \Omega)(H) = &\sum_{k_l \leq r} {\rm tr} \Big(2 P_{k_l}^{T} {\cal B}(H) P_{k_l}^c (\sigma_{k_l}(X_0)I-\varLambda_{k_l})^{-1} {P_{k_l}^c}^{T} {\cal B}(H)P_{k_l} \Big) \\
			 &+ \delta_{ K_{\Psi_n}(X_0,\Omega)}(H)-2\langle \Omega, H V_{\alpha} \Sigma_{\alpha}^{-1} U_{\alpha}^{T} H \rangle,
		\end{align*}
	    where $K_{\Psi_n}(X_0,\Omega)$ is the critical cone of $\Psi_n$ at $X_0$ for $\Omega$.
	 \end{cor}

	\section{Subderivatives of Orthogonally Invariant Matrix Functions}
	 In this section, we present a chain rule for the subderivative of orthogonally invariant matrix functions in Theorem \ref{Subderivatives of Orthogonally Invariant Matrix Functions}. Similarly, we also give that of absolutely symmetric function in Theorem \ref{Subderivatives of Absolutely Symmetric Functions}. These important results are central to our developments in this paper. 
	 Recall that $Q_{\pm}$ is a signed permutation matrix if all its components are either $0$ or $\pm 1$ and each row and each column has exactly one nonzero element. Let $\mathbf{P}^n_{\pm}$ denote the set of all $n \times n$ signed permutation matrices. 
	 Recall also that a function $f \colon \R^n \rightarrow [-\infty,+\infty]$ is called absolutely symmetric if for every $x \in \R^n$ and every $n \times n$ signed permutation matrix $Q_{\pm}$, we have $f(Q_{\pm} x)=f(x)$. 
	 It is well-known (cf. Lewis \cite[Proposition 5.1]{LSI}) that for any orthogonally invariant matrix function $F \colon {\mathbb M}_{m,n} \rightarrow [-\infty,+\infty]$, there exists an absolutely symmetric function $f \colon \R^n \rightarrow [-\infty,+\infty]$ satisfying
	 \begin{equation}\label{composite form of orthogonally invariant matrix function F}
	 	F = f \circ \sigma.
	 \end{equation} 
     Indeed, $f$ can be chosen as a composite function of $F$ and the linear mapping $x \mapsto {\rm diag}(x)$ with $x \in \R^n$, namely,
	 \begin{equation}\label{composite form of absolutely symmetric function}
	 	f =F \circ {\rm diag}.
	 \end{equation}
     A set $\Gamma \subset {\mathbb M}_{m,n}$ is called an orthogonally invariant matrix set if $\delta_\Gamma$ is an orthogonally invariant matrix function. Similarly, $\Delta \subset \R^n$ is called an  absolutely symmetric set if $\delta_\Delta$ is an absolutely symmetric function. Therefore, it is easy to see that for any orthogonally invariant matrix set $\Gamma \subset {\mathbb M}_{m,n}$, there exists an absolutely symmetric set $\Delta \subset \R^n$ such that 
     \begin{equation}\label{orthogonally invariant matrix set Gamma}
     	\Gamma = \{X \in {\mathbb M}_{m,n} \colon \sigma(X) \in \Delta\},
     \end{equation} 
     where $\Delta$ can be chosen as 
     \begin{equation}\label{absolutely symmetric set Delta}
     	\Delta=\{x \in \R^n \colon {\rm diag}(x) \in \Gamma\}.     	
     \end{equation}
	 Next, we are going to justify a similar result as \cite[Proposition 2.3]{DLMS} for orthogonally invariant matrix sets, which allows us to obtain a chain rule for the subderivative of orthogonally invariant matrix functions via the established theory for composite functions in \cite[Theorem 3.4]{MMS}.

	 \begin{prop}\label{MSC of orthogonally invariant matrix function}
	  Let $K$ be an absolutely symmetric subset of $\R^n$. Then the distance function ${\rm d}_{\sigma^{-1}(K)}$ to the orthogonally invariant matrix set $\sigma^{-1}(K)$ satisfies:
	   \begin{equation}\label{dist to SVS}
	 		{\rm d}_{\sigma^{-1}(K)} = {\rm d}_K \circ \sigma.
	   \end{equation}
	  In particular, if $F$ is an orthogonally invariant matrix function, then for any $X \in {\mathbb M}_{m,n}$, we have
	   \begin{equation}\label{distance function to domF}
	  	    {\rm dist}(X, {\rm dom}F) = {\rm dist}(\sigma(X), {\rm dom}f),
	   \end{equation}
	  where $f$ is taken by \eqref{composite form of orthogonally invariant matrix function F}.
	 \end{prop}
	 
	 {\it Proof} To see that ${\rm d}_{\sigma^{-1}(K)}$ is an orthogonally invariant matrix function, we fix $X \in {\mathbb M}_{m,n}$ and $(U,V) \in  {\cal O}^{m,n}$ such that $X = U {\rm diag}(\sigma(X)) V^{T}$. Then we have
	  \begin{align*}
    	{\rm d}_K(\sigma(X))  &=     \inf \limits_{y \in K} \| \sigma(X) - y\|          \\
  	                          &=     \inf \limits_{y \in K} \| {\rm diag}(\sigma(X)) - {\rm diag}(y)\|     \\
  	                          &\geq  \inf \limits_{Y \in {\sigma^{-1}(K)}} \| {\rm diag}(\sigma(X)) -Y\|             \\ 
  	                          &=     \inf \limits_{Y \in {\sigma^{-1}(K)}} \|U {\rm diag}(\sigma(X)) V^{T} - U Y V^{T} \|     \\
  	                          &=     \inf \limits_{U^{T} \bar Y V \in {\sigma^{-1}(K)}} \|X- \bar Y\|     \\
  	                          &=     \inf \limits_{\bar Y  \in {\sigma^{-1}(K)}} \|X- \bar Y\|    \\
  	                          &=     {\rm d}_{\sigma^{-1}(K)}(X).
      \end{align*}
      On the other hand, it follows from \eqref{Von Neumann’s Trace inequality} that 
	  \begin{align*}
	 	{\rm d}_{\sigma^{-1}(K)}(X) &=    \inf \limits_{Y \in {\sigma^{-1}(K)}} \|X-Y\|  \\
	 	                            &\geq \inf \limits_{Y \in {\sigma^{-1}(K)}} \| \sigma(X) - \sigma(Y) \|  \\
	 	                            &=    \inf \limits_{y \in K} \| \sigma(X) - y\|  \\
	 	                            &=    {\rm d}_K(\sigma(X)).    	                               
	  \end{align*}
	  In particular, if $F$ is an orthogonally invariant matrix function, there exists an absolutely symmetric function $f \colon \R^n \rightarrow [-\infty,+\infty]$ satisfying \eqref{composite form of orthogonally invariant matrix function F}. Combining this with \eqref{composite form of orthogonally invariant matrix function F}, we obtain that
	  \begin{equation} \label{def of domF}
	    {\rm dom} F = \{ X \in {\mathbb M}_{m,n} \colon \sigma(X) \in {\rm dom} f \} = \sigma^{-1}({\rm dom} f)
	  \end{equation}
	  is an orthogonally invariant matrix set and ${\rm dom} f$ is an absolutely symmetric set. From \eqref{dist to SVS}, we can deduce that ${\rm dist}(X, {\rm dom}F) = {\rm dist}(\sigma(X), {\rm dom}f)$ for any $X \in {\mathbb M}_{m,n}$.
	  \qed
	
	 To prove a chain rule for subderivatives of orthogonally invariant matrix functions, let's first recall a useful characterization of the subdifferential of orthogonally invariant matrix functions.
	 
	 \begin{prop}\label{subdifferential of orthogonally invariant matrix functions}
	  Assume that $f \colon \R^n \rightarrow [-\infty,+\infty]$ is a proper, convex, lsc, and absolutely symmetric function. Then the following properties are equivalent: 
	  
	  {\rm (i)}
	  $Y \in \partial(f \circ \sigma)(X)$;   
	  
	  {\rm (ii)} 
	  $\sigma(Y) \in \partial f(\sigma(X))$ and the matrices $X$ and $Y$ have a simultaneous ordered singular value decomposition.
	 \end{prop}

	 {\it Proof}  According to Lewis \cite[Theorem 2.4, Corollary 2.6]{Lewis 1995}, $f \circ \sigma$ is convex and lsc if and only if $f$ is convex and lsc. The claimed equivalence then follows from Lewis \cite[Corollary 2.5]{Lewis 1995}.
	 \qed

	 \begin{thm}\label{Subderivatives of Orthogonally Invariant Matrix Functions} (Subderivatives of Orthogonally Invariant Matrix Functions).
	  Let $f \colon \R^n \rightarrow [-\infty,+\infty]$ be an absolutely symmetric function and let $X \in {\mathbb M}_{m,n}$ with $(f \circ \sigma)(X)$ finite. If $f$ satisfies one of the following conditions:\\
	   {\rm (a)}
	   $f$ is  lsc and convex with $\partial f(\sigma(X)) \neq \emptyset$;\\
	   {\rm (b)} 
	   $f$ is locally Lipschitz continuous around $\sigma(X)$ relative to its domain.\\
	  Then for all $H \in {\mathbb M}_{m,n}$, we have
	  \begin{equation}\label{a chain rule for subderivatives of OIMF}
	  	{\rm d}(f \circ \sigma)(X)(H) = {\rm d}f(\sigma(X))(\sigma^{\prime}(X;H)). 
	  \end{equation}
	 \end{thm}
 
     {\it Proof}  Choose any $H \in {\mathbb M}_{m,n}$. According to Proposition \ref{fdd of SV M}, $\sigma^{\prime}(X;\cdot)$ is a Lipschitz-continuous and positively homogeneous function. Moreover, $\sigma^{\prime}(X;E)+O(t^2\|E\|^2)/t \rightarrow \sigma^{\prime}(X;H)$ as $t \downarrow 0$ and $E \rightarrow H$. Combining this with the definition of subderivative, we derive
     \begin{align}\label{geq in subderivatives of OIMF}
     	{\rm d}(f \circ \sigma)(X)(H) &=    \liminf_{ {t \downarrow 0} \atop{ E \to H}} \frac{ f(\sigma(X+tE)) - f(\sigma(X)) }{t}  \nonumber  \\
     	                              &=    \liminf_{ {t \downarrow 0} \atop{ E \to H}} \frac{f(\sigma(X)+t\sigma^{\prime}(X;E)+O(t^2 \|E\|^2))-f(\sigma(X))}{t}  \nonumber  \\
     	                              &=     \liminf_{ {t \downarrow 0} \atop{ E \to H}} \frac{f(\sigma(X)+t(\sigma^{\prime}(X;E)+\frac{O(t^2 \|E\|^2)}{t}))-f(\sigma(X))}{t}   \nonumber  \\
     	                              &\geq  {\rm d}f(\sigma(X))(\sigma^{\prime}(X;H)).
     \end{align}
     This establishes the inequality "$\geq$" in \eqref{a chain rule for subderivatives of OIMF}. 
     For the opposite inequality, let's consider two cases. First, if ${\rm d}f(\sigma(X))(\sigma^{\prime}(X;H)) = \infty$, the latter inequality clearly holds. So, we'll assume ${\rm d}f(\sigma(X))(\sigma^{\prime}(X;H)) < \infty$. 
     If $f$ is convex and lsc, $f \circ \sigma$ is also convex and lsc according to \cite[Corollary 2.6]{Lewis 1995}. Furthermore, Lewis and Sendov \cite[Theorem 7.1]{LSI} and the non-emptiness of $\partial f(\sigma(X))$ imply that $\partial(f \circ \sigma)(X) \neq \emptyset$. 
     Thus, it follows from Bonnans and Shapiro \cite[Proposition 2.126]{Bonnans Shapiro} that 
      ${\rm d}(f \circ \sigma)(X)(H) = \sup\limits_{Y \in \partial(f \circ \sigma)(X)}\langle Y,H \rangle$.
     Let $\epsilon >0$ and choose $Y \in \partial(f \circ \sigma)(X)$ such that ${\rm d}(f \circ \sigma)(X)(H) \leq \langle Y,H \rangle + \epsilon$. 
     Since $Y \in \partial(f \circ \sigma)(X)$, it follows from Proposition \ref{subdifferential of orthogonally invariant matrix functions} that there exists $(U,V) \in {\cal O}^{m,n}$ such that $\sigma(Y) \in \partial f(\sigma(X))$. 
     Setting $\Sigma(Y):= U^T Y V$ and applying Proposition \ref{fdd of SV M} and the fact that 
     $\left\langle \Sigma(Y)_{\alpha_i \alpha_i}, \frac{1}{2}(V_{\alpha_i}^T H^T U_{\alpha_i} - U_{\alpha_i}^T H V_{\alpha_i}) \right\rangle = 0$ for any $i=1, \ldots, t$, we get
     \begin{align*}
    	  & {\rm d}(f \circ \sigma)(X)(H)  \leq \left\langle Y,H \right\rangle + \epsilon = \left\langle \Sigma(Y), U^T H V \right\rangle + \epsilon \\
       =  & \sum_{i=1}^{t} \left\langle \Sigma(Y)_{\alpha_i \alpha_i}, U_{\alpha_i}^T H V_{\alpha_i} \right\rangle + \left\langle \Sigma(Y)_{\hat{\beta} \beta}, U_{\hat{\beta}}^T H V_{\beta} \right\rangle + \epsilon \\
       =  & \sum_{i=1}^{t} \left\langle \Sigma(Y)_{\alpha_i \alpha_i}, U_{\alpha_i}^T H V_{\alpha_i} \right\rangle + \sum_{i=1}^{t} \left\langle \Sigma(Y)_{\alpha_i \alpha_i}, \frac{1}{2} \big(V_{\alpha_i}^T H^T U_{\alpha_i} - U_{\alpha_i}^T H V_{\alpha_i} \big) \right\rangle + \left\langle \Sigma(Y)_{\hat{\beta} \beta}, U_{\hat{\beta}}^T H V_{\beta} \right\rangle + \epsilon \\
       =  & \sum_{i=1}^{t} \left\langle \Sigma(Y)_{\alpha_i \alpha_i}, \frac{1}{2} \big( U_{\alpha_i}^T H V_{\alpha_i} + V_{\alpha_i}^T H^T U_{\alpha_i} \big) \right\rangle + \left\langle \Sigma(Y)_{\hat{\beta} \beta}, U_{\hat{\beta}}^T H V_{\beta} \right\rangle + \epsilon \\
     \leq & \sum_{i=1}^{t} \left\langle \lambda \big( \Sigma(Y)_{\alpha_i \alpha_i} \big), \frac{1}{2} \lambda \big(U_{\alpha_i}^T H V_{\alpha_i} + V_{\alpha_i}^T H^T U_{\alpha_i} \big) \right\rangle + \left\langle \sigma(\Sigma(Y)_{\hat{\beta} \beta}), \sigma(U_{\hat{\beta}}^T H V_{\beta}) \right\rangle + \epsilon \\
       =  & \sum_{i=1}^{t}  \sum_{s \in\alpha_i} \left\langle \sigma_s(Y),\frac{1}{2} \lambda_{l_s} \big(U_{\alpha_i}^T H V_{\alpha_i} + V_{\alpha_i}^T H^T U_{\alpha_i} \big) \right\rangle + \sum_{s \in \beta} \left\langle \sigma_s(Y), \sigma_{l_s} \big(U_{\hat{\beta}}^T H V_{\beta} \big) \right\rangle + \epsilon \\
       =  & \sum_{i=1}^{t}  \sum_{s \in\alpha_i} \left\langle \sigma_s(Y),\sigma_s^{\prime}(X;H) \right\rangle + \sum_{s \in \beta} \left\langle \sigma_s(Y), \sigma_s^{\prime}(X;H) \right\rangle + \epsilon \\
       =  & \left\langle \sigma(Y), \sigma^{\prime}(X;H) \right\rangle + \epsilon \\
      \leq & {\rm d}f(\sigma(X))(\sigma^{\prime}(X;H)) + \epsilon.
     \end{align*}
     This second inequality follows from Fan’s inequality \eqref{Fan’s inequality} and von Neumann’s trace inequality \eqref{Von Neumann’s Trace Inequality} and the last inequality stems from $f$ being lsc and convex and  $\sigma(Y) \in \partial f(\sigma(X))$. 
     As $\epsilon \downarrow 0$, we obtain the opposite inequality in \eqref{geq in subderivatives of OIMF}, thus proving \eqref{a chain rule for subderivatives of OIMF} in this case. 
     Assume that $f$ is locally Lipschitz continuous around $\sigma(X)$ relative to its domain. To prove the opposite inequality in \eqref{geq in subderivatives of OIMF},  we need only consider the case where ${\rm d}f(\sigma(X))(\sigma^{\prime}(X;H)) < \infty$. By the definition of subderivatives, there exist sequences $t_k \downarrow 0$ and $v_k \to \sigma^{\prime}(X;H)$ such that
     \begin{equation}\label{subderi of f non inf}
     	{\rm d}f(\sigma(X))(\sigma^{\prime}(X;H)) = \lim\limits_{k \to \infty} \frac{f(\sigma(X)+t_kv_k)-f(\sigma(X))}{t_k}.
     \end{equation}
     Without loss of generality, we assume that $\sigma(X)+t_kv_k \in {\rm dom}f$ for all $k \in {\mathbb N}$. From \eqref{distance function to domF}, we get
     $$
       {\rm dist}(X+t_k H, {\rm dom}F) = {\rm dist}(\sigma(X+t_k H), {\rm dom}f) \,\, \mbox{for all}\,\, k \in {\mathbb N},
     $$
     which, together with Proposition \ref{fdd of SV}, implies that for all  $k \in {\mathbb N}$,
     \begin{align*}
     	{\rm dist} \bigg( H, \frac{{\rm dom}F - X}{t_k} \bigg) &= \frac{1}{t_k}{\rm dist}(\sigma(X)+ t_k \sigma^{\prime}(X;H)+O(t_k^2), {\rm dom}f) \\
     	&= \frac{1}{t_k} \|\sigma(X)+ t_k \sigma^{\prime}(X;H)+O(t_k^2)- (\sigma(X)+t_kv_k)\|  \\
     	&= \bigg\| \sigma^{\prime}(X;H)- v_k+\frac{O(t_k^2)}{t_k} \bigg\|.
     \end{align*}
     Then for each $k \in {\mathbb N}$, we can find a matrix $H_k \in  \frac{{\rm dom}F - X}{t_k}$ satisfying 
     $$
     \| H-H_k\| < \bigg\| \sigma^{\prime}(X;H)- v_k+\frac{O(t_k^2)}{t_k} \bigg\| + \frac{1}{k}.
     $$
     This allows us to have $X+t_k H_k \in {\rm dom}F$ for all $k$ and $H_k \to H$ as $k \to \infty$. 
     This, coupled with \eqref{subderi of f non inf}, Proposition \ref{fdd of SV} and the imposed 
     assumption on $f$, leads us to obtain
     \begin{align*}
     	{\rm d}f(\sigma(X))(\sigma^{\prime}(X;H)) &= \lim\limits_{k \to \infty} \bigg[ \frac{f(\sigma( X+t_k H_k) ) - f(\sigma( X) )}{t_k} + \frac{f(\sigma(X)+t_kv_k) - f(\sigma( X+t_k H_k) )}{t_k} \bigg] \\
     	                                          &\geq \liminf_{k \to \infty} \frac{f(\sigma( X+t_k H_k) ) - f(\sigma( X) )}{t_k} - l \lim\limits_{k \to \infty} \bigg\|  \frac{\sigma( X+t_k H_k ) - \sigma(X)}{t_k} -v_k \bigg\| \\
     	                                          &\geq  {\rm d}(f \circ \sigma)(X)(H) - l \lim\limits_{k \to \infty} \bigg\|  \sigma^{\prime}( X;H_k)+ \frac{O(t_k^2)}{t_k} -v_k \bigg\| \\
     	                                          &=  {\rm d}(f \circ \sigma)(X)(H),
     	\end{align*}
     where $l \geq 0$ is a Lipschitz constant of $f$ around $\sigma(X)$ relative to its domain.  This verifies the inequality "$\leq$" in \eqref{a chain rule for subderivatives of OIMF}, and ends the proof.
     \qed
    
     As a direct result of Theorem \ref{Subderivatives of Orthogonally Invariant Matrix Functions}, we derive a straightforward representation of tangent cones for orthogonally invariant matrix sets. Note that reducing \eqref{distance function to domF}, which was stated for the orthogonally invariant matrix function in \eqref{composite form of orthogonally invariant matrix function F}, to the orthogonally invariant matrix set $\Gamma$ in \eqref{orthogonally invariant matrix set Gamma} leads us the estimate
     \begin{equation}\label{distance function to Gamma}
     	{\rm dist}(X, \Gamma) = {\rm dist}(\sigma(X), \Delta) \,\, \mbox{for any}\,\, X \in {\mathbb M}_{m,n},
     \end{equation}
     where $\Delta$ is taken from \eqref{absolutely symmetric set Delta}.
     
     \begin{cor}\label{Tangent Cone to the Orthogonally Invariant Matrix Sets} (Tangent Cone to the Orthogonally Invariant Matrix Sets).
       Let $\Gamma$ be an orthogonally invariant matrix set represented by \eqref{orthogonally invariant matrix set Gamma}. Then for any $X \in \Gamma$, we have 
      	\begin{equation*}
      		T_{\Gamma}(X) = \{ H \in {\mathbb M}_{m,n} \colon \sigma^{\prime}(X;H) \in T_{\Delta} (\sigma(X))\}.
      	\end{equation*}
     \end{cor}
 
     {\it Proof} Using the absolutely symmetric set $\Delta$ from \eqref{orthogonally invariant matrix set Gamma}, we apply Theorem \ref{Subderivatives of Orthogonally Invariant Matrix Functions} to the absolutely symmetric function $\delta_{\Delta}$ and then 
     $$
     {\rm d}\delta_{\Gamma}(X)(H) = {\rm d}\delta_{\Delta}(\sigma(X))(\sigma^{\prime}(X;H)).
     $$
     The representation of the tangent cone to $\Gamma$ at $X \in \Gamma$ follows from two key facts: ${\rm d}\delta_{\Gamma}(X) = \delta_{T_{\Gamma}(X)}$ and ${\rm d}\delta_{\Delta}(\sigma(X)) = \delta_{T_{\Delta}(\sigma(X))}$.
     \qed
	
	 In what follows, we are going to present the second-order tangent sets of orthogonally invariant matrix sets. To this end, we begin by demonstrating that certain second-order approximations of orthogonally invariant matrix sets possess an outer Lipschitzian property. Define the set-valued mapping $G_H \colon \R^n \rightrightarrows {\mathbb M}_{m,n}$ via the second-order tangent set to the absolutely symmetric set $\Delta$ in \eqref{absolutely symmetric set Delta} by
	 \begin{equation}\label{def of G_H}
	 	G_H(b) = \{ W \in {\mathbb M}_{m,n} \colon \sigma^{\prime \prime}(X;H,W) + b \in T_{\Delta}^2 (\sigma(X),\sigma^{\prime}(X;H))\}.
	 \end{equation}

	 \begin{prop}
	 	Assume that $\Gamma$ is an orthogonally invariant matrix set represented by \eqref{orthogonally invariant matrix set Gamma} and that $X \in \Gamma$ and $H \in T_{\Gamma}(X)$. Then, the mapping $G_H$ in \eqref{def of G_H} enjoys the following uniform outer Lipschitzian property at the origin:
        \begin{equation}\label{outer Lipschitzian}
        	G_H(b) \subset G_H(0)+\|b\| {\mathbb B} \,\, \mbox{for any}\,\, b \in \R^n.
        \end{equation}
	 \end{prop}
 
	 {\it Proof} Let $b \in \R^n$ and pick any $W \in G_H(b)$. It follows from \eqref{def of G_H} that $\sigma^{\prime \prime}(X;H,W) + b \in T_{\Delta}^2 (\sigma(X),\sigma^{\prime}(X;H))$. By the definition of second-order tangent set, there exists a sequence $t_k \downarrow 0$ such that
	 $$
	 \sigma(X) + t_k \sigma^{\prime}(X;H) + \frac{1}{2}t_k^2\sigma^{\prime \prime}(X;H,W) + \frac{1}{2}t_k^2 b +o(t_k^2) \in \Delta \,\,\mbox{ for any} \,\, k \in {\mathbb N}.
	 $$
	 For any $k$ sufficiently large, we conclude from \eqref{estimate for the SVF} that
	 \begin{equation*}
	 	\sigma \left(X + t_k H+\frac{1}{2}t_k^2 W \right) = \sigma(X) + t_k \sigma^{\prime}(X;H) + \frac{1}{2}t_k^2\sigma^{\prime \prime}(X;H,W) +o(t_k^2),
	 \end{equation*}
	 which, coupled with \eqref{distance function to Gamma}, leads us to
	 $$
	 {\rm dist}\left(X + t_k H+\frac{1}{2}t_k^2 W, \Gamma \right) = {\rm dist}\left(\sigma(X + t_k H+\frac{1}{2}t_k^2 W), \Delta \right) \leq \frac{1}{2}t_k^2 \|b\| +o(t_k^2).
	 $$
	 Thus, there exists $Y_k \in \Gamma$ such that 
	 $$\|E_k\| \leq \frac{1}{2} \|b\| + \frac{o(t_k^2)}{t_k^2},$$
	 where $E_k := \frac{X + t_k H+\frac{1}{2}t_k^2 W - Y_k}{t_k^2}$. 
	 Passing to a subsequence, if necessary, implies that there exists $E \in {\mathbb M}_{m,n}$ such that $E_k \to E$ as $k \to \infty$. This brings us to 
	 \begin{equation}\label{2E leq b}
	 	\|E\| \leq \frac{1}{2} \|b\|.
	 \end{equation}
	 Since $X + t_k H+\frac{1}{2}t_k^2 W - t_k^2 E_k = Y_k \in \Gamma$, it follows from \eqref{orthogonally invariant matrix set Gamma} and \eqref{estimate for the SVF} that 
	 $$
	 \sigma(X) + t_k\sigma^{\prime}(X;H) + \frac{1}{2}t_k^2 \sigma^{\prime\prime}(X;H,W-2E) +o(t_k^2) = \sigma(X + t_k H+\frac{1}{2}t_k^2 W - t_k^2 E_k ) \in \Delta.
	 $$ 
	  By the definition of the second-order tangent set, we get
	  $$
	  \sigma^{\prime\prime}(X;H,W-2E) \in T_{\Delta}^2 (\sigma(X), \sigma^{\prime}(X;H)),
	  $$
	 which leads us to $W-2E \in G_H(0)$. This, combined with \eqref{2E leq b}, implies \eqref{outer Lipschitzian} and thus ends the proof. 
	 \qed

	 \begin{prop}\label{prop of SOTS}(Second-Order Tangent Sets of Orthogonally Invariant Matrix Sets).
	  Assume that $\Gamma$ is an orthogonally invariant matrix set represented by \eqref{orthogonally invariant matrix set Gamma} and that $X \in \Gamma$ and $H \in T_{\Gamma}(X)$. Then, we have 
		\begin{equation}\label{formula of SOTS}
			T_{\Gamma}^2 (X,H) = \{  W \in {\mathbb M}_{m,n} \colon \sigma^{\prime \prime}(X;H,W) \in T_{\Delta}^2 (\sigma(X),\sigma^{\prime}(X;H))\},
		\end{equation}
	  where $\Delta$ is taken from \eqref{absolutely symmetric set Delta}. Moreover, if the absolutely symmetric set $\Delta$ is parabolically derivable at $\sigma(X)$ for $\sigma^{\prime}(X;H)$, then $\Gamma$ is parabolically derivable at $X$ for $H$.
	 \end{prop}

     {\it Proof} Pick any $W \in T_{\Gamma}^2 (X,H)$. By the definition of second-order tangent set, there exists a sequence $t_k \downarrow 0$ such that
     $$
      X + t_kH + \frac{1}{2}t_k^2 W +o(t_k^2) \in \Gamma.
     $$
     It follows from \eqref{orthogonally invariant matrix set Gamma} and \eqref{estimate for the SVF} that 
     $$
     \sigma(X) + t_k\sigma^{\prime}(X;H) + \frac{1}{2}t_k^2 \sigma^{\prime\prime}(X;H,W) +o(t_k^2) = \sigma(X + t_k H+\frac{1}{2}t_k^2 W +o(t_k^2)) \in \Delta.
     $$ 
     Thus we get $\sigma^{\prime \prime}(X;H,W) \in T_{\Delta}^2 (\sigma(X),\sigma^{\prime}(X;H))$  and prove the inclusion ``$\subset$'' in \eqref{formula of SOTS}. 
     To prove the opposite inclusion in \eqref{formula of SOTS}, we take any $W \in {\mathbb M}_{m,n}$ satisfying $\sigma^{\prime \prime}(X;H,W) \in T_{\Delta}^2 (\sigma(X),\sigma^{\prime}(X;H))$. By the definition of second-order tangent set, there exists a sequence $t_k \downarrow 0$ such that
     $$
       \sigma(X)+ t_k\sigma^{\prime}(X;H) + \frac{1}{2}t_k^2 \sigma^{\prime \prime}(X;H,W)  +o(t_k^2) \in \Delta.
     $$
     Combining \eqref{estimate for the SVF} with \eqref{distance function to Gamma}, we obtain that for any $t>0$ sufficiently small, 
     \begin{align*}
     	{\rm dist}\left(X + t H+\frac{1}{2}t^2 W, \Gamma \right) &= {\rm dist}\left(\sigma \left(X + t H+\frac{1}{2}t^2 W \right), \Delta \right)\\
     	                                                         &= {\rm dist}\left(\sigma(X) + t \sigma^{\prime}(X;H) + \frac{1}{2}t^2\sigma^{\prime \prime}(X;H,W) +o(t^2), \Delta \right).
     \end{align*}
     Thus we get that $X + t_k H+\frac{1}{2}t_k^2 W \in \Gamma$, hence that $W \in T_{\Gamma}^2 (X,H)$. This verifies the inclusion "$\supset$" in \eqref{formula of SOTS}, and then \eqref{formula of SOTS} is proved.
     
     To prove the parabolic derivability of $\Gamma$ at $X$ for $H$, we first show that $T_{\Gamma}^2 (X,H)$ is a non-empty set. 
     Suppose that the absolutely symmetric set $\Delta$ is parabolically derivable at $\sigma(X)$ for $\sigma^{\prime}(X;H)$, it follows from the definition of parabolically derivable that $T_{\Delta}^2(\sigma(X), \sigma^{\prime}(X;H))$ is nonempty. 
     Thus there exists $p \in T_{\Delta}^2(\sigma(X), \sigma^{\prime}(X;H))$ and $\widetilde{W} \in {\mathbb M}_{m,n}$ such that 
     $$
      \sigma^{\prime \prime}(X;H, \widetilde{W}) + b \in T_{\Delta}^2 (\sigma(X),\sigma^{\prime}(X;H)),
     $$
     where $b:= p-\sigma^{\prime \prime}(X;H, \widetilde{W})$. It is easy to obtain $\widetilde{W} \in G_H(b)$ by \eqref{def of G_H}. According to the outer Lipschitzian property of the set-valued mapping $G_H$, there exists $W \in G_H(0)$ such that $\|W-\widetilde{W}\| \leq b$. This leads us to $\sigma^{\prime \prime}(X;H,W) \in T_{\Delta}^2 (\sigma(X),\sigma^{\prime}(X;H))$. By \eqref{formula of SOTS}, we have $W \in T_{\Gamma}^2 (X,H)$, which justifies $T_{\Gamma}^2 (X,H) \neq \emptyset$. 
     It remains to prove that for each $W \in T_{\Gamma}^2 (X,H)$, there exists $\epsilon >0$ and an arc $\xi \colon [0,\epsilon] \to \Gamma$ with $\xi(0)=X$, $\xi^{\prime}_+(0)=H$, and $\xi^{\prime \prime}_+(0)=W$.  It suffices to show that for all $t \in [0,\epsilon]$,
     \begin{equation}\label{pd of Gamma}
     	X + t H+\frac{1}{2}t^2 W +o(t^2)\in \Gamma.
     \end{equation}
     Pick any $W \in T_{\Gamma}^2 (X,H)$, we have $\sigma^{\prime \prime}(X;H,W) \in T_{\Delta}^2(\sigma(X),\sigma^{\prime}(X;H))$ by \eqref{formula of SOTS}. 
     Since the absolutely symmetric set $\Delta$ is parabolically derivable at $\sigma(X)$ for $\sigma^{\prime}(X;H)$, there exists $\epsilon^{\prime} >0$ and an arc $\zeta \colon [0,\epsilon^{\prime}] \to \Delta$ with $\zeta(0)=\sigma(X)$, $\zeta^{\prime}_+(0)=\sigma^{\prime}(X;H)$, and $\zeta^{\prime \prime}_+(0)=\sigma^{\prime \prime}(X;H,W)$. 
     Thus we get $\zeta(t)= \sigma(X)+ t\sigma^{\prime}(X;H) + \frac{1}{2}t^2 \sigma^{\prime \prime}(X;H,W) +o(t^2) \in \Delta$ for all $t \in [0,\epsilon^{\prime}]$. 
     This, combined with \eqref{estimate for the SVF} and \eqref{orthogonally invariant matrix set Gamma}, implies \eqref{pd of Gamma} for some $\epsilon \leq \epsilon^{\prime}$. This justifies that $\Gamma$ is parabolically derivable at $X$ for $H$, and completes the proof. 
     \qed

     We proceed by proving a chain rule for subderivatives of absolutely symmetric function. We begin with presenting a counterpart of the estimate in \eqref{distance function to domF} for domains of absolutely symmetric functions.
    
     \begin{prop}\label{MSC of absolutely symmetric function}
    	Let $F \colon {\mathbb M}_{m,n} \rightarrow [-\infty,+\infty]$ be an orthogonally invariant matrix function, represented by \eqref{composite form of orthogonally invariant matrix function F}. Then for any $x \in \R^n$, we have
    	\begin{equation}\label{dist to domf}
    		{\rm dist}(x, {\rm dom}f) = {\rm dist}({\rm diag} (x), {\rm dom}F),
    	\end{equation}
    	where $f$ is taken by \eqref{composite form of orthogonally invariant matrix function F}.
     \end{prop}
    
     {\it Proof} For any $x \in \R^n$, there exists a signed permutation matrix $Q_{\pm} \in \mathbf{P}^n_{\pm}$ such that $\sigma({\rm diag} (x)) = Q_{\pm} x$. As pointed out before, the absolutely symmetric function $f$ can be represented as a composite function of the orthogonally invariant matrix function $F$ and the linear mapping $x \mapsto {\rm diag}(x)$ with $x \in \R^n$. Thus, it follows from \eqref{composite form of absolutely symmetric function} that   
     \begin{equation}\label{def of domf}
     	{\rm dom} f = \{ x \in \R^n \colon {\rm diag} (x) \in {\rm dom} F \}
     \end{equation}
     is an absolutely symmetric set. 
     This, combined with \eqref{def of domF}, implies that $Q^T_{\pm} \sigma(X) \in {\rm dom}f$ for any $X \in {\rm dom}F$. 
     Then by von Neumann’s trace inequality \eqref{Von Neumann’s Trace Inequality}, we obtain
     \begin{align*}
    	{\rm dist}(x, {\rm dom}f) \leq \| x - Q^T_{\pm} \sigma(X) \|= \| Q_{\pm}x - \sigma(X) \| 
    	                          =    \| \sigma({\rm diag} (x)) - \sigma(X) \| \leq  \| {\rm diag} (x) - X \|,
     \end{align*}
     for all $X \in {\rm dom}F$, which, in turn, leads us to
     \begin{equation}\label{leq in dist to domf}
      {\rm dist}(x, {\rm dom}f) \leq {\rm dist}({\rm diag} (x), {\rm dom}F).
     \end{equation}
     To prove the opposite inequality, pick any $y \in {\rm dom}f$. By \eqref{def of domf}, we have ${\rm diag} (y) \in {\rm dom} F$ and then 
     $$
      {\rm dist}({\rm diag} (x), {\rm dom}F) \leq \|{\rm diag} (x)-{\rm diag} (y)  \| = \|x-y\|,
     $$
     which implies the opposite inequality in \eqref{leq in dist to domf}, and completes the proof.
     \qed
     
     Note that the linear mapping $x \mapsto {\rm diag}(x)$ with $x \in \R^n$ is twice continuously differentiable. This, coupled with \eqref{dist to domf} and \cite[Theorem 3.4]{MMS}, brings us to the following theorem.
    
     \begin{thm}\label{Subderivatives of Absolutely Symmetric Functions} (Subderivatives of Absolutely Symmetric Functions).
     	Let $F \colon {\mathbb M}_{m,n} \rightarrow [-\infty,+\infty]$ be an orthogonally invariant matrix function, represented by \eqref{composite form of orthogonally invariant matrix function F}, and let its associated absolutely symmetric function $f$ be locally Lipschitz continuous relative to its domain. Then, for any $X \in {\mathbb M}_{m,n}$ with $F(X)$ finite and any $z \in \R^n$, we have
    	\begin{equation}\label{a chain rule for subderivatives of ASF}
    		{\rm d}f(\sigma(X))(z) = {\rm d}F( {\rm diag}(\sigma(X)) )({\rm diag}(z)). 
    	\end{equation}
     \end{thm}
   
	 {\it Proof} Since $f$ is locally Lipschitz continuous relative to its domain, by \eqref{Von Neumann’s Trace inequality} we obtain that $F$ is locally Lipschitz continuous relative to its domain. According to \cite[Theorem 3.4]{MMS}, we  justify \eqref{a chain rule for subderivatives of ASF}.
	 \qed

	 Similar to Corollary \ref{Tangent Cone to the Orthogonally Invariant Matrix Sets}, we proceed to derive a refined representation of tangent cones for absolutely symmetric sets according to Theorem \ref{Subderivatives of Absolutely Symmetric Functions}.
	
	 \begin{cor}\label{Tangent Cone to the Absolutely Symmetric Sets} (Tangent Cone to the Absolutely Symmetric Sets).
		Let $\Delta$ be an absolutely symmetric set represented by \eqref{absolutely symmetric set Delta}. Then for any $X \in \Gamma$, we have 
		\begin{equation*}
			T_{\Delta}(\sigma(X)) = \{ z \in \R^n \colon {\rm diag}(z) \in T_{\Gamma} ( {\rm diag}(\sigma(X)) )\}.
		\end{equation*}
	 \end{cor}
   
     We close this section by revealing that the subderivative of $f$ at $\sigma(X)$ is a symmetric function with respect to $\mathbf{P}^n_{\pm}(X)$, which is a subset of $\mathbf{P}^n_{\pm}$ consisting of all $n \times n$ block diagonal matrices in the form $Q_{\pm} = {\rm diag}(Q_1, \ldots,Q_t,Q_{t+1})$, where $Q_i \in \R^{|\alpha_i| \times |\alpha_i|}$ is a signed permutation matrix for any $i=1, \ldots,t$ with $\alpha_i$ taken from \eqref{alpha_i}, $t$ being the number of distinct nonzero singular values of $X$ and $Q_{t+1} \in \R^{|\beta| \times |\beta|}$ is also a signed permutation matrix. Obviously, we have that $\mathbf{P}^n_{\pm}(X) \subset \mathbf{P}^n_{\pm}$ and that if $Q_{\pm} \in \mathbf{P}^n_{\pm}(X)$, then $Q_{\pm} \sigma(X) = \sigma(X)$.
    
     \begin{prop}\label{symmetric property of subderivative function}
    	Assume that $f \colon \R^n \rightarrow [-\infty,+\infty]$ is an absolutely symmetric function and $X \in {\mathbb M}_{m,n}$ with $f(\sigma(X))$ finite. Then, for any $v \in \R^n$ and any signed permutation matrix $Q_{\pm} \in \mathbf{P}^n_{\pm}(X)$, we have
    	\begin{equation*}
    		{\rm d}f(\sigma(X))(Q_{\pm} v) = {\rm d}f(\sigma(X))(v),
    	\end{equation*}
    	which means that the subderivative $v \mapsto {\rm d}f(\sigma(X))(v)$ is symmetric with respect to $\mathbf{P}^n_{\pm}(X)$.
     \end{prop}
    
     {\it Proof} For any $v \in \R^n$ and $Q_{\pm} \in \mathbf{P}^n_{\pm}(X)$, it follows from the absolutely symmetric property of $f$ that 
     \begin{align*}
     	{\rm d}f(\sigma(X))(v) &= \liminf_{ {\tau \downarrow 0} \atop{ v^{\prime} \to v}} \frac{ f(\sigma(X) + \tau  v^{\prime}) - f(\sigma(X))}{\tau} \\
     	                       &= \liminf_{ {\tau \downarrow 0} \atop{ v^{\prime} \to v}} \frac{ f(\sigma(X) + \tau Q_{\pm} v^{\prime}) - f(\sigma(X))}{\tau} \\
     	                       &\geq  \liminf_{ {\tau \downarrow 0} \atop{ w \to Q_{\pm} v}} \frac{ f(\sigma(X) + \tau w) - f(\sigma(X))}{\tau} \\
     	                       &=  {\rm d}f(\sigma(X))(Q_{\pm} v).
     \end{align*}
     According to $Q_{\pm}^{-1} = {\rm diag}(Q_1^{-1}, \ldots,Q_t^{-1}, Q_{t+1}^{-1}) \in \mathbf{P}^n_{\pm}(X)$, we can similarly show the opposite inequality above, and the proof is complete.
     \qed

	\section{Second Subderivatives of Orthogonally Invariant Matrix Functions}	
	 This section is dedicated to the study of second order variational analysis of orthogonally invariant matrix functions. Our objective is to compute the second subderivatives of such functions when the corresponding absolutely symmetric functions are convex. Additionally, we derive second-order optimality conditions for a  class of matrix optimization problems. We begin our investigation of the second subderivative of orthogonally invariant matrix functions by establishing a lower bound for it.
    
     \begin{prop}
      Let $f \colon \R^n \rightarrow [-\infty,+\infty]$ be lsc, convex, and absolutely symmetric. Assume that $\mu_1> \cdots > \mu_t$ are the distinct nonzero singular values of $X \in {\mathbb M}_{m,n}$, and that  $(U,V) \in {\cal O}^{m,n}(X) \cap {\cal O}^{m,n}(Y)$. Then, for any $H \in {\mathbb M}_{m,n}$, we have 
    	\begin{align}\label{a lower estimate for second subderivative}
    		\begin{split}
    			{\rm d^2} (f \circ \sigma)(X \mid Y)(H) &\geq {\rm d}^2 f(\sigma(X) \mid \sigma(Y))(\sigma^{\prime}(X;H))   \\
    			&+2\sum_{i=1}^t \left\langle \Sigma(Y)_{\alpha_i\alpha_i}, P_{\alpha_i}^{T}{\cal B}(H) P_{\alpha_i}^c (\mu_iI-\varLambda_{\alpha_i})^{-1} {P_{\alpha_i}^c}^{T}{\cal B}(H)P_{\alpha_i} \right\rangle \\
    			&+2 \left\langle \Sigma(Y)_{\hat{\beta} \beta}, -U_{\hat{\beta}}^{T} H V_{\alpha} \Sigma_{\alpha}^{-1}(X) U_{\alpha}^{T} H V_{\beta} \right\rangle,
    		\end{split}
    	\end{align}
      where the columns of $P_{\alpha_i}$ form an orthonormal basis of eigenvectors of ${\cal B}(X) $ associated with $\mu_i$, $P_{\alpha_i}^c$ is the submatrix of $P$ obtained by removing all the columns of $P_{\alpha_i}$, and $\varLambda_{\alpha_i} \in \R^{(m+n-|\alpha_i |)\times (m+n-|\alpha_i |)}$ is a diagonal matrix whose diagonal elements are the eigenvalues of $P^{T} {\cal B}(X) P$ that are not equal to $\mu_i$ for all $i \in \{ 1, \ldots,t\}$.
     \end{prop}
	
	 {\it Proof} Let $H \in {\mathbb M}_{m,n}$ and pick sequences $H_K \to H$ and $t_k \downarrow 0$. 
	 Setting $\bigtriangleup_{t_k} \sigma(X)(H_k):= \frac{\sigma(X+t_kH_k)-\sigma(X)}{t_k}$, we obtain 
	 \begin{align*}
	   &\bigtriangleup_{t_k}^2 (f \circ \sigma)(X \mid Y)(H_k) = \frac{f(\sigma(X+t_kH_k)) - f(\sigma(X)) - t_k \langle Y,H_k \rangle}{\frac{1}{2}t_k^2} \\
	  =&\frac{f(\sigma(X)+t_k \bigtriangleup_{t_k} \sigma(X)(H_k)) - f(\sigma(X)) - t_k \langle \sigma(Y),\bigtriangleup_{t_k} \sigma(X)(H_k) \rangle}{\frac{1}{2}t_k^2} + \frac{\langle \sigma(Y),\bigtriangleup_{t_k} \sigma(X)(H_k) \rangle - \langle Y,H_k \rangle}{\frac{1}{2}t_k} \\
	  =& \bigtriangleup_{t_k}^2 f(\sigma(X) \mid \sigma(Y))(\bigtriangleup_{t_k} \sigma(X)(H_k)) + \frac{\langle \sigma(Y),\bigtriangleup_{t_k} \sigma(X)(H_k) \rangle - \langle Y,H_k \rangle}{\frac{1}{2}t_k}.		
	 \end{align*}
     It follows from $Y = U \Sigma(Y) V^T$ that
	 $$
	  \left\langle Y, H_k \right\rangle = \left\langle U \Sigma(Y) V^T, H_k \right\rangle = \left\langle \Sigma(Y) , U^T H_k V \right\rangle = \sum_{i=1}^t \left\langle \Sigma(Y)_{\alpha_i\alpha_i}, U_{\alpha_i}^T H_k V_{\alpha_i} \right\rangle + \left\langle \Sigma(Y)_{\hat{\beta} \beta}, U_{\hat{\beta}}^T H_k V_{\beta} \right\rangle.
	 $$
     Moreover, it results from Fan’s inequality \eqref{Fan’s inequality}, von Neumann’s trace inequality \eqref{Von Neumann’s Trace Inequality}, and Proposition \ref{fdd of SV} that
	 \begin{align*}
	    &   \langle \sigma(Y),\bigtriangleup_{t_k} \sigma(X)(H_k) \rangle \\
	   =&    \sum_{i=1}^t \sum_{s \in \alpha_i } \frac{\sigma_s(Y) (\sigma_s(X+t_kH_k) - \sigma_s(X))}{t_k} + \sum_{s \in \beta } \frac{\sigma_s(Y) (\sigma_s(X+t_kH_k) - \sigma_s(X))}{t_k} \\
	   =&    \sum_{i=1}^t \sum_{s \in \alpha_i } \sigma_s(Y) \lambda_{l_s} \big(P_{\alpha_i}^{T} {\cal B}(H_k) P_{\alpha_i} + t_k P_{\alpha_i}^{T}{\cal B}(H_k) P_{\alpha_i}^c (\mu_i I-\varLambda_{\alpha_i})^{-1} {P_{\alpha_i}^c}^{T}{\cal B}(H_k)P_{\alpha_i} \big) + O(t_k^2) \\
	    &  + \sum_{s \in \beta } \sigma_s(Y) \sigma_{l_s} \big(U_{\hat{\beta}}^{T} H_k V_{\beta} - t_k U_{\hat{\beta}}^{T} H_k V_{\alpha} \Sigma_{\alpha}^{-1}(X) U_{\alpha}^{T} H_k V_{\beta} \big) + O(t_k^2) \\
    \geq&  \sum_{i=1}^t \left\langle \Sigma(Y)_{\alpha_i\alpha_i}, P_{\alpha_i}^{T} {\cal B}(H_k) P_{\alpha_i} + t_k P_{\alpha_i}^{T}{\cal B}(H_k) P_{\alpha_i}^c (\mu_i I-\varLambda_{\alpha_i})^{-1} {P_{\alpha_i}^c}^{T}{\cal B}(H_k)P_{\alpha_i} \right\rangle \\
        &  + \left\langle \Sigma(Y)_{\hat{\beta} \beta}, U_{\hat{\beta}}^{T} H_k V_{\beta} - t_k U_{\hat{\beta}}^{T} H_k V_{\alpha} \Sigma_{\alpha}^{-1}(X) U_{\alpha}^{T} H_k V_{\beta} \right\rangle + O(t_k^2).
	 \end{align*}
     Thus, we have 
	 \begin{align*}
		&  \frac{ \left\langle \sigma(Y),\bigtriangleup_{t_k} \sigma(X)(H_k) \right\rangle - \left\langle Y,H_k \right\rangle}{\frac{1}{2}t_k} \\
	\geq&  2\sum_{i=1}^t \left\langle \Sigma(Y)_{\alpha_i \alpha_i}, \frac{1}{2} \big(V_{\alpha_i}^T H_k^T U_{\alpha_i} - U_{\alpha_i}^T H_k V_{\alpha_i} \big)  + P_{\alpha_i}^{T}{\cal B}(H_k) P_{\alpha_i}^c (\mu_i I-\varLambda_{\alpha_i})^{-1} {P_{\alpha_i}^c}^{T}{\cal B}(H_k)P_{\alpha_i} \right\rangle \\
		&  + 2 \left\langle \Sigma(Y)_{\hat{\beta} \beta},  - U_{\hat{\beta}}^{T} H_k V_{\alpha} \Sigma_{\alpha}^{-1}(X) U_{\alpha}^{T} H_k V_{\beta} \right\rangle + O(t_k^2) \\
	   =&  2\sum_{i=1}^t \left\langle \Sigma(Y)_{\alpha_i \alpha_i},  P_{\alpha_i}^{T}{\cal B}(H_k) P_{\alpha_i}^c (\mu_i I-\varLambda_{\alpha_i})^{-1} {P_{\alpha_i}^c}^{T}{\cal B}(H_k)P_{\alpha_i} \right\rangle \\
		&  + 2 \left\langle \Sigma(Y)_{\hat{\beta} \beta},  - U_{\hat{\beta}}^{T} H_k V_{\alpha} \Sigma_{\alpha}^{-1}(X) U_{\alpha}^{T} H_k V_{\beta} \right\rangle + O(t_k^2),
	 \end{align*}
     and further
     \begin{align*}
    	  \bigtriangleup_{t_k}^2 (f \circ \sigma)(X \mid Y)(H_k) &\geq \bigtriangleup_{t_k}^2 f(\sigma(X) \mid \sigma(Y))(\bigtriangleup_{t_k} \sigma(X)(H_k)) \\
    	&  + 2\sum_{i=1}^t \left\langle \Sigma(Y)_{\alpha_i \alpha_i},  P_{\alpha_i}^{T}{\cal B}(H_k) P_{\alpha_i}^c (\mu_i I-\varLambda_{\alpha_i})^{-1} {P_{\alpha_i}^c}^{T}{\cal B}(H_k)P_{\alpha_i} \right\rangle  \\
    	&  + 2 \left\langle \Sigma(Y)_{\hat{\beta} \beta},  - U_{\hat{\beta}}^{T} H_k V_{\alpha} \Sigma_{\alpha}^{-1}(X) U_{\alpha}^{T} H_k V_{\beta} \right\rangle + O(t_k^2),
     \end{align*}
     which, in turn, brings us to the lower estimate in \eqref{a lower estimate for second subderivative} for the second subderivative of $f \circ \sigma$ at $X$ for $Y$ because $\bigtriangleup_{t_k} \sigma(X)(H_k) \to \sigma^{\prime}(X;H)$ as $k \to \infty$.
     \qed
	
	 We now advance to a discussion concerning the critical cone of orthogonally invariant matrix functions.
	 \begin{prop}\label{Critical Cone of Orthogonally Invariant Matrix Functions}(Critical Cone of Orthogonally Invariant Matrix Functions).
	  Let $f \colon \R^n \rightarrow [-\infty,+\infty]$ be lsc, convex, and absolutely symmetric. Assume that $\mu_1> \cdots > \mu_t$ are the distinct nonzero singular values of $X \in {\mathbb M}_{m,n}$, and that $Y \in \partial (f \circ \sigma)(X)$. Then, we have $H \in  K_{f \circ \sigma}(X,Y)$ if and only if $\sigma^{\prime}(X;H) \in  K_f(\sigma(X),\sigma(Y))$, the matrices $\Sigma(Y)_{\alpha_i \alpha_i}$ and $\frac{1}{2} \big( U_{\alpha_i}^T H V_{\alpha_i} + V_{\alpha_i}^T H^T U_{\alpha_i} \big)$ have a simultaneous ordered spectral decomposition, and the matrices $\Sigma(Y)_{\hat{\beta} \beta}$ and $U_{\hat{\beta}}^T H V_{\beta}$ have a simultaneous ordered singular value decomposition  for any $i = 1,\ldots,t$ and $(U,V) \in {\cal O}^{m,n}(X) \cap {\cal O}^{m,n}(Y)$.
	 \end{prop}
	
	 {\it Proof} It follows from $Y \in \partial (f \circ \sigma)(X)$ and Proposition \ref{subdifferential of orthogonally invariant matrix functions} that $\sigma(Y) \in \partial f(\sigma(X))$ and we can find $(U,V) \in {\cal O}^{m,n}(X) \cap {\cal O}^{m,n}(Y)$. From Theorem \ref{Subderivatives of Orthogonally Invariant Matrix Functions}, we can conclude that
	 \begin{align*}
		   &  \left\langle Y,H \right\rangle = \left\langle \Sigma(Y), U^T H V \right\rangle \\
		=  & \sum_{i=1}^{t} \left\langle \Sigma(Y)_{\alpha_i \alpha_i}, U_{\alpha_i}^T H V_{\alpha_i} \right\rangle + \left\langle \Sigma(Y)_{\hat{\beta} \beta}, U_{\hat{\beta}}^T H V_{\beta} \right\rangle  \\
		=  & \sum_{i=1}^{t} \left\langle \Sigma(Y)_{\alpha_i \alpha_i}, U_{\alpha_i}^T H V_{\alpha_i} \right\rangle + \sum_{i=1}^{t} \left\langle \Sigma(Y)_{\alpha_i \alpha_i}, \frac{1}{2}\big(V_{\alpha_i}^T H^T U_{\alpha_i} - U_{\alpha_i}^T H V_{\alpha_i} \big) \right\rangle + \left\langle \Sigma(Y)_{\hat{\beta} \beta}, U_{\hat{\beta}}^T H V_{\beta} \right\rangle  \\
		=  & \sum_{i=1}^{t} \left\langle \Sigma(Y)_{\alpha_i \alpha_i}, \frac{1}{2} \big( U_{\alpha_i}^T H V_{\alpha_i} + V_{\alpha_i}^T H^T U_{\alpha_i} \big) \right\rangle + \left\langle \Sigma(Y)_{\hat{\beta} \beta}, U_{\hat{\beta}}^T H V_{\beta} \right\rangle  \\
	  \leq & \sum_{i=1}^{t} \left\langle \lambda \big( \Sigma(Y)_{\alpha_i \alpha_i} \big), \frac{1}{2} \lambda \big(U_{\alpha_i}^T H V_{\alpha_i} + V_{\alpha_i}^T H^T U_{\alpha_i} \big) \right\rangle + \left\langle \sigma \big(\Sigma(Y)_{\hat{\beta} \beta} \big), \sigma \big(U_{\hat{\beta}}^T H V_{\beta} \big) \right\rangle  \\
		=  & \sum_{i=1}^{t}  \sum_{s \in\alpha_i} \left\langle \sigma_s(Y),\frac{1}{2} \lambda_{l_s} \big(U_{\alpha_i}^T H V_{\alpha_i} + V_{\alpha_i}^T H^T U_{\alpha_i} \big) \right\rangle + \sum_{s \in \beta} \left\langle \sigma_s(Y), \sigma_{l_s} \big(U_{\hat{\beta}}^T H V_{\beta} \big) \right\rangle  \\
		=  & \sum_{i=1}^{t}  \sum_{s \in\alpha_i} \left\langle \sigma_s(Y),\sigma_s^{\prime}(X;H) \right\rangle + \sum_{s \in \beta} \left\langle \sigma_s(Y), \sigma_s^{\prime}(X;H) \right\rangle  \\
		=  & \left\langle \sigma(Y), \sigma^{\prime}(X;H) \right\rangle  
		\leq  {\rm d}f(\sigma(X))(\sigma^{\prime}(X;H)) 
		={\rm d}(f \circ \sigma)(X)(H).
	 \end{align*}
	 The condition $H \in  K_{f \circ \sigma}(X,Y)$ is equivalent to equality between the left and right hand sides (and hence throughout), and the claimed equivalence then results from $\sigma^{\prime}(X;H) \in  K_f(\sigma(X),\sigma(Y))$, Fan’s inequality \eqref{Fan’s inequality}, and von Neumann’s trace theorem (Lemma \ref{Von Neumann’s Trace Theorem}).
	 \qed

	 The following result provides sufficient conditions for the parabolic epi-differentiability of orthogonally invariant matrix functions. Furthermore, it derives a valuable formula for the parabolic subderivatives of this class of functions.
	
	 \begin{thm}\label{PS of OIMF}(Parabolic Subderivatives of Orthogonally Invariant Matrix Functions).
	  Let $f \colon \R^n \rightarrow [-\infty,+\infty]$ be an absolutely symmetric function and let $X \in {\mathbb M}_{m,n}$ with $(f \circ \sigma)(X)$ finite. Assume that $H \in T_{{\rm dom}(f \circ \sigma)} (X)$ and that $f$ is locally Lipschitz continuous relative to its domain and parabolically epi-differentiable at $\sigma(X)$ for $\sigma^{\prime}(X;H)$. Then the parabolic subderivative of orthogonally invariant matrix function $f \circ \sigma$ at $X$ for $H$ with respect to $W \in {\mathbb M}_{m,n}$ is
	  \begin{equation}\label{formula of PS}
	  	{\rm d}^2(f \circ \sigma)(X)(H \mid W) = {\rm d}^2 f(\sigma(X))(\sigma^{\prime}(X;H) \mid \sigma^{\prime \prime}(X;H,W)),
	  \end{equation}
	  and $f \circ \sigma$ is parabolically epi-differentiable at $X$ for $H$.
	 \end{thm}

     {\it Proof} Let $F:= f \circ \sigma$ and pick $W \in {\mathbb M}_{m,n}$. Since $f$ is parabolically epi-differentiable at $\sigma(X)$ for $\sigma^{\prime}(X;H)$, ${\rm dom}f$ is parabolically derivable at $\sigma(X)$ for $\sigma^{\prime}(X;H)$. We conclude from \eqref{def of domF}, Proposition \ref{prop of SOTS}, Corollary \ref{Tangent Cone to the Orthogonally Invariant Matrix Sets}, and \eqref{formula of SOTS} that ${\rm dom}F$ is parabolically derivable at $X$ for $H$ and the tangent cone of ${\rm dom}F$ at $X$ and second-order tangent set of ${\rm dom}F$ at $X$ for $H \in T_{{\rm dom}F}(X)$ are, respectively, 
     \begin{equation*}
    	T_{{\rm dom}F}(X) = \{ H \in {\mathbb M}_{m,n} \colon \sigma^{\prime}(X;H) \in T_{{\rm dom}f} (\sigma(X))\}
     \end{equation*}
     and
     \begin{equation}\label{SOTS of domF}
    	T_{{\rm dom}F}^2 (X,H) = \{ W \in {\mathbb M}_{m,n} \colon \sigma^{\prime \prime}(X;H,W) \in  T_{{\rm dom}f}^2(\sigma(X),\sigma^{\prime}(X;H))\} \neq \emptyset.
     \end{equation}
     Since $H \in T_{{\rm dom}F}(X)$ amounts to $\sigma^{\prime}(X;H) \in T_{{\rm dom}f} (\sigma(X))$, we deduce from the imposed assumption on $f$ that
     \begin{equation}\label{dom d^2f}
    	{\rm dom}\,{\rm d}^2 f(\sigma(X)) (\sigma^{\prime}(X;H)\mid \cdot)= T_{{\rm dom}f}^2 (\sigma(X), \sigma^{\prime}(X;H)).
     \end{equation}
     The proof falls naturally into two cases. 
     Let us first assume $W \notin T_{{\rm dom}F}^2 (X,H)$. Combining \eqref{SOTS of domF} with \eqref{dom d^2f}, we obtain that 
     $$
      {\rm d}^2 f(\sigma(X)) (\sigma^{\prime}(X;H)\mid \sigma^{\prime \prime}(X;H,W)) = \infty.
     $$ 
     On the other hand, by definition, it is not hard to see that the inclusion ${\rm dom}\,{\rm d}^2 F(X)(H \mid \cdot) \subset T_{{\rm dom}F}^2 (X,H)$ always holds for any $H \in T_{{\rm dom}F}(X)$. Thus we have ${\rm d}^2 F(X)(H \mid W) = \infty$. 
     This implies \eqref{formula of PS} for every $W \notin T_{{\rm dom}F}^2 (X,H)$. Now consider an arbitrary sequence $t_k \downarrow 0$,  set $W_k:=W$ for all $k \in {\mathbb N}$, and obtain that 
     \begin{align*}
    	 {\rm d}^2 F(X)(H \mid W) &\leq \liminf_{k \to \infty}\frac{F(X+ t_k H + \frac{1}{2} t_k^2 W_k)-F(X)-t_k {\rm d}F(X)(H)}{\frac{1}{2} t_k^2} \\
    	&\leq \limsup_{k \to \infty}\frac{F(X+ t_k H + \frac{1}{2} t_k^2 W_k)-F(X)-t_k {\rm d}F(X)(H)}{\frac{1}{2} t_k^2}.
     \end{align*}
     Since $ {\rm d}^2 F(X)(H \mid W) = \infty $ for all $W \notin T_{{\rm dom}F}^2 (X,H)$, the inequality in the above expression can take equality.
    
     Turn to the case $W \in T_{{\rm dom}F}^2 (X,H)$. Combining \eqref{SOTS of domF} with \eqref{dom d^2f}, we obtain that 
     \begin{equation}\label{d^2f in second case}
    	{\rm d}^2 f(\sigma(X)) (\sigma^{\prime}(X;H)\mid \sigma^{\prime \prime}(X;H,W)) < \infty.
     \end{equation}
     Since ${\rm dom}F$ is parabolically derivable at $X$ for $H$, for arbitrary sequence $t_k \downarrow 0$, there exists a sequence $W_k \to W$ as $k \to \infty$ such that
     $$
	  X_k := X + t_kH + \frac{1}{2}t_k^2 W_k = X + t_kH + \frac{1}{2}t_k^2 W +o(t_k^2) \in {\rm dom}F.
     $$
     Because $f$ is parabolically epi-differentiable at $\sigma(X)$ for $\sigma^{\prime}(X;H)$, for the vector $\sigma^{\prime \prime}(X;H,W) \in \R^n$, corresponding to the above $t_k$, we find a sequence $w_k \to \sigma^{\prime \prime}(X;H,W)$ such that
     \begin{equation*}
    	{\rm d}^2 f(\sigma(X))(\sigma^{\prime}(X;H)  \mid \sigma^{\prime \prime}(X;H,W)) = \lim_{k \to \infty} \frac{f(y_k)-f(\sigma(X))-t_k {\rm d}f(\sigma(X))(\sigma^{\prime}(X;H))}{\frac{1}{2} t_k^2},
     \end{equation*}
     where $y_k := \sigma(X)+ t_k \sigma^{\prime}(X;H) + \frac{1}{2} t_k^2 w_k$. It follows from \eqref{d^2f in second case} that $y_k \in {\rm dom}f$ for all $k$ sufficiently large. We deduce from \eqref{a chain rule for subderivatives of OIMF} and \eqref{estimate for the SVF} that 
     \begin{align*}
    	{\rm d}^2F(X)(H \mid W) &\leq \liminf_{k \to \infty}\frac{F(X + t_kH + \frac{1}{2}t_k^2 W_k)-F(X)-t_k {\rm d}F(X)(H)}{\frac{1}{2} t_k^2} \\
    	                        &\leq \limsup_{k \to \infty}\frac{F(X_k)-F(X)-t_k {\rm d}F(X)(H)}{\frac{1}{2} t_k^2}  \\
    	                        &=    \limsup_{k \to \infty}\frac{f(\sigma(X_k))-f(\sigma(X))-t_k {\rm d}f(\sigma(X))(\sigma^{\prime}(X;H))}{\frac{1}{2} t_k^2}  \\
    	                        &\leq \limsup_{k \to \infty}\frac{f(y_k)-f(\sigma(X))-t_k {\rm d}f(\sigma(X))(\sigma^{\prime}(X;H))}{\frac{1}{2} t_k^2} + \limsup_{k \to \infty}\frac{f(\sigma(X_k))-f(y_k)}{\frac{1}{2} t_k^2} \\
    	                        &\leq {\rm d}^2 f(\sigma(X))(\sigma^{\prime}(X;H)  \mid \sigma^{\prime \prime}(X;H,W)) + c \limsup_{k \to \infty} \left\| \sigma^{\prime \prime}(X;H,W)+\frac{o(t_k^2)}{\frac{1}{2}t_k^2} - w_k \right\|  \\
    	                        &=    {\rm d}^2 f(\sigma(X))(\sigma^{\prime}(X;H)  \mid \sigma^{\prime \prime}(X;H,W)),
     \end{align*}
     where $c \geq 0$ is a Lipschitz constant of $f$ around $\sigma(X)$ relative to its domain. On the other hand, we conclude from the definition of parabolic subderivative that for any sequence $t_k \downarrow 0$ and any sequence $W_k \to W$ as $k \to \infty$, one has
     \begin{align*}
    	 &   \liminf_{k \to \infty}\frac{F(X + t_kH + \frac{1}{2}t_k^2 W_k)-F(X)-t_k {\rm d}F(X)(H)}{\frac{1}{2} t_k^2} \\
    	=&    \liminf_{k \to \infty}\frac{f(\sigma(X) + t_k\sigma^{\prime}(X;H) + \frac{1}{2}t_k^2 \sigma^{\prime\prime}(X;H,W) +o(t_k^2) )-f(\sigma(X))-t_k {\rm d}f(\sigma(X))(\sigma^{\prime}(X;H))}{\frac{1}{2} t_k^2}  \\
     \geq& \liminf_{t\downarrow 0 \atop{w^{\prime} \to \sigma^{\prime\prime}(X;H,W)}}\frac{f(\sigma(X) + t\sigma^{\prime}(X;H) + \frac{1}{2}t^2 w^{\prime} )-f(\sigma(X))-t {\rm d}f(\sigma(X))(\sigma^{\prime}(X;H))}{\frac{1}{2} t^2}  \\
    	=&   {\rm d}^2 f(\sigma(X))(\sigma^{\prime}(X;H)  \mid \sigma^{\prime \prime}(X;H,W)),
     \end{align*}
     which verifies the inequality "$\geq$" in \eqref{formula of PS}. 
     In summary, for any $W \in T_{{\rm dom}F}^2 (X,H)$ and any sequence $t_k \downarrow 0$, there exists a sequence $W_k \to W$ such that 
     $$
     {\rm d}^2F(X)(H \mid W) =\lim_{k \to \infty}\frac{F(X + t_kH + \frac{1}{2}t_k^2 W_k)-F(X)-t_k {\rm d}F(X)(H)}{\frac{1}{2} t_k^2},
     $$
     and
     $$
     {\rm d}^2F(X)(H \mid W) = {\rm d}^2 f(\sigma(X))(\sigma^{\prime}(X;H)  \mid \sigma^{\prime \prime}(X;H,W)),
     $$
     which coupled with \eqref{d^2f in second case}, leads us to $T_{{\rm dom}F}^2 (X,H) \subset {\rm dom}\,{\rm d}^2F(X)(H \mid \cdot)$. We conclude from $T_{{\rm dom}F}^2 (X,H)\neq \emptyset$ that ${\rm dom}\,{\rm d}^2F(X)(H \mid \cdot) \neq \emptyset$, and thus $F$ is parabolically epi-differentiable at $X$ for $H$. 
     \qed

     Recall from Proposition \ref{symmetric property of subderivative function} that the subderivative of $f$ at $\sigma(X)$ is a symmetric function with respect to a subset of $\mathbf{P}^n_{\pm}$. The following result shows the parabolic subderivative of $f$ at $\sigma(X)$ for $\sigma^{\prime}(X;H)$ is also a symmetric function with respect to a subset of $\mathbf{P}^n_{\pm}$. 
     Now suppose that $N_i$ is the number of distinct eigenvalues of $\frac{1}{2} \big(U_{\alpha_i}^T H V_{\alpha_i} + V_{\alpha_i}^T H^T U_{\alpha_i} \big)$ for $i=1, \ldots, t$ and that $N_{t+1}$ is the number of distinct nonzero singular values of $U_{\hat{\beta}}^T H V_{\beta}$. Pick the index sets $\beta^i_j$ for $i=1, \ldots, t$, $j=1, \ldots, N_i$ from \eqref{def of beta^i_j} and $\beta^{t+1}_k$ for $k=1, \ldots, N_{t+1}, N_{t+1}+1$ from \eqref{def of beta^{t+1}_j}. 
     Denote by $\mathbf{P}^n_{\pm}(X,H)$ a subset of $\mathbf{P}^n_{\pm}$ consisting of all $n \times n$ block diagonal matrices in the form $Q_{\pm} = {\rm diag}(Q_1, \ldots,Q_t,Q_{t+1})$ such that for each $i=1, \ldots, t$, the $|\alpha_i| \times |\alpha_i|$ signed permutation matrix $Q_i$ has a block diagonal representation 
     $$
     Q_i = {\rm diag} \left(B_1^i, \ldots,B_{N_i}^i \right),
     $$ 
     and $|\beta| \times |\beta|$ signed permutation matrix $Q_{t+1}$ has a block diagonal representation 
     $$
     Q_{t+1} = {\rm diag} \left( B_1^{t+1}, \ldots,B_{N_{t+1}}^{t+1}, B_{N_{t+1}+1}^{t+1} \right), 
     $$
     where $ B_j^i \in \R^{|\beta^i_j|\times |\beta^i_j|} $ is a signed permutation matrix for any $i=1, \ldots, t$, $j=1, \ldots, N_i$ and $ B_k^{t+1} \in \R^{|\beta^{t+1}_k|\times |\beta^{t+1}_k|} $ is also a signed permutation matrix for any $k=1, \ldots, N_{t+1}, N_{t+1}+1$. 
     It is obvious that $\mathbf{P}^n_{\pm}(X,H) \subset \mathbf{P}^n_{\pm}(X) \subset \mathbf{P}^n_{\pm}$ and that if $Q_{\pm} \in \mathbf{P}^n_{\pm}(X,H)$, then $Q_{\pm} \sigma(X) = \sigma(X)$ and $Q_{\pm} \sigma^{\prime}(X;H) = \sigma^{\prime}(X;H)$. 
     \begin{prop}\label{symmetric property of parabolic subderivative function}
    	Let $f \colon \R^n \rightarrow [-\infty,+\infty]$ be an absolutely symmetric function and $X,H \in {\mathbb M}_{m,n}$ with $f(\sigma(X))$ and ${\rm d}f(\sigma(X))(\sigma^{\prime}(X;H))$ finite. Then, for any $u \in \R^n$ and any signed permutation matrix $Q_{\pm} \in \mathbf{P}^n_{\pm}(X;H)$, we have
    	\begin{equation*}
    		{\rm d}^2 f(\sigma(X))(\sigma^{\prime}(X;H) \mid Q_{\pm} u) = {\rm d}^2 f(\sigma(X))(\sigma^{\prime}(X;H) \mid u),
    	\end{equation*}
    	which means that the parabolic subderivative $u \mapsto {\rm d}^2 f(\sigma(X))(\sigma^{\prime}(X;H) \mid u)$ is symmetric with respect to $\mathbf{P}^n_{\pm}(X;H)$.
     \end{prop}
	  
	 {\it Proof} This follows by the same method as in Proposition \ref{symmetric property of subderivative function}. 
	 \qed

	 As mentioned before, we can partition any vector $z \in \R^n$ into $z = (z_{\alpha_1}, \ldots, z_{\alpha_t}, z_{\beta})$ with $\alpha_i$, $i = 1, \ldots, t$, taken from \eqref{alpha_i}. 
	 Pick $i \in \{1, \ldots, t\}$ and recall from \eqref{def of beta^i_j} that the index set $\alpha_i = \cup_{j=1}^{N_i} \beta^i_j$. This allows us to partition further $z_{\alpha_i}$ and $z_{\beta}$, respectively, into $z_{\alpha_i} = (z_{\beta^i_1}, \ldots, z_{\beta^i_{N_i}})$ and $z_{\beta} = (z_{\beta^{t+1}_1}, \ldots, z_{\beta^{t+1}_{N_{t+1}}}, z_{\beta^{t+1}_{N_{t+1}+1}})$, where $z_{\beta^i_j} \in \R^{|\beta^i_j|}$ for $i \in \{1, \ldots, t\}, j \in \{1, \ldots, N_i\}$, and $z_{\beta^{t+1}_k} \in \R^{|\beta^{t+1}_k|}$ for $k \in \{1, \ldots, N_{t+1}, N_{t+1}+1\}$. In short, we can equivalently write $z$ as
	 \begin{equation}\label{representation of z}
		 \left(z_{\beta^1_1}, \ldots, z_{\beta^1_{N_1}}, \ldots, z_{\beta^t_1}, \ldots, z_{\beta^t_{N_t}}, z_{\beta^{t+1}_1}, \ldots, z_{\beta^{t+1}_{N_{t+1}}}, z_{\beta^{t+1}_{N_{t+1}+1}} \right),
	 \end{equation}
	 where $t$, $N_{t+1}$, taken from \eqref{alpha_i} and \eqref{def of beta^{t+1}_j}, and $N_i, i = 1, \ldots, t$, taken from \eqref{def of beta^i_j}, stand for the number of distinct nonzero singular values of $X \in {\mathbb M}_{m,n}$ and $U_{\hat{\beta}}^T H V_{\beta}$, and the number of distinct eigenvalues of $\frac{1}{2} \big(U_{\alpha_i}^T H V_{\alpha_i} + V_{\alpha_i}^T H^T U_{\alpha_i} \big)$, respectively. 
	 Thus, the representation of $z$ in \eqref{representation of z} is associated with the signed permutation matrices with representation 
	 $$
	   {\rm diag} \left(B_1^1, \ldots, B_{N_1}^1, \ldots, B_1^t, \ldots,B_{N_t}^t, B_1^{t+1}, \ldots,B_{N_{t+1}}^{t+1}, B_{N_{t+1}+1}^{t+1}\right)
	 $$
	 in $\mathbf{P}^n_{\pm}(X,H)$, where $ B_j^i \in \R^{|\beta^i_j|\times |\beta^i_j|} $ is an signed permutation matrix for any $i=1, \ldots, t$, $j=1, \ldots, N_i$ and $ B^{t+1}_k \in \R^{|\beta^{t+1}_k|\times |\beta^{t+1}_k|} $ is also an signed permutation matrix for any $k=1, \ldots, N_{t+1}, N_{t+1}+1$. Denote by $\R^n_{\downarrow}$ the set of all vectors $(v_1, \ldots, v_n)$ such that $v_1 \geq \ldots \geq v_n$.

     \begin{prop}\label{hat W and bar z}
     	 Let $f \colon \R^n \rightarrow [-\infty,+\infty]$ be lsc, convex, and absolutely symmetric. Assume that $Y \in \partial (f \circ \sigma)(X)$ and $H \in  K_{f \circ \sigma}(X,Y)$. If $f$ is parabolically regular at $\sigma(X)$ for $\sigma(Y)$, then the following assertions hold:
    	
    	{\rm (i)} There exists $\bar z \in \R^n$ with representation \eqref{representation of z}, where $z_{\beta^i_j} \in \R^{|\beta^i_j|}_{\downarrow}$ for $i \in \{1, \ldots, t\}, j \in \{1, \ldots, N_i\}$ and  $z_{\beta^{t+1}_k} \in \R^{|\beta^{t+1}_k|}_{\downarrow}$ for $k \in \{1, \ldots, N_{t+1}, N_{t+1}+1\}$, satisfying 
    	\begin{equation}\label{z satisfing min}
    	    {\rm d}^2 f(\sigma(X) \mid \sigma(Y))(\sigma^{\prime}(X;Y)) = {\rm d}^2 f(\sigma(X))(\sigma^{\prime}(X;Y)  \mid \bar z) - \langle \sigma(Y), \bar z \rangle.
    	\end{equation}
    
        {\rm (ii)} There exists a matrix $\widehat{W} \in {\mathbb M}_{m,n}$ such that $\sigma^{\prime \prime}(X;H,\widehat{W}) = \bar z$, where $\bar z$ satisfies the above condition {\rm (i)}.
     \end{prop}
 
     {\it Proof} Because $Y \in \partial (f \circ \sigma)(X)$ and $H \in  K_{f \circ \sigma}(X,Y)$, $\sigma(Y) \in \partial f(\sigma(X))$ and $\sigma^{\prime}(X;H) \in K_f(\sigma(X), \sigma(Y))$ by Proposition \ref{subdifferential of orthogonally invariant matrix functions} and Proposition \ref{Critical Cone of Orthogonally Invariant Matrix Functions}, respectively. It follows from the parabolic regularity of $f$ at $\sigma(X)$ for $\sigma(Y) \in \partial f(\sigma(X))$ that 
     \begin{equation*}
    	{\rm d}^2 f(\sigma(X), \sigma(Y))(\sigma^{\prime}(X;H)) = \inf_z \{ {\rm d}^2 f(\sigma(X))(\sigma^{\prime}(X;H)  \mid z) - \langle \sigma(Y), z \rangle \}.
     \end{equation*}
     As explained above, there exists $ \hat z \in \R^n$ with representation \eqref{representation of z} satisfying \eqref{z satisfing min}. It is sufficient to show that the components of each $z_{\beta^i_j} \in \R^{|\beta^i_j|}$ for $i \in \{1, \ldots, t\}, j \in \{1, \ldots, N_i\}$ and  $z_{\beta^{t+1}_k} \in \R^{|\beta^{t+1}_k|}$ for $k \in \{1, \ldots, N_{t+1}, N_{t+1}+1\}$ have nonincreasing order. 
     We choose $|\beta^i_j|\times |\beta^i_j|$ signed permutation matrix $ B_j^i $ for any $i=1, \ldots, t$, $j=1, \ldots, N_i$ and $|\beta^{t+1}_k|\times |\beta^{t+1}_k|$ signed permutation matrix $ B^{t+1}_k $ for $k=1, \ldots, N_{t+1}, N_{t+1}+1$ such that ${\bar z}_{\beta^i_j} = B_j^i z_{\beta^i_j} \in \R^{|\beta^i_j|}_{\downarrow}$ and ${\bar z}_{\beta^{t+1}_k} = B^{t+1}_k z_{\beta^{t+1}_k} \in \R^{|\beta^{t+1}_k|}_{\downarrow}$. 
     Let 
     $$
      Q_{\pm} := {\rm diag}\left( B_1^1, \ldots, B_{N_1}^1, \ldots, B_1^t, \ldots,B_{N_t}^t, B_1^{t+1}, \ldots,B_{N_{t+1}}^{t+1}, B_{N_{t+1}+1}^{t+1} \right),
     $$
     then $Q_{\pm} \in \mathbf{P}^n_{\pm}(X,H)$. 
     Set 
     \begin{equation}\label{bar z}
    	\bar z = \left({\bar z}_{\beta^1_1}, \ldots, {\bar z}_{\beta^1_{N_1}}, \ldots, {\bar z}_{\beta^t_1}, \ldots, {\bar z}_{\beta^t_{N_t}}, {\bar z}_{\beta^{t+1}_1}, \ldots, {\bar z}_{\beta^{t+1}_{N_{t+1}}}, {\bar z}_{\beta^{t+1}_{N_{t+1}+1}} \right)
     \end{equation}
     and observe that $\bar z = Q_{\pm} \hat z$. 
     We deduce from Proposition \ref{symmetric property of parabolic subderivative function} that
     \begin{equation}\label{application about symmetric property of parabolic subderivative function}
    	{\rm d}^2 f(\sigma(X))(\sigma^{\prime}(X;H) \mid \bar z) = {\rm d}^2 f(\sigma(X))(\sigma^{\prime}(X;H) \mid \hat z).
     \end{equation}
     On the other hand, suppose that 
     $$
      \left(\sigma(Y)_{\beta^1_1}, \ldots, \sigma(Y)_{\beta^1_{N_1}}, \ldots, \sigma(Y)_{\beta^t_1}, \ldots, \sigma(Y)_{\beta^t_{N_t}}, \sigma(Y)_{\beta^{t+1}_1}, \ldots, \sigma(Y)_{\beta^{t+1}_{N_{t+1}}}, \sigma(Y)_{\beta^{t+1}_{N_{t+1}+1}} \right)
     $$
     is a partition of the vector $\sigma(Y)$ corresponding to \eqref{representation of z} and observe that $\sigma(Y)_{\beta^i_j} \in \R^{|\beta^i_j|}_{\downarrow}$ for $i \in \{1, \ldots, t\}, j \in \{1, \ldots, N_i\}$ and  $\sigma(Y)_{\beta^{t+1}_k} \in \R^{|\beta^{t+1}_k|}_{\downarrow}$ for $k \in \{1, \ldots, N_{t+1}, N_{t+1}+1\}$. It is not hard to get that
     \begin{align*}
    	      \left \langle  \sigma(Y), \hat z  \right \rangle 
    	&=    \sum_{i=1}^t \sum_{j=1}^{N_i} \left \langle  \sigma(Y)_{\beta^i_j}, z_{\beta^i_j}  \right \rangle +   \sum_{k=1}^{N_{t+1}+1} \left\langle  \sigma(Y)_{\beta^{t+1}_k}, z_{\beta^{t+1}_k}  \right \rangle \\
    	&\leq \sum_{i=1}^t \sum_{j=1}^{N_i} \left \langle  \sigma(Y)_{\beta^i_j}, {\bar z}_{\beta^i_j}  \right \rangle +   \sum_{k=1}^{N_{t+1}+1} \left \langle  \sigma(Y)_{\beta^{t+1}_k}, {\bar z}_{\beta^{t+1}_k}  \right \rangle \\
    	&=    \left \langle  \sigma(Y), \bar z  \right \rangle.
     \end{align*}
     This, together with \eqref{application about symmetric property of parabolic subderivative function} and $\hat z$ with representation \eqref{representation of z} satisfying \eqref{z satisfing min}, justifies
     $$
      {\rm d}^2 f(\sigma(X) \mid \sigma(Y))(\sigma^{\prime}(X;Y)) \geq  {\rm d}^2 f(\sigma(X))(\sigma^{\prime}(X;Y)  \mid \bar z) - \langle \sigma(Y), \bar z \rangle.
     $$
     We conclude from the inequality \eqref{rs of ss and ps} that \eqref{z satisfing min} holds for $\bar z$.
    
     Turning now to proof {\rm (ii)}, pick the vector $\bar z \in \R^n$ from \eqref{bar z}. We can equivalently write via the index sets $\alpha_i$, $i = 1, \ldots,t$, from \eqref{alpha_i} that
     $$
      \bar z = \left(\bar z_{\alpha_1}, \ldots, \bar z_{\alpha_t}, \bar z_{\beta} \right) \,\, \mbox{with}\,\, \bar z_{\alpha_i} = \left({\bar z}_{\beta^i_1}, \ldots, {\bar z}_{\beta^i_{N_i}} \right)  \,\, \mbox{and}\,\, \bar z_{\beta} = \left( {\bar z}_{\beta^{t+1}_1}, \ldots, {\bar z}_{\beta^{t+1}_{N_{t+1}}}, {\bar z}_{\beta^{t+1}_{N_{t+1}+1}} \right),
     $$
     where ${\bar z}_{\beta^i_j} \in \R^{|\beta^i_j|}_{\downarrow}$  for $i \in \{1, \ldots, t\}, j \in \{1, \ldots, N_i\}$, and ${\bar z}_{\beta^{t+1}_k} \in \R^{|\beta^{t+1}_k|}_{\downarrow}$ for $k \in \{1, \ldots, N_{t+1}, N_{t+1}+1\}$. 
     It suffices to prove that there exists a matrix $\widehat{W} \in {\mathbb M}_{m,n}$ satisfying $\sigma^{\prime \prime}(X;H,\widehat{W}) = \bar z$. Pick the $ |\alpha_i| \times |\alpha_i|$ matrix  $Q^i \in {\cal O}^{|\alpha_i|} \Big( \frac{1}{2} \big(U_{\alpha_i}^T H V_{\alpha_i} + V_{\alpha_i}^T H^T U_{\alpha_i} \big) \Big)$, $i = 1, \ldots, t$, and take the two matrices $( Q_{\hat{\beta} \hat{\beta}}, \widehat{Q}_{\beta \beta} ) \in {\cal O}^{|\hat{\beta}|, |\beta|} \big(U_{\hat{\beta}}^T H V_{\beta} \big)$. 
     Consider the $|\alpha| \times |\alpha|$ block diagonal matrix 
     \begin{equation}\label{A_{alpha alpha}}
     	A_{\alpha \alpha} = {\rm diag}\left( Q^1 {\rm diag}(\bar z_{\alpha_1} ) (Q^1)^T, \ldots, Q^t {\rm diag}(\bar z_{\alpha_t} ) (Q^t)^T \right),
     \end{equation}
     and the $|\hat{\beta}| \times |\beta|$ matrix 
     \begin{equation}\label{A_{hat{beta} beta}}
     	A_{\hat{\beta} \beta} = Q_{ \hat{\beta} \hat{\beta}} \left( {\rm diag}(\bar z_\beta) \right)_{\hat{\beta} \beta} \widehat{Q}_{\beta \beta}^T.
     \end{equation}
     It is clear that $A_{\hat{\beta} \beta } = Q_{ \hat{\beta} \beta} \left( {\rm diag}(\bar z_\beta) \right)_{\beta \beta} \widehat{Q}_{\beta \beta}^T$. 
     We claim that there exists a matrix $\widehat{W} \in {\mathbb M}_{m,n}$ satisfying
     \begin{equation}\label{widehat{W} in case 1}
     	P_{\alpha_i}^T {\cal B}(\widehat{W}) P_{\alpha_i} = P_{\alpha_i}^T \Big[ 2{\cal B}(H) P_{\alpha_i}^c (\varLambda_s- \mu_i I)^{-1} {P_{\alpha_i}^c}^{T} {\cal B}(H) +  P_{\alpha} A_{\alpha \alpha} P_{\alpha} ^T \Big]  P_{\alpha_i}
     \end{equation}
     and
     \begin{equation}\label{widehat{W} in case 2}
     	U_{\hat{\beta}}^T \widehat{W} V_{\beta} = U_{\hat{\beta}}^{T} \Big[ 2 H V_{\alpha} \Sigma_{\alpha}^{-1}(X) U_{\alpha}^{T} H + U_{\hat{\beta}} A_{\hat{\beta} \beta} V_{\beta}^T \Big] V_{\beta}.
     \end{equation}
     Indeed, let $W = {\rm diag}(W_{\alpha_1 \alpha_1}, \ldots, W_{\alpha_t \alpha_t}, W_{\hat{\beta} \beta})\in {\mathbb M}_{m,n}$  satisfy that 
     \begin{equation}\label{W_1}
     	W_{\alpha_i \alpha_i} + W_{\alpha_i \alpha_i} ^T 
     	= P_{\alpha_i}^T \Big[ 2{\cal B}(H) P_{\alpha_i}^c (\varLambda_s- \mu_i I)^{-1} {P_{\alpha_i}^c}^{T} {\cal B}(H) +  P_{\alpha} A_{\alpha \alpha} P_{\alpha} ^T \Big]  P_{\alpha_i},
     \end{equation}
     for all $i = 1, \ldots, t$, and 
     \begin{equation}\label{W_2}
      W_{\hat{\beta} \beta} =  U_{\hat{\beta}}^{T} \Big[ 2 H V_{\alpha} \Sigma_{\alpha}^{-1}(X) U_{\alpha}^{T} H + U_{\hat{\beta}} A_{\hat{\beta} \beta} V_{\beta}^T \Big] V_{\beta}.
     \end{equation}
     Set $\widehat{W} = U W V^T$. Then by some elementary calculations, we can obtain that
     $${\cal B}(\widehat{W}) =  P_{\alpha} [U_{\alpha}^T \widehat{W} V_{\alpha} + V_{\alpha}^T \widehat{W}^T U_{\alpha}] P_{\alpha}^T = P_{\alpha} [W_{\alpha \alpha} + W_{\alpha \alpha} ^T] P_{\alpha}^T$$
     and
     $$
      W_{\hat{\beta} \beta} = U_{\hat{\beta}}^T  U W V^T V_{\beta} = U_{\hat{\beta}}^T \widehat{W} V_{\beta}.
     $$
     This, together with \eqref{W_1} and \eqref{W_2}, justifies that \eqref{widehat{W} in case 1} and \eqref{widehat{W} in case 2} hold. Next we show that $ \sigma^{\prime\prime}(X;H,\widehat{W})= \bar z$. 
     For any $s \in \{1, \ldots, n\}$, we first consider the case $s \in \alpha$. There exists $i \in \{ 1, \ldots, t\}$ and $j \in \{1, \ldots,N_i\}$ such that  $s \in \alpha_i$ and $l_s \in \beta^i_j$. It follows from \eqref{sdd of alpha} that 
     \begin{align*}
    	\sigma^{\prime\prime}_s(X;H,\widehat{W}) 
     &= \lambda_{\tilde{l}_s} \bigg( \Big(Q^i_{\beta^i_j} \Big)^T P_{\alpha_i}^T \Big[{\cal B}(\widehat{W})-2{\cal B}(H) P_{\alpha_i}^c (\varLambda_s- \mu_i I)^{-1} {P_{\alpha_i}^c}^{T} {\cal B}(H) \Big] P_{\alpha_i} Q^i_{\beta^i_j} \bigg)\\
     &=  \lambda_{\tilde{l}_s} \bigg( \Big(Q^i_{\beta^i_j} \Big)^T \Big[ P_{\alpha_i}^T P_ {\alpha} A_{\alpha \alpha} P_{\alpha} ^T P_{\alpha_i} \Big] Q^i_{\beta^i_j} \bigg)\\
     &= \lambda_{\tilde{l}_s} \bigg( \Big(Q^i_{\beta^i_j} \Big)^T Q^i {\rm diag}(\bar z_{\alpha_i} ) (Q^i)^T Q^i_{\beta^i_j} \bigg)\\
     &= \lambda_{\tilde{l}_s} \bigg( {\rm diag}\Big({\bar z}_{\beta^i_j}\Big) \bigg) = \Big({\bar z}_{\beta^i_j}\Big)_{\tilde{l}_s} = {\bar z}_s.
     \end{align*}
     We now turn to the case $s \in \beta$. If $l_s(X) \in \beta^{t+1}_k$ for some $k \in \{1, \ldots,N_{t+1}\}$, we deduce from \eqref{sdd of beta >0} that 
     \begin{align*}
    	 \sigma^{\prime\prime}_s(X;H,\widehat{W}) 
        = &\frac{1}{2} \lambda_{\tilde{l}_s} \bigg( \Big(Q_{\hat{\beta} \beta^{t+1}_k} \Big)^T \Big[ U_{\hat{\beta}}^T \widehat{W} V_{\beta}-2 U_{\hat{\beta}}^{T} H V_{\alpha} \Sigma_{\alpha}^{-1}(X_0) U_{\alpha}^{T} H V_{\beta} \Big] \widehat{Q}_{\beta \beta^{t+1}_k} \\
    	  &+ \Big(\widehat{Q}_{\beta \beta^{t+1}_k} \Big)^T \Big[ U_{\hat{\beta}}^T \widehat{W} V_{\beta}-2 U_{\hat{\beta}}^{T} H V_{\alpha} \Sigma_{\alpha}^{-1}(X_0) U_{\alpha}^{T} H V_{\beta} \Big]^T Q_{\hat{\beta} \beta^{t+1}_k}  \bigg) \\
    	= &\frac{1}{2} \lambda_{\tilde{l}_s} \bigg( \Big(Q_{\hat{\beta} \beta^{t+1}_k} \Big)^T \Big[U_{\hat{\beta}}^T U_{\hat{\beta}} A_{\hat{\beta} \beta} V_{\beta}^T V_{\beta} \Big] \widehat{Q}_{\beta \beta^{t+1}_k} 
    	   + \Big(\widehat{Q}_{\beta \beta^{t+1}_k} \Big)^T \Big[ U_{\hat{\beta}}^T U_{\hat{\beta}} A_{\hat{\beta} \beta} V_{\beta}^T V_{\beta} \Big]^T Q_{\hat{\beta} \beta^{t+1}_k}  \bigg) \\
    	= &\frac{1}{2} \lambda_{\tilde{l}_s} \bigg( \Big(Q_{\hat{\beta} \beta^{t+1}_k} \Big)^T Q_{ \hat{\beta} \beta} \left( {\rm diag}(\bar z_\beta) \right)_{\beta \beta} \widehat{Q}_{\beta \beta}^T \widehat{Q}_{\beta \beta^{t+1}_k} 
    	   + \Big(\widehat{Q}_{\beta \beta^{t+1}_k} \Big)^T \widehat{Q}_{\beta \beta} \left( {\rm diag}(\bar z_\beta) \right)_{\beta \beta}Q_{ \hat{\beta} \beta}^T  Q_{\hat{\beta} \beta^{t+1}_k}  \bigg) \\
        = &\lambda_{\tilde{l}_s} \left( {\rm diag}\Big({\bar z}_{\beta^{t+1}_k} \Big) \right)= \Big({\bar z}_{\beta^{t+1}_k}\Big)_{\tilde{l}_s} = {\bar z}_s.
     \end{align*}
     On the other hand, if $l_s(X) \in \beta^{t+1}_{N_{t+1}+1}$,  we conclude from \eqref{sdd of beta =0} that 
     \begin{align*}
    	\sigma^{\prime\prime}_s(X;H,\widehat{W}) 
       = &\sigma_{\tilde{l}_s} \left( 
    	   \begin{bmatrix}
    		 Q_{\hat{\beta} \beta^{t+1}_{N_{t+1}+1}}  & Q_{\hat{\beta} \beta_0}
    	   \end{bmatrix}^T
    	  \Big[ U_{\hat{\beta}}^T \widehat{W} V_{\beta}-2 U_{\hat{\beta}}^{T} H V_{\alpha} \Sigma_{\alpha}^{-1}(X_0) U_{\alpha}^{T} H V_{\beta} \Big] \widehat{Q}_{\beta \beta^{t+1}_{N_{t+1}+1}}  \right)\\
       = &\sigma_{\tilde{l}_s} \left( 
            \begin{bmatrix}
       	      Q_{\hat{\beta} \beta^{t+1}_{N_{t+1}+1}}  & Q_{\hat{\beta} \beta_0}
            \end{bmatrix}^T
            \Big[U_{\hat{\beta}}^T U_{\hat{\beta}} A_{\hat{\beta} \beta} V_{\beta}^T V_{\beta} \Big]
            \widehat{Q}_{\beta \beta^{t+1}_{N_{t+1}+1}}  \right)\\
       = &\sigma_{\tilde{l}_s} \left( 
           \begin{bmatrix}
        	Q_{\hat{\beta} \beta^{t+1}_{N_{t+1}+1}}  & Q_{\hat{\beta} \beta_0}
           \end{bmatrix}^T
           Q_{ \hat{\beta} \hat{\beta}} \left( {\rm diag}(\bar z_\beta) \right)_{\hat{\beta} \beta} \widehat{Q}_{\beta \beta}^T
           \widehat{Q}_{\beta \beta^{t+1}_{N_{t+1}+1}}  \right)\\
       = &\sigma_{\tilde{l}_s} \left( \bigg( {\rm diag}\Big({\bar z}_{\beta^{t+1}_{N_{t+1}+1}} \Big) \bigg)_{(\beta^{t+1}_{N_{t+1}+1} \cup \beta_0) \beta^{t+1}_{N_{t+1}+1}}\right)= \Big({\bar z}_{\beta^{t+1}_{N_{t+1}+1}} \Big)_{\tilde{l}_s} = {\bar z}_s.
     \end{align*}
     In summary, we can find a matrix $\widehat{W} \in {\mathbb M}_{m,n}$ satisfying simultaneously \eqref{widehat{W} in case 1} and \eqref{widehat{W} in case 2} such that $\sigma^{\prime \prime}(X;H,\widehat{W}) = \bar z$, where $\bar z$ comes from {\rm (i)}. 
     \qed

     Finally, we present a exact formula for second subderivative of orthogonally invariant matrix functions.
	 \begin{thm}\label{second subderivative of OIMF}
	 	Let $f \colon \R^n \rightarrow [-\infty,+\infty]$ be an absolutely symmetric function, and be lsc, convex, and locally Lipschitz continuous relative to its domain. Assume that $\mu_1> \cdots > \mu_t$ are the distinct nonzero singular values of $X \in {\mathbb M}_{m,n}$ and $Y \in \partial (f \circ \sigma)(X)$, and that $f$ is parabolically epi-differentiable at $\sigma(X)$ and parabolically regular at $\sigma(X)$ for $\sigma(Y)$. Then $f \circ \sigma$ is parabolically regular at $X$ for $Y$, and for any $H \in  K_{f \circ \sigma}(X,Y)$ we have 
		\begin{align*}
				{\rm d^2} (f \circ \sigma)(X \mid Y)(H) &= {\rm d}^2 f(\sigma(X) \mid \sigma(Y))(\sigma^{\prime}(X;H))   \\
				&+2\sum_{i=1}^t \left\langle \Sigma(Y)_{\alpha_i\alpha_i}, P_{\alpha_i}^{T}{\cal B}(H) P_{\alpha_i}^c (\mu_iI-\varLambda_{\alpha_i})^{-1} {P_{\alpha_i}^c}^{T}{\cal B}(H)P_{\alpha_i} \right\rangle \\
				&+2 \left\langle \Sigma(Y)_{\hat{\beta} \beta}, -U_{\hat{\beta}}^{T} H V_{\alpha} \Sigma_{\alpha}^{-1}(X) U_{\alpha}^{T} H V_{\beta} \right\rangle,
		\end{align*}
	    where $(U,V) \in {\cal O}^{m,n}(X) \cap {\cal O}^{m,n}(Y)$, the columns of $P_{\alpha_i}$ form an orthonormal basis of eigenvectors of ${\cal B}(X) $ associated with $\mu_i$, $P_{\alpha_i}^c$ is the submatrix of $P$ obtained by removing all the columns of $P_{\alpha_i}$, and $\varLambda_{\alpha_i} \in \R^{(m+n-|\alpha_i |)\times (m+n-|\alpha_i |)}$ is a diagonal matrix whose diagonal elements are the eigenvalues of $P^{T} {\cal B}(X) P$ that are not equal to $\mu_i$ for all $i \in \{ 1, \ldots,t\}$.
	 \end{thm}
	
	 {\it Proof} Set $F:= f \circ \sigma$. Since $H \in  K_F(X,Y)$, it follows from Proposition \ref{Critical Cone of Orthogonally Invariant Matrix Functions} that there exists  
	 $$
	 {\bar Q}^i \in  {\cal O}^{|\alpha_i|} \Big( \Sigma(Y)_{\alpha_i \alpha_i} \Big)  \bigcap  {\cal O}^{|\alpha_i|} \Big( \frac{1}{2} \big(U_{\alpha_i}^T H V_{\alpha_i} + V_{\alpha_i}^T H^T U_{\alpha_i} \big) \Big),
	 $$
	 and 
	 $$
	  ( {\bar Q}_{\hat{\beta} \hat{\beta}}, \widehat{{\bar Q}}_{\beta \beta} ) \in {\cal O}^{|\hat{\beta}|, |\beta|} \big(\Sigma(Y)_{\hat{\beta} \beta}\big)  \bigcap  {\cal O}^{|\hat{\beta}|, |\beta|} \big(U_{\hat{\beta}}^T H V_{\beta} \big)
	 $$
 	 for any $i \in \{1, \ldots, t\}$ and $(U,V) \in {\cal O}^{m,n}(X) \cap {\cal O}^{m,n}(Y)$ such that
	 $$
	  \Sigma(Y)_{\alpha_i \alpha_i} = {\bar Q}^i \Sigma(Y)_{\alpha_i \alpha_i} \big( {\bar Q}^i \big)^T, 
	 $$
	 $$
	  \frac{1}{2} \big(U_{\alpha_i}^T H V_{\alpha_i} + V_{\alpha_i}^T H^T U_{\alpha_i}\big) = {\bar Q}^i {\rm diag}\left(\lambda \Big(\frac{1}{2} \big(U_{\alpha_i}^T H V_{\alpha_i} + V_{\alpha_i}^T H^T U_{\alpha_i} \big) \Big) \right) \big( {\bar Q}^i \big)^T,
	 $$
	 $$
	  \Sigma(Y)_{\hat{\beta} \beta} = {\bar Q}_{\hat{\beta} \hat{\beta}}    \Sigma(Y)_{\hat{\beta} \beta}   \big( \widehat{{\bar Q}}_{\beta \beta} \big)^T,
	 $$
	 and
	 $$
	  U_{\hat{\beta}}^T H V_{\beta} = {\bar Q}_{\hat{\beta} \hat{\beta}}       \Sigma \big( U_{\hat{\beta}}^T H V_{\beta} \big)        \big( \widehat{{\bar Q}}_{\beta \beta} \big)^T.
	 $$
	 Consider the $m \times n$ block diagonal matrix 
	 $$
	 {\bar A} = {\rm diag}\left(  {\bar Q}^1 {\rm diag}(\bar z_{\alpha_1} ) ( {\bar Q}^1)^T, \ldots,  {\bar Q}^t {\rm diag}(\bar z_{\alpha_t} ) ( {\bar Q}^t)^T,  {\bar Q}_{ \hat{\beta} \hat{\beta}} \left( {\rm diag}(\bar z_\beta) \right)_{\hat{\beta} \beta} \widehat{ {\bar Q}}_{\beta \beta}^T \right),
	 $$
	 where ${\bar Q}^i$ replace the matrices $Q^i$ in the definition of the matrix $A_{\alpha \alpha}$ in \eqref{A_{alpha alpha}} for $i = 1, \ldots, t$, and ${\bar Q}_{\hat{\beta} \hat{\beta}}, \widehat{{\bar Q}}_{\beta \beta}$ replace the matrices $ Q_{\hat{\beta} \hat{\beta}}, \widehat{Q}_{\beta \beta}$ in the definition of the matrix $A_{\hat{\beta} \beta}$ in \eqref{A_{hat{beta} beta}}, respectively. It is not hard to see that the same conclusion can be achieved as the one in Proposition \ref{hat W and bar z} for the above ${\bar A}$. 
	 We can deduce from Proposition \ref{hat W and bar z} that there exists $\bar z \in \R^n$ with representation \eqref{representation of z} satisfying \eqref{z satisfing min}, where $z_{\beta^i_j} \in \R^{|\beta^i_j|}_{\downarrow}$ for $i \in \{1, \ldots, t\}, j \in \{1, \ldots, N_i\}$ and    $z_{\beta^{t+1}_k} \in \R^{|\beta^{t+1}_k|}_{\downarrow}$ for $k \in \{1, \ldots, N_{t+1}, N_{t+1}+1\}$, and further there exists a matrix   $\widehat{W} \in {\mathbb M}_{m,n}$ such that $\sigma^{\prime \prime}(X;H,\widehat{W}) = \bar z$. 
	 This, together with $Y = U \Sigma(Y) V^T$ and the fact that 
	 $\left\langle \Sigma(Y)_{\alpha_i \alpha_i}, \frac{1}{2}(V_{\alpha_i}^T \widehat{W}^T U_{\alpha_i} - U_{\alpha_i}^T \widehat{W} V_{\alpha_i}) \right\rangle = 0$ for any $i=1, \ldots, t$, we obtain that
	 \begin{align*}
		& {\rm d}^2 f(\sigma(X))(\sigma^{\prime}(X;H) \mid \sigma^{\prime \prime}(X;H,\widehat{W})) - \langle Y, \widehat{W} \rangle  \\
	   =& {\rm d}^2 f(\sigma(X))(\sigma^{\prime}(X;H) \mid \bar z ) - \sum_{i=1}^t \left\langle \Sigma(Y)_{\alpha_i\alpha_i}, U_{\alpha_i}^T \widehat{W} V_{\alpha_i} \right\rangle - \left\langle \Sigma(Y)_{\hat{\beta} \beta}, U_{\hat{\beta}}^T \widehat{W} V_{\beta} \right\rangle \\
	   =& {\rm d}^2 f(\sigma(X) \mid \sigma(Y))(\sigma^{\prime}(X;Y)) + \langle \sigma(Y), \bar z \rangle - \sum_{i=1}^t \left\langle \Sigma(Y)_{\alpha_i\alpha_i},P_{\alpha_i}^T {\cal B}(\widehat{W}) P_{\alpha_i} \right\rangle - \left\langle \Sigma(Y)_{\hat{\beta} \beta}, U_{\hat{\beta}}^T \widehat{W} V_{\beta} \right\rangle \\
	   =& {\rm d}^2 f(\sigma(X) \mid \sigma(Y))(\sigma^{\prime}(X;Y)) + \langle \sigma(Y), \bar z \rangle \\
	    & - \sum_{i=1}^t \left\langle \Sigma(Y)_{\alpha_i\alpha_i}, P_{\alpha_i}^T \Big[ 2{\cal B}(H) P_{\alpha_i}^c (\varLambda_s- \mu_i I)^{-1} {P_{\alpha_i}^c}^{T} {\cal B}(H) +  P_{\alpha} {\bar A}_{\alpha \alpha} P_{\alpha} ^T \Big]  P_{\alpha_i} \right\rangle \\
	    & - \left\langle \Sigma(Y)_{\hat{\beta} \beta}, U_{\hat{\beta}}^{T} \Big[ 2 H V_{\alpha} \Sigma_{\alpha}^{-1}(X) U_{\alpha}^{T} H + U_{\hat{\beta}} {\bar A}_{\hat{\beta} \beta} V_{\beta}^T \Big] V_{\beta} \right\rangle \\
	   =& {\rm d}^2 f(\sigma(X) \mid \sigma(Y))(\sigma^{\prime}(X;H)) + \langle \sigma(Y), \bar z \rangle  \\
	    &+2\sum_{i=1}^t \left\langle \Sigma(Y)_{\alpha_i\alpha_i}, P_{\alpha_i}^{T}{\cal B}(H) P_{\alpha_i}^c (\mu_iI-\varLambda_{\alpha_i})^{-1} {P_{\alpha_i}^c}^{T}{\cal B}(H)P_{\alpha_i} \right\rangle 
	     +2 \left\langle \Sigma(Y)_{\hat{\beta} \beta}, -U_{\hat{\beta}}^{T} H V_{\alpha} \Sigma_{\alpha}^{-1}(X) U_{\alpha}^{T} H V_{\beta} \right\rangle \\
	    & - \sum_{i=1}^t \left\langle \Sigma(Y)_{\alpha_i\alpha_i}, P_{\alpha_i}^T \Big[ P_{\alpha} {\bar A}_{\alpha \alpha} P_{\alpha} ^T \Big]  P_{\alpha_i} \right\rangle 
	      - \left\langle \Sigma(Y)_{\hat{\beta} \beta}, U_{\hat{\beta}}^{T} \Big[ U_{\hat{\beta}} {\bar A}_{\hat{\beta} \beta} V_{\beta}^T \Big] V_{\beta} \right\rangle.
	 \end{align*}
     From the definition of $\bar A$ and the representation of $\bar z$, we deduce that
     \begin{align*}
          \langle \sigma(Y), \bar z \rangle 
      & = \sum_{i=1}^t \left\langle \Sigma(Y)_{\alpha_i\alpha_i}, {\rm diag} \big(\bar z_{\alpha_i} \big) \right\rangle + 
          \left\langle  \Sigma(Y)_{\hat{\beta} \beta}, \Big( {\rm diag} \big(\bar z_{\beta} \big) \Big)_{\hat{\beta} \beta} \right\rangle  \\
      & = \sum_{i=1}^t \left\langle {\bar Q}^i \Sigma(Y)_{\alpha_i\alpha_i} ({\bar Q}^i)^T, {\bar Q}^i {\rm diag} \big(\bar z_{\alpha_i} \big)({\bar Q}^i)^T \right\rangle + 
          \left\langle {\bar Q}_{ \hat{\beta} \hat{\beta}} \Sigma(Y)_{\hat{\beta} \beta} \widehat{ {\bar Q}}_{\beta \beta}^T, {\bar Q}_{ \hat{\beta} \hat{\beta}} \Big( {\rm diag} \big(\bar z_{\beta} \big) \Big)_{\hat{\beta} \beta} \widehat{ {\bar Q}}_{\beta \beta}^T \right\rangle  \\
      & = \sum_{i=1}^t \left\langle \Sigma(Y)_{\alpha_i\alpha_i}, P_{\alpha_i}^T \Big[ P_{\alpha} {\bar A}_{\alpha \alpha} P_{\alpha} ^T \Big]      P_{\alpha_i} \right\rangle  + \left\langle \Sigma(Y)_{\hat{\beta} \beta}, U_{\hat{\beta}}^{T} \Big[ U_{\hat{\beta}} {\bar A}_{\hat{\beta} \beta} V_{\beta}^T \Big] V_{\beta} \right\rangle.
     \end{align*}
	 Thus, from \eqref{a lower estimate for second subderivative}, Proposition \ref{prop rs of ss and ps} and \eqref{formula of PS} we get
	 \begin{align*}
	        & \, {\rm d}^2 f(\sigma(X) \mid \sigma(Y))(\sigma^{\prime}(X;H)) 
	          +2\sum_{i=1}^t \left\langle \Sigma(Y)_{\alpha_i\alpha_i}, P_{\alpha_i}^{T}{\cal B}(H) P_{\alpha_i}^c (\mu_iI-\varLambda_{\alpha_i})^{-1} {P_{\alpha_i}^c}^{T}{\cal B}(H)P_{\alpha_i} \right\rangle \\
	        & \,  +2 \left\langle \Sigma(Y)_{\hat{\beta} \beta}, -U_{\hat{\beta}}^{T} H V_{\alpha} \Sigma_{\alpha}^{-1}(X) U_{\alpha}^{T} H V_{\beta} \right\rangle \\
	  \leq  & \, {\rm d^2} F(X \mid Y)(H) \\
	  \leq  & \,  \inf_{W \in {\mathbb M}_{m,n}} \{ {\rm d}^2 F(X)(H \mid W) - \langle Y, W \rangle \} \\
	    =   & \,  \inf_{W \in {\mathbb M}_{m,n}} \{ {\rm d}^2 f(\sigma(X))(\sigma^{\prime}(X;H) \mid \sigma^{\prime \prime}(X;H,W)) - \langle Y, W \rangle\} \\
	  \leq  & \,  {\rm d}^2 f(\sigma(X))(\sigma^{\prime}(X;H) \mid \sigma^{\prime \prime}(X;H,\widehat{W})) - \langle Y, \widehat{W} \rangle  \\
	    =   & \,  {\rm d}^2 f(\sigma(X) \mid \sigma(Y))(\sigma^{\prime}(X;H)) 
	          +2\sum_{i=1}^t \left\langle \Sigma(Y)_{\alpha_i\alpha_i}, P_{\alpha_i}^{T}{\cal B}(H) P_{\alpha_i}^c (\mu_iI-\varLambda_{\alpha_i})^{-1} {P_{\alpha_i}^c}^{T}{\cal B}(H)P_{\alpha_i} \right\rangle \\
	        & \, +2 \left\langle \Sigma(Y)_{\hat{\beta} \beta}, -U_{\hat{\beta}}^{T} H V_{\alpha} \Sigma_{\alpha}^{-1}(X) U_{\alpha}^{T} H V_{\beta} \right\rangle.
	 \end{align*}
     It follows immediately that ${\rm d^2} F(X \mid Y)(H) = \inf\limits_{W \in {\mathbb M}_{m,n}} \{ {\rm d}^2 F(X)(H \mid W) - \langle Y, W \rangle \} $ for any $H \in  K_F(X,Y)$, which implies that $F$ is parabolically regular at $X$ for $Y$.
     Moreover, this verifies the proposed formula for the second subderivative of $F$ at $X$ for $Y$ for any $H \in  K_F(X,Y)$, thereby completing the proof. 
     \qed

     Combining Theorem \ref{second subderivative of OIMF} with \cite[Theorem 3.8]{MS}, we provide sufficient conditions for twice epi-differentiability of orthogonally invariant matrix function.
     
     \begin{prop}
     	Let $f \colon \R^n \rightarrow [-\infty,+\infty]$ be an absolutely symmetric function, and be lsc, convex, and locally Lipschitz continuous relative to its domain. Assume that $Y \in \partial (f \circ \sigma)(X)$ and $f$ is parabolically epi-differentiable at $\sigma(X)$ and parabolically regular at $\sigma(X)$ for $\sigma(Y)$. Then $f \circ \sigma$ is twice epi-differentiable at $X$ for $Y$.
     \end{prop}

     As an application, we present second-order optimality conditions for the following optimization problem according to \cite[Theorem 13.24]{Rockafellar Wets}. Given a twice differentiable function $\psi \colon {\mathbb M}_{m,n} \rightarrow \R$ and an absolutely symmetric function $f \colon \R^n \rightarrow [-\infty,+\infty]$, consider the optimization problem
     $$
     \min_{X \in {\mathbb M}_{m,n}} \psi(X) + (f \circ \sigma)(X).  \leqno(P)
     $$
     From \cite[Exercise 13.18]{Rockafellar Wets}, we can derive the following conclusion.
     \begin{thm}
       Let $f \colon \R^n \rightarrow [-\infty,+\infty]$ be an absolutely symmetric function, and be lsc, convex, and locally Lipschitz continuous relative to its domain. Assume that $X_0$ is a feasible solution to $(P)$ with $f$ being parabolically epi-differentiable at $\sigma(X_0)$ and parabolically regular at $\sigma(X_0)$ for $\sigma(Y)$ and $Y \in \partial (f \circ \sigma)(X_0)$. If $-\nabla \psi(X_0) \in \partial (f \circ \sigma)(X_0) $, then the following second-order optimality conditions for $(P)$ hold.
    	
    	{\rm (i)} If $X_0$ is a local minimizer of $(P)$, then the second-order necessary condition
    	\begin{align*}
    		&\nabla^2 \psi(X_0)(H,H) + {\rm d}^2 f(\sigma(X) \mid \sigma(Y))(\sigma^{\prime}(X;H)) 
    		+ 2\sum_{i=1}^t \left\langle \Sigma(Y)_{\alpha_i\alpha_i}, P_{\alpha_i}^{T}{\cal B}(H) P_{\alpha_i}^c (\mu_iI-\varLambda_{\alpha_i})^{-1} {P_{\alpha_i}^c}^{T}{\cal B}(H)P_{\alpha_i} \right\rangle \\
    		&+ 2 \left\langle \Sigma(Y)_{\hat{\beta} \beta}, -U_{\hat{\beta}}^{T} H V_{\alpha} \Sigma_{\alpha}^{-1}(X) U_{\alpha}^{T} H V_{\beta} \right\rangle
    		\geq 0
    	\end{align*}
    	holds for all $H \in  K_{f \circ \sigma}\big(X_0,-\nabla \psi(X_0) \big)$.
    	
    	{\rm (ii)}  Having the second-order sufficient condition
    	\begin{align*}
    		\begin{split}
    			&\nabla^2 \psi(X_0)(H,H) +{\rm d}^2 f(\sigma(X) \mid \sigma(Y))(\sigma^{\prime}(X;H)) 
    			+ 2\sum_{i=1}^t \left\langle \Sigma(Y)_{\alpha_i\alpha_i}, P_{\alpha_i}^{T}{\cal B}(H) P_{\alpha_i}^c (\mu_iI-\varLambda_{\alpha_i})^{-1} {P_{\alpha_i}^c}^{T}{\cal B}(H)P_{\alpha_i} \right\rangle \\
    			&+ 2 \left\langle \Sigma(Y)_{\hat{\beta} \beta}, -U_{\hat{\beta}}^{T} H V_{\alpha} \Sigma_{\alpha}^{-1}(X) U_{\alpha}^{T} H V_{\beta} \right\rangle
    			> 0
    		\end{split}
    	\end{align*}
    	holds for all $H \in  K_{f \circ \sigma}\big(X_0,-\nabla \psi(X_0) \big)$ is equivalent to having the existence $\epsilon > 0$ and $c > 0$ such that 
    	$$
    	\psi(X) + (f \circ \sigma)(X) \geq \psi(X_0) + (f \circ \sigma)(X_0) + c \|X-X_0\|^2  \,\,\mbox{when}\,\, X \in {\mathbb B}_{\epsilon}(X_0).
    	$$
     \end{thm}

     In particular, when absolutely symmetric function $f$ is a polyhedral function, whose epigraph is a polyhedral convex set, all the assumptions in Theorem \ref{second subderivative of OIMF} are satisfied. In this case, the second subderivative of $f \circ \sigma$ has a simple representation as follow.
     
     \begin{cor}\label{cor in polyhedral case}
    	Let an absolutely symmetric function $f \colon \R^n \rightarrow [-\infty,+\infty]$ be polyhedral. Assume that $\mu_1> \cdots > \mu_t$ are the distinct nonzero singular values of $X \in {\mathbb M}_{m,n}$ and that $Y \in \partial (f \circ \sigma)(X)$. Then for any $H \in {\mathbb M}_{m,n}$, we have 
    	\begin{align}\label{second subderivative in polyhedral case}
    		\begin{split}
    			{\rm d^2} (f \circ \sigma) (X \mid Y)(H) &= \delta_{K_{f \circ \sigma}(X,Y)}(H)  
    			+2\sum_{i=1}^t \left\langle \Sigma(Y)_{\alpha_i\alpha_i}, P_{\alpha_i}^{T}{\cal B}(H) P_{\alpha_i}^c (\mu_iI-\varLambda_{\alpha_i})^{-1} {P_{\alpha_i}^c}^{T}{\cal B}(H)P_{\alpha_i} \right\rangle \\
    			&+2 \left\langle \Sigma(Y)_{\hat{\beta} \beta}, -U_{\hat{\beta}}^{T} H V_{\alpha} \Sigma_{\alpha}^{-1}(X) U_{\alpha}^{T} H V_{\beta} \right\rangle,
    		\end{split}
    	\end{align}
       	where $(U,V) \in {\cal O}^{m,n}(X) \cap {\cal O}^{m,n}(Y)$, the columns of $P_{\alpha_i}$ form an orthonormal basis of eigenvectors of ${\cal B}(X) $ associated with $\mu_i$, $P_{\alpha_i}^c$ is the submatrix of $P$ obtained by removing all the columns of $P_{\alpha_i}$, and $\varLambda_{\alpha_i} \in \R^{(m+n-|\alpha_i |)\times (m+n-|\alpha_i |)}$ is a diagonal matrix whose diagonal elements are the eigenvalues of $P^{T} {\cal B}(X) P$ that are not equal to $\mu_i$ for all $i \in \{ 1, \ldots,t\}$.
     \end{cor}
 
     {\it Proof}  It follows from \cite[Theorem 2.49]{Rockafellar Wets} that $f$ is a convex piecewise linear function. By the definition of ‘piecewise linear’ in \cite[Definition 2.47]{Rockafellar Wets}, the convex set $D= {\rm dom}f$ is the union of finitely many polyhedral sets $D_i$, relative to each of which $f(x)$ is given by an expression of the form $\langle a_i, x \rangle + \alpha_i$ for some vector $a_i \in \mathbb{R}^n$ and $\alpha_i \in \mathbb{R}$. Fixing $\bar x$, we can take $f(\bar x) = 0$ and suppose that $\bar x$ belongs to every $D_i$ (since this can be arranged by adding to $f$ the indicator of a polyhedral neighborhood of $\bar x$, which won’t affect the local properties of $f$ at $\bar x$). Then
     $$
     f_i(\bar x + \tau w) - f_i(\bar x) = \langle a_i, \bar x + \tau w \rangle + \alpha_i-(\langle a_i, \bar x \rangle + \alpha_i)= \tau \langle a_i,w\rangle \,\, \mbox{ when } \,\,\bar x + \tau w \in D_i,
     $$ 
     and $T_D(\bar x) = \bigcup\limits_i T_{D_i}(\bar x)$, so
     \begin{equation*}
     	\bigtriangleup_{\tau} f(\bar x)(w) = \frac{f(\bar x+\tau w)-f(\bar x)}{\tau}  
     	=  \left\{
     	\begin{array}{ll}
     		\langle a_i, w \rangle   & \mbox{ when } w \in \frac{D_i-\bar x}{\tau},\\
     		+\infty & \mbox{ when } w \notin \frac{D-\bar x}{\tau}, \\
     	\end{array}
     	\right.
     \end{equation*}
     \begin{equation*}
     	{\rm d}f(\bar x)(w)= \left\{
     	\begin{array}{ll}
     		\langle a_i, w \rangle   & \mbox{if } w \in T_{D_i}(\bar x),\\
     		+\infty & \mbox{if } w \notin T_{D}(\bar x).\\
     	\end{array}
     	\right.
     \end{equation*}
     Consider any $\bar v \in \partial f(\bar x)$. We have $\langle \bar v,w \rangle \leq {\rm d}f(\bar x)(w)$ for all $w$, and
     \begin{equation*}
     	\bigtriangleup_{\tau}^2 f(\bar x \mid \bar v)(w) = \frac{f(\bar x+\tau w)-f(\bar x)-\tau \langle \bar v,w\rangle}{\frac{1}{2} \tau^2}  
     	=  \left\{
     	\begin{array}{ll}
     		\frac{\langle a_i- \bar v, w \rangle}{\frac{1}{2} \tau}   & \mbox{ when } w \in \frac{D_i-\bar x}{\tau},\\
     		+\infty & \mbox{ when } w \notin \frac{D-\bar x}{\tau}, \\
     	\end{array}
     	\right.
     \end{equation*}
     where $\langle a_i- \bar v, w \rangle = {\rm d}f(\bar x)(w)- \langle \bar v, w \rangle \geq 0$ when $w \in T_{D_i}(\bar x)$. Therefore, we have
     \begin{equation*}
     	{\rm d}^2 f(\bar x \mid \bar v)(w)= \left\{
     	\begin{array}{ll}
     		0   & \mbox{if } w \in T_{D_i}(\bar x) \cap K_f(\bar x, \bar v),\\
     		+\infty & \mbox{if } w \notin T_{D}(\bar x) \cap K_f(\bar x, \bar v), \\
     	\end{array}
     	\right.
     \end{equation*}
     where $K_f(\bar x, \bar v)= \{w \in \R^n \colon {\rm d}f(\bar x)(w) = \langle \bar v,w \rangle\}$ is the critical cone of a function $f$ at a point $\bar x$ for $\bar v$. When $w \notin T_D(\bar x)$, we have ${\rm d}f(\bar x)(w) = \infty$ and $w \notin K_f(\bar x, \bar v)$, so $ K_f(\bar x, \bar v) \subset T_D(\bar x)$. Then 
     $$
     \bigcup\limits_i [T_{D_i}(\bar x) \cap K_f(\bar x, \bar v)] = T_{D}(\bar x) \cap K_f(\bar x, \bar v) =K_f(\bar x, \bar v),
     $$
     and further
     \begin{equation}\label{rem d2f}
     	{\rm d}^2 f(\bar x \mid \bar v)(w) = \delta_{K_f(\bar x, \bar v)}(w).
     \end{equation}
     Since $f$ is a convex piecewise linear function, we obtain from \cite[Exercise 13.61]{Rockafellar Wets} and \cite[Example 3.2]{MS} that $f$ is parabolically epi-differentiable and parabolically regular at $\sigma(X)$. Moreover, it follows from Theorem \ref{PS of OIMF} and \cite[Proposition 3.4]{MS} that 
     \begin{equation}\label{dom d2f eq kf}
     	{\rm dom}\,{\rm d}^2 (f \circ \sigma) (X \mid Y) = K_{f \circ \sigma}(X,Y)
     \end{equation}
     Take $H \in K_{f \circ \sigma}(X,Y)$ and observe from Proposition \ref{Critical Cone of Orthogonally Invariant Matrix Functions} that $\sigma^{\prime}(X;H) \in  K_f(\sigma(X),\sigma(Y))$.
     Thus, we obtain from \eqref{rem d2f} that 
     $$
     {\rm d}^2 f(\sigma(X) \mid \sigma(Y))(\sigma^{\prime}(X;H)) = \delta_{K_f(\sigma(X),\sigma(Y))}(\sigma^{\prime}(X;H)) = 0.
     $$
     It follows from Theorem \ref{second subderivative of OIMF} that
     \begin{align*}
     	{\rm d^2} (f \circ \sigma)(X \mid Y)(H) &= 2\sum_{i=1}^t \left\langle \Sigma(Y)_{\alpha_i\alpha_i}, P_{\alpha_i}^{T}{\cal B}(H) P_{\alpha_i}^c (\mu_iI-\varLambda_{\alpha_i})^{-1} {P_{\alpha_i}^c}^{T}{\cal B}(H)P_{\alpha_i} \right\rangle \\
     	&+2 \left\langle \Sigma(Y)_{\hat{\beta} \beta}, -U_{\hat{\beta}}^{T} H V_{\alpha} \Sigma_{\alpha}^{-1}(X) U_{\alpha}^{T} H V_{\beta} \right\rangle.
     \end{align*}
     On the other hand, if $H \notin K_{f \circ \sigma}(X,Y)$, we have ${\rm d^2} (f \circ \sigma)(X \mid Y)(H) = \infty$ by \eqref{dom d2f eq kf}. This proves the claimed formula for the second subderivative of $f \circ \sigma$ at $X$ for $Y$. 
     \qed

	 \begin{rem}
	  If the absolutely symmetric function $f$ is $\ell_1$-norm, then $f$ is a convex piecewise linear function and orthogonally invariant matrix function $f \circ \sigma$ is nuclear norm. Thus, $f$ has the following expression of the form
	  \begin{equation}\label{the expression of l1 norm}
	   	\begin{split}
	   		f(x) = \|x\|_1 = \left\{
	   		\begin{array}{cl}
	   			x_1 +x_2 + \sum \limits_{i=3}^n x_i  =:\langle a_1, x \rangle        & \mbox{if } x \in \mathbb{R}_{+} \times \mathbb{R}_{+} \times \mathbb{R}^{n-2}_{+} =: D_1,  \\
	   			x_1 -x_2 + \sum \limits_{i=3}^n x_i  =:\langle a_2, x \rangle        & \mbox{if } x \in \mathbb{R}_{+} \times \mathbb{R}_{-} \times \mathbb{R}^{n-2}_{+} =: D_2,  \\
	   			-x_1 +x_2 + \sum \limits_{i=3}^n x_i  =:\langle a_3, x \rangle        & \mbox{if } x \in \mathbb{R}_{-} \times \mathbb{R}_{+} \times \mathbb{R}^{n-2}_{+} =: D_3,   \\
	   			-x_1 -x_2 + \sum \limits_{i=3}^n x_i  =:\langle a_4, x \rangle        & \mbox{if } x \in \mathbb{R}_{-} \times \mathbb{R}_{-} \times \mathbb{R}^{n-2}_{+} =: D_4,   \\
	   			\vdots                                                               &     \ \ \  \ \ \  \ \ \   \vdots               \\
	   			-\sum \limits_{i=1}^{n-1} x_i +x_n    =:\langle a_{2^n-1}, x \rangle & \mbox{if } x \in \mathbb{R}^{n-1}_{-} \times \mathbb{R}_{+} =: D_{2^n-1},   \\
	   			-\sum \limits_{i=1}^{n} x_i     =:\langle a_{2^n}, x \rangle         & \mbox{if } x \in \mathbb{R}^{n}_{-} =: D_{2^n},   \\
	   		\end{array}
	   		\right.
	   	\end{split}
	  \end{equation}
      where $D:= {\rm dom}f = \bigcup \limits_{i=1}^{2^n} D_i$, $D_i$ is a polyhedral set, and $a_i \in \mathbb{R}^n$  for any $i=1, \ldots,2^n$. We can express $D_i$ as $D_i= K^i_1 \times K^i_2 \times \cdots \times K^i_n$, where $K^i_j$ is $\mathbb{R}_{+}$ or $\mathbb{R}_{-}$ for $i \in \{1, \cdots, 2^n\} , \, j\in \{1, \cdots, n\}$. 
      For any $\bar x\in \mathbb{R}^n$, we denote $I^{+}(\bar x):=\{i \colon {\bar x}_i >0 \},$  $I^{-}(\bar x):=\{i \colon {\bar x}_i <0 \}$, and $I^0(\bar x) := \{i \colon {\bar x}_i =0 \}$. For any $\bar v \in \partial f(\bar x)$, we obtain that
      \begin{equation*}
       	\bar v_i = \left\{
       	\begin{array}{cl}
       		1               & \mbox{if } i \in I^{+}(\bar x),  \\
       		-1              & \mbox{if } i \in I^{-}(\bar x),  \\
       		{ [-1,1] }      & \mbox{if } i \in I^{0}(\bar x).  \\
       	\end{array}
       	\right.
      \end{equation*}
      For any $ w \in \mathbb{R}^n$, we have
      \begin{eqnarray*}
       	{\rm d}f(\bar x)(w) &=& \lim_{t \downarrow 0} \frac{\|\bar x + t w \|_1 -\|\bar x\|_1}{t} \\
       	&=& \lim_{t\downarrow 0} \frac{ \sum \limits_{i \in I^{+}(\bar x) \cup I^{-}(\bar x) } |(\bar x + t w)_i | + \sum \limits_{i \in I^0(\bar x)} t| w_i | - \sum \limits_{i \in I^{+}(\bar x) \cup I^{-}(\bar x)} |{\bar x}_i| }{t}\\
       	&=&  \lim_{t\downarrow 0} \frac{ \sum \limits_{i \in I^{+}(\bar x) \cup I^{-}(\bar x)}[|(\bar x + t w)_i|  - |{\bar x}_i|] + \sum \limits_{i \in I^0(\bar x)} t|w_i|} {t}\\
       	&=& \sum\limits_{i \in I^{+}(\bar x)}w_i  - \sum\limits_{i \in I^{-}(\bar x)}w_i  + \sum\limits_{i \in I^0(\bar x)} |w_i|.
      \end{eqnarray*} 	
	  Then the critical cone of $f$ at a point $\bar x$ for $\bar v$ can be calculated as
	  \begin{equation}\label{Kf in eg}
	  	 K_f(\bar x, \bar v) = \left\{ w \in \mathbb{R}^n \colon
	  	\begin{array}{ll}
	   	  w_i \in \mathbb{R}       & \mbox{when } i \in I^{+}(\bar x) \cup I^{-}(\bar x),  \\
	  	  w_i \in \mathbb{R}_{+}   & \mbox{when } i \in I^{0}(\bar x), \bar v_i=1,         \\
	  	  w_i \in \mathbb{R}_{-}   & \mbox{when } i \in I^{0}(\bar x), \bar v_i=-1,        \\
	  	  w_i = 0                  & \mbox{when } i \in I^{0}(\bar x), |\bar v_i| < 1      \\
	  	\end{array}
	   \right\},
	  \end{equation}
	  and the tangent cone of $K^i_j$ at ${\bar x}_i$ can be expressed as
	  \begin{equation*}
	  	T_{K^i_j} ({\bar x}_i) = \left\{
	  	\begin{array}{cl}
	  		\mathbb{R}              & \mbox{if } i \in I^{+}(\bar x) \cup I^{-}(\bar x),  \\
	  		   K^i_j                & \mbox{if } i \in I^{0}(\bar x).  \\
	  	\end{array}
	  	\right.
	  \end{equation*}
	  This, together with \eqref{the expression of l1 norm}, justifies ${\rm d}f(\bar x)(w)= \langle a_i, w \rangle$ for any $w \in T_{D_i}(\bar x) = T_{K^i_1}({\bar x}_1) \times T_{K^i_2}({\bar x}_2) \times \cdots \times T_{K^i_n}({\bar x}_n)$ and ${\rm d}f(\bar x)(w)= +\infty $ for any $w \notin T_{D}(\bar x) = \mathbb{R}^n$.

	  Combining \eqref{subderivative of Psi} with 
	  \begin{equation*}
	  	{\rm d} (f \circ \sigma)(X)(H) = {\rm tr}((U^{(1)})^T   H V^{(1)}) + \|(U^{(2)})^T  HV^{(2)}\|_*,
	  \end{equation*} 
	  we obtain that $K_{f \circ \sigma}(X,Y) = K_{\Psi_n}(X,Y)$, where $n$ satisfies the condition of Theorem \ref{twice epi-differentiable of Psi_t}. Take $H \in K_{f \circ \sigma}(X,Y)$ and observe from Proposition \ref{Critical Cone of Orthogonally Invariant Matrix Functions} that $\sigma^{\prime}(X;H) \in  K_f(\sigma(X),\sigma(Y))$. It follows from \eqref{rem d2f}, \eqref{dom d2f eq kf} and Theorem \ref{second subderivative of OIMF} that for any $H \in {\mathbb M}_{m,n}$, 
	  \begin{align*}
	 	{\rm d^2} (f \circ \sigma) (X \mid Y)(H) &= \delta_{ K_{\Psi_n}(X,Y)}(H)  
	 	+2\sum_{i=1}^t \left\langle \Sigma(Y)_{\alpha_i\alpha_i}, P_{\alpha_i}^{T}{\cal B}(H) P_{\alpha_i}^c (\mu_iI-\varLambda_{\alpha_i})^{-1} {P_{\alpha_i}^c}^{T}{\cal B}(H)P_{\alpha_i} \right\rangle \\
	 	&+2 \left\langle \Sigma(Y)_{\hat{\beta} \beta}, -U_{\hat{\beta}}^{T} H V_{\alpha} \Sigma_{\alpha}^{-1}(X) U_{\alpha}^{T} H V_{\beta} \right\rangle,
	  \end{align*}
	  which is the second subderivative of nuclear norm as that in Corollary \ref{twice epi-differentiable of nuclear norm} using a different approach.	    
     \end{rem}
	
     In what follows, we present second-order optimality conditions for the optimization problem $(P)$ with $f$ being polyhedral.
	
	 \begin{cor}
	  Assume that $X_0$ is a feasible solution to $(P)$ with $f$ being polyhedral. If $-\nabla \psi(X_0) \in \partial (f \circ \sigma)(X_0) $, then the following second-order optimality conditions for $(P)$ hold.
		
	  {\rm (i)} If $X_0$ is a local minimizer of $(P)$, then the second-order necessary condition
		\begin{align*}
		  &\nabla^2 \psi(X_0)(H,H) 
		  + 2\sum_{i=1}^t \left\langle \Sigma(Y)_{\alpha_i\alpha_i}, P_{\alpha_i}^{T}{\cal B}(H) P_{\alpha_i}^c (\mu_iI-\varLambda_{\alpha_i})^{-1} {P_{\alpha_i}^c}^{T}{\cal B}(H)P_{\alpha_i} \right\rangle \\
		  &+ 2 \left\langle \Sigma(Y)_{\hat{\beta} \beta}, -U_{\hat{\beta}}^{T} H V_{\alpha} \Sigma_{\alpha}^{-1}(X) U_{\alpha}^{T} H V_{\beta} \right\rangle
		  \geq 0
		\end{align*}
		holds for all $H \in  K_{f \circ \sigma}\big(X_0,-\nabla \psi(X_0) \big)$.
		
	  {\rm (ii)}  Having the second-order sufficient condition
		\begin{align}\label{second-order condition in application}
		  \begin{split}
		  	&\nabla^2 \psi(X_0)(H,H) 
		  	+ 2\sum_{i=1}^t \left\langle \Sigma(Y)_{\alpha_i\alpha_i}, P_{\alpha_i}^{T}{\cal B}(H) P_{\alpha_i}^c (\mu_iI-\varLambda_{\alpha_i})^{-1} {P_{\alpha_i}^c}^{T}{\cal B}(H)P_{\alpha_i} \right\rangle \\
		  	&+ 2 \left\langle \Sigma(Y)_{\hat{\beta} \beta}, -U_{\hat{\beta}}^{T} H V_{\alpha} \Sigma_{\alpha}^{-1}(X) U_{\alpha}^{T} H V_{\beta} \right\rangle
		  	> 0
		  \end{split}
		\end{align}
		holds for all $H \in  K_{f \circ \sigma}\big(X_0,-\nabla \psi(X_0) \big)$ is equivalent to having the existence $\epsilon > 0$ and $c > 0$ such that 
		$$
		  \psi(X) + (f \circ \sigma)(X) \geq \psi(X_0) + (f \circ \sigma)(X_0) + c \|X-X_0\|^2  \,\,\mbox{when}\,\, X \in {\mathbb B}_{\epsilon}(X_0).
		$$
	 \end{cor}
 
     In the case of $f$ being polyhedral, the second-order sufficient condition \eqref{second-order condition in application} is consistent with the ``no-gap" second order sufficient condition in \cite[Theorem 3.5]{CDZ}, and we have provided the specific expression for the second term in the latter.

	\section{Conclusion}
	 The primary objective of this work was to explore several second-order properties of orthogonally invariant matrix functions, with a focus on parabolic epi-differentiability, parabolic regularity, and twice epi-differentiability. 
	 Our results are based on the concept of metric subregularity constraint qualification, which, in this context, is automatically satisfied. For a convex orthogonally invariant matrix function, we derive the exact formula for its second subderivatives and establish sufficient conditions for twice epi-differentiability. 
	 Furthermore, we present second-order optimality conditions for a class of matrix optimization problems. 
	 In particular, inspired by Torki \cite{Torki 2001}, we adopt a different approach to calculate the first- and second-order epi-derivatives of the nuclear norm.

     \section*{Acknowledgements}
     The research of the second author was supported by the National Natural Sciences Grant of China (No.11701126). The research of the third author was supported by the National Natural Sciences Grant of China (No. 11871182).

\end{document}